# New Archimedean Bivariate Copulas Modeling Positive Dependency

Iman M. Attia *

Imanattiathesis1972@gmail.com ,imanattia1972@gmail.com

**Department of Mathematical Statistics, Faculty of Graduate Studies for Statistical Research, Cairo University, Egypt*

## *Abstract*

In this paper, the author introduces new methods to construct Archimedean copulas. The generator of each copula fulfills the sufficient conditions as regards the boundary and being continuous, decreasing, and convex. Each inverse generator also fulfills the necessary conditions as regards the boundary conditions, marginal uniformity, and 2-increasing properties. Although these copulas satisfy these conditions, they have some limitations. They do not cover the entire dependency spectrum, ranging from perfect negative dependency to perfect positive dependency, passing through the state of independence. Both copulas exhibit positive dependency and upper tail dependency, but neither has lower tail dependency. The product copula is present for each of them. For each copula, the author discusses the derivation, the properties, whether it has a singular part or not, the Kendall tau measure for dependency, and the upper and lower tail dependency. The article shows figures for depicting the generator function, joint CDF, and joint PDF for each copula. For each copula, inference is supported by real data analysis using the inference function for margins (IFM) for estimation.

### Introduction

A copula can express the nonlinear dependency between two or more variables, in comparison with correlation which can capture linear dependency between variables. It is a technique used to implement the joint distribution of the many variables (two or more). This has many wide applications in modern era of science like hydrology and geophysics (Valle & Kaplan, 2019) (Liu et al., 2018) (Salvadori et al., 2007), transportation research(Ma et al., 2017) (Bhat & Eluru, 2009) , medicine (Kuss et al., 2014) (Lapuyade-Lahorgue et al., 2017) (Gomes et al., 2019) , engineering (Kilgore & Thompson, 2011), biology (Konigorski et al., 2014) (Dokuzoğlu & Purutçuoğlu, 2017), climate research (Schölzel & Friederichs, 2008) (Oppenheimer et al., 2016), economics (Oh & Patton, 2018) (Kole et al., 2007) (De Lira Salvatierra & Patton, 2015). Sklar (Sklar, 1973) showed that a copula exists for any multivariate distribution, such that the joint distribution equals the copula applied to the marginal.

**Definition** 1 (**Sklar's theorem 1**): For any bivariate distribution function F with marginal $F_1$ and $F_2$ , there exists a copula such that

$$F(x,y) = C(F_1(x), F_2(y)), \qquad x, y \in \overline{R} = R \cup \{-\infty, \infty\}$$

If F is continuous, then the copula C is unique. Otherwise, it is uniquely determined on the $F_1(\overline{R}) \times F_2(\overline{R})$ . The converse is also true, for any copula and univariate distribution functions $F_1$ and $F_2$ , the function $C(F_1(x), F_2(y))$ is a bivariate distribution function with margins $F_1$ and $F_2$ .

A bivariate distribution function $C(u, v)$ with marginal $u = F_X(x)$ and $v = F_Y(y)$ is said to be generated by an Archimedean copula, if it can be written in the form of equation (1)

$$C(u,v) = \varphi^{-1}\big(\varphi(u) + \varphi(v)\big) \qquad (1)$$

The generator $\varphi(u)$ should satisfy the following conditions (sufficient conditions):

1) $\varphi(u)$ is continuous, strictly decreasing and convex function mapping $[0,1]$ onto $[0,\infty)$ , in other words , $\varphi'(u) < 0$ ( decreasing) and $\varphi''(u) > 0$ ( convex)

2) $\varphi(0) = \infty$

3) $\varphi(1) = 0$

**Definition 2**: a bivariate copula is a non-decreasing and a right-continuous bivariate function mapping

$C: [0,1] \times [0,1] \rightarrow [0,1]$ which satisfies the following conditions

1) $C(0, v) = C(u, 0) = 0,\ \forall (u, v) \in [0,1]^2$ ( Grounded)

2) $C(1, v) = v$ and $C(u, 1) = u\ \forall (u, v) \in [0,1]^2$ ( uniform marginal)

3) $C(u_1, v_1) - C(u_1, v_2) - C(u_2, v_1) + C(u_2, v_2) \geq 0, for\ u_1 > u_2\ and\ v_1 > v_2$

$\forall\ (u_1, u_2\ , v_1, v_2) \in [0,1]^4$ (2- increasing) . This last property is equal to $\frac{\partial C(u,v)}{\partial u \partial v} \geq 0$

$\frac{\partial C(u,v)}{\partial u \partial v} = c(u, v)$ is the joint PDF of the copula.

So the copula links the marginal distributions through the dependency parameter.

The singular part of the copula is defined by (Genest & Mackay, 1986) in the following theorem.

**Theorem 2**: the ratio between the generator and the first derivative of the generator evaluated at zero does not equal zero. The theorem states that the distribution $C(u, v)$ generated by $\varphi$ has a singular component iff $\frac{\varphi(0)}{\varphi'(0)} \neq 0$ in that case, $\varphi(u) + \varphi(v) = \varphi(0)$ with probability = $-\frac{\varphi(0)}{\varphi'(0)}$

This paper is structured in 5 sections. Section 1 discusses the first copula, its derivation, properties and figures for its PDF and CDF. Section 2 explains the second copula with the same items. Section 3 illustrates the method of estimations. Section 4 discusses real data analysis for copula fitting. Section 5 comprehends the conclusion. Section 6 mentions the future work

## Section 1: (First Copula)

Let's say the inverse generator is $C\big(\varphi(u) + \varphi(v)\big) = \Phi^{-1}(t) = z = exp\left(\frac{-t^{\alpha}}{\alpha}\right)$ . The generator is $\ln \Phi^{-1}(t) = \ln\left\{exp\left(\frac{-t^{\alpha}}{\alpha}\right)\right\} = \frac{-t^{\alpha}}{\alpha}$ and hence, $-\alpha \ln \Phi^{-1}(t) = t^{\alpha}$. The generator $\varphi(z) = \big(-\alpha \ln \Phi^{-1}(t)\big)^{\frac{1}{\alpha}} = \big(-\alpha \ln(z)\big)^{\frac{1}{\alpha}} = t$ , $z \in (0,1)$. The generator should fulfill the sufficient conditions: $\varphi(0) = (-\alpha \ln 0)^{\frac{1}{\alpha}} = \infty$ ,$\varphi(1) = (-\alpha \ln 1)^{\frac{1}{\alpha}} = 0$, and $\varphi'(z) = \frac{1}{\alpha}\,(-\alpha \ln z)^{\frac{1}{\alpha}-1}\left(\frac{-\alpha}{z}\right) = (-\alpha \ln z)^{\frac{1}{\alpha}-1}\left(\frac{-1}{z}\right) < 0$ ensures that the generator is a decreasing function. While; $\varphi''(z) = \left(\frac{1}{\alpha}-1\right)\ (-\alpha \ln z)^{\frac{1}{\alpha}-2}\left(\frac{-\alpha}{z}\right)\left(\frac{-1}{z}\right) + \ (-\alpha \ln z)^{\frac{1}{\alpha}-1}\left(\frac{1}{z^2}\right) > 0$, ensures that the generator is convex for $0 < \alpha \leq 1$.

for bivariate distribution: $\varphi(u) = (-\alpha \ln u)^{\frac{1}{\alpha}}$ and $\varphi(v) = (-\alpha \ln v)^{\frac{1}{\alpha}}$

$$\varphi(w) = \varphi(u) + \varphi(v) = (-\alpha \ln u)^{\frac{1}{\alpha}} + (-\alpha \ln v)^{\frac{1}{\alpha}}$$

$$C\big(\varphi(u) + \varphi(v)\big) = exp\left(-\frac{1}{\alpha}\left[(-\alpha \ln u)^{\frac{1}{\alpha}} + (-\alpha \ln v)^{\frac{1}{\alpha}}\right]^{\alpha}\right)$$

$$\mathrm{C}(u,v) = C\big(F_X(x), F_Y(y)\big)$$

For a copula to be valid copula it should fulfill the boundary condition, the marginal uniformity and 2-increasing condition. And this is equivalent for the inverse generator the necessary conditions to be fulfilled are the following:

**Proposition 1**: The boundary conditions are $\mathrm{C}(0,v) = \mathrm{C}(u,0) = 0$ and marginal uniformity are $\mathrm{C}(1,v) = v$ *and* $\mathrm{C}(u,1) = u$

**Proof:**

$$C(1,v) = exp\left(-\frac{1}{\alpha}\left[(-\alpha \ln(1))^{\frac{1}{\alpha}} + (-\alpha \ln v)^{\frac{1}{\alpha}}\right]^{\alpha}\right) = exp\left(-\frac{1}{\alpha}\left[0 + (-\alpha \ln v)^{\frac{1}{\alpha}}\right]^{\alpha}\right) = exp\left(-\frac{1}{\alpha}\left[(-\alpha \ln v)^{\frac{1}{\alpha}}\right]^{\alpha}\right) = v$$

$$C(u,1) = exp\left(-\frac{1}{\alpha}\left[(-\alpha \ln u)^{\frac{1}{\alpha}}\right]^{\alpha}\right) = exp\left(\frac{-\alpha \ln u}{-\alpha}\right) = u$$

**Proposition 2**: for copula to be valid it should be 2-increasing, in other words, $\frac{\partial^2 \mathrm{C}(u,v)}{\partial u \partial v} \geq 0$

**Proof:**

$$\frac{\partial \mathrm{C}(u,v)}{\partial u} = exp\left(-\frac{1}{\alpha}\left[(-\alpha \ln(u))^{\frac{1}{\alpha}} + (-\alpha \ln v)^{\frac{1}{\alpha}}\right]^{\alpha}\right)\left(\frac{-1}{\alpha}\right)(\alpha)\left[(-\alpha \ln(u))^{\frac{1}{\alpha}} + (-\alpha \ln v)^{\frac{1}{\alpha}}\right]^{\alpha-1}\left(\frac{1}{\alpha}(-\alpha \ln(u))^{\frac{1}{\alpha}-1}\left(\frac{-\alpha}{u}\right)\right)$$

$$\frac{\partial \mathrm{C}(u,v)}{\partial u} = exp\left(-\frac{1}{\alpha}\left[(-\alpha \ln(u))^{\frac{1}{\alpha}} + (-\alpha \ln v)^{\frac{1}{\alpha}}\right]^{\alpha}\right)\left[(-\alpha \ln(u))^{\frac{1}{\alpha}} + (-\alpha \ln v)^{\frac{1}{\alpha}}\right]^{\alpha-1}\left((-\alpha \ln(u))^{\frac{1}{\alpha}-1}\left(\frac{1}{u}\right)\right)$$

$$\frac{\partial}{\partial v}\left(\frac{\partial \mathrm{C}(u,v)}{\partial u}\right) = exp\left(-\frac{1}{\alpha}\left[(-\alpha \ln(u))^{\frac{1}{\alpha}} + (-\alpha \ln v)^{\frac{1}{\alpha}}\right]^{\alpha}\right)\left(\frac{-1}{\alpha}\right)(\alpha)\left[(-\alpha \ln(u))^{\frac{1}{\alpha}} + (-\alpha \ln v)^{\frac{1}{\alpha}}\right]^{\alpha-1}$$

$$\left(\frac{1}{\alpha}(-\alpha \ln(v))^{\frac{1}{\alpha}-1}\left(\frac{-\alpha}{v}\right)\right) \times \left[(-\alpha \ln(u))^{\frac{1}{\alpha}} + (-\alpha \ln v)^{\frac{1}{\alpha}}\right]^{\alpha-1}\left((-\alpha \ln(u))^{\frac{1}{\alpha}-1}\left(\frac{1}{u}\right)\right) +$$

$$exp\left(-\frac{1}{\alpha}\left[(-\alpha \ln(u))^{\frac{1}{\alpha}} + (-\alpha \ln v)^{\frac{1}{\alpha}}\right]^{\alpha}\right)\left((-\alpha \ln(u))^{\frac{1}{\alpha}-1}\left(\frac{1}{u}\right)\right)(\alpha-1)\left[(-\alpha \ln(u))^{\frac{1}{\alpha}} + (-\alpha \ln v)^{\frac{1}{\alpha}}\right]^{\alpha-2}\left(\frac{1}{\alpha}(-\alpha \ln(v))^{\frac{1}{\alpha}-1}\left(\frac{-\alpha}{v}\right)\right)$$

$$= exp\left(-\frac{1}{\alpha}\left[(-\alpha \ln(u))^{\frac{1}{\alpha}} + (-\alpha \ln v)^{\frac{1}{\alpha}}\right]^{\alpha}\right)\left[(-\alpha \ln(u))^{\frac{1}{\alpha}} + (-\alpha \ln v)^{\frac{1}{\alpha}}\right]^{2\alpha-2}\left((-\alpha \ln(v))^{\frac{1}{\alpha}-1}\left(\frac{1}{v}\right)\right) \times$$

$$\left((-\alpha \ln(u))^{\frac{1}{\alpha}-1}\left(\frac{1}{u}\right)\right) + exp\left(-\frac{1}{\alpha}\left[(-\alpha \ln(u))^{\frac{1}{\alpha}} + (-\alpha \ln v)^{\frac{1}{\alpha}}\right]^{\alpha}\right)\left((-\alpha \ln(u))^{\frac{1}{\alpha}-1}\left(\frac{1}{u}\right)\right) \times$$

$$(\alpha-1)\left[(-\alpha \ln(u))^{\frac{1}{\alpha}} + (-\alpha \ln v)^{\frac{1}{\alpha}}\right]^{\alpha-2}\left((-\alpha \ln(v))^{\frac{1}{\alpha}-1}\left(\frac{-1}{v}\right)\right)$$

$$= exp\left(-\frac{1}{\alpha}\left[(-\alpha\ln(u))^{\frac{1}{\alpha}} + (-\alpha\ln v)^{\frac{1}{\alpha}}\right]^{\alpha}\right)\left[(-\alpha\ln(u))^{\frac{1}{\alpha}} + (-\alpha\ln v)^{\frac{1}{\alpha}}\right]^{\alpha-2}\left[\frac{(\alpha^2\ln v\ln u)^{\frac{1}{\alpha}-1}}{vu}\right]\left\{\left[(-\alpha\ln(u))^{\frac{1}{\alpha}} + (-\alpha\ln v)^{\frac{1}{\alpha}}\right]^{\alpha} - (\alpha-1)\right\}$$

$$= exp\left(-\frac{1}{\alpha}\left[(-\alpha\ln(u))^{\frac{1}{\alpha}} + (-\alpha\ln v)^{\frac{1}{\alpha}}\right]^{\alpha}\right)\left[(-\alpha\ln(u))^{\frac{1}{\alpha}} + (-\alpha\ln v)^{\frac{1}{\alpha}}\right]^{\alpha-2}\left[\frac{(\alpha^2\ln v\ln u)^{\frac{1}{\alpha}-1}}{vu}\right]\left\{\left[(-\alpha\ln(u))^{\frac{1}{\alpha}} + (-\alpha\ln v)^{\frac{1}{\alpha}}\right]^{\alpha} + 1 - \alpha\right\}$$

The second derivative is positive at alpha defined on the interval $0 < \alpha \le 1$ and this copula density integrate to one using MCMC simulation from zero to one, so it is a valid density function.

**Proposition 3**: this copula is absolutely continuous copula and it has no singular part.

**Proof** : to test for singularity: $\frac{\varphi(0)}{\varphi'(0)} = 0$

$$\frac{\varphi(u)}{\varphi'(u)} = \frac{(-\alpha\ln u)^{\frac{1}{\alpha}}}{\frac{1}{\alpha}(-\alpha\ln u)^{\frac{1}{\alpha}-1}\left(\frac{-\alpha}{u}\right)} = \frac{-\alpha\ln u}{\frac{-1}{u}} \quad and \quad \lim_{u\to 0}\frac{\varphi(u)}{\varphi'(u)} = \frac{(-\alpha\ln u)^{\frac{1}{\alpha}}}{\frac{1}{\alpha}(-\alpha\ln u)^{\frac{1}{\alpha}-1}\left(\frac{-\alpha}{u}\right)} = \frac{-\alpha\ln u}{\frac{-1}{u}}$$

Using L'Hopital : $\lim_{u\to 0}\frac{\varphi(u)}{\varphi'(u)} = \lim_{u\to 0}\frac{(-\alpha\ln u)^{\frac{1}{\alpha}}}{\frac{1}{\alpha}(-\alpha\ln u)^{\frac{1}{\alpha}-1}\left(\frac{-\alpha}{u}\right)} = \lim_{u\to 0}\frac{\frac{\alpha}{u}}{\frac{-1}{u^2}} = \lim_{u\to 0}(-\alpha u) = -\alpha(0) = 0$

As long as this limit is zero at u=0 so it has no singular part and it is absolutely continuous copula.

**Proposition 4**: Kendall tau for this copula is $\tau = 1 - \alpha$ and $0 < \alpha \le 1$

**Proof:**

$$\tau = 4\int_0^1 \frac{\varphi(u)}{\varphi'(u)}du + 1 \qquad and \qquad \int_0^1 \frac{\varphi(u)}{\varphi'(u)}\, du = \int_0^1 \alpha\, u\ln(u)\, du = \frac{-\alpha}{4}\ , integration\ by\ parts$$

$$\tau = 4\int_0^1 \frac{\varphi(u)}{\varphi'(u)}du + 1 = 4\left(\frac{-\alpha}{4}\right) + 1 = 1 - \alpha$$

If $\alpha = 1$ , *then* $\tau = 0$ indicating independence.

**Preposition 5**: This copula is product copula at alpha =1, in other word, if alpha=1, $C(u,v) = uv$

**Proof:**

$$C(u,v) = exp\left(-\frac{1}{\alpha}\left[(-\alpha\ln u)^{\frac{1}{\alpha}} + (-\alpha\ln v)^{\frac{1}{\alpha}}\right]^{\alpha}\right) = exp(-[-\ln u - \ln v]) = exp([\ln u + \ln v]) = exp([\ln u\, v]) = uv$$

Figures (1-3) illustrate the contour plot and the surface of the copula density and distribution with different values of the dependency parameter.

The limitation for this copula is that it does not cover the negative dependency between variables. This copula models only positive dependency.

**Preposition 6:** This copula has no Lower tail dependency and the maximum upper tail dependency coefficient is one.

**Proof:**

$$\lambda_U = \lim_{u\to 1}\frac{1 - 2u + C(u,u)}{1-u} = \frac{1 - 2u + exp\left(-\frac{1}{\alpha}\left[2(-\alpha\ln u)^{\frac{1}{\alpha}}\right]^{\alpha}\right)}{1-u} = \frac{1 - 2u + exp\left(-\frac{1}{\alpha}\ 2^{\alpha}(-\alpha\ln u)\right)}{1-u} = \frac{1 - 2u + u^{2^{\alpha}}}{1-u}$$

Using L'Hopital rule

$$\lambda_U = \lim_{u\to 1}\frac{1-2u+u^{2^\alpha}}{1-u} = \lim_{u\to 1}\frac{-2+2^\alpha u^{(2^\alpha-1)}}{-1} = \quad 2-2^\alpha \in (0,1] \quad for\ \forall\alpha \in (0,1) \qquad (2)$$

$$\lambda_L = \lim_{u\to 0}\frac{C(u,u)}{u} = \lim_{u\to 0}\frac{exp\left(-\frac{1}{\alpha}\ 2^\alpha(-\alpha\ln u)\right)}{u} = \lim_{u\to 0}\frac{u^{2^\alpha}}{u} = u^{(2^\alpha-1)} = 0 \quad for\ \forall\alpha \in (0,1) \qquad (3)$$

So the maximum upper tail dependence coefficient for $\forall\alpha \in (0,1)$ of this copula is 1, however the copula has no lower tail dependency as the lower tail dependence coefficient is zero for $\forall\alpha \in (0,1)$ . The $\lambda_U = 1$ , when dependency parameter approaches 0. While $\lambda_U = 0$ , when this parameter is 1, which is the product copula (state of independence). The $\lambda_U \in (0,1]$ when the dependency parameter $\alpha \in (0,1)$ . The copula has upper tail dependency for all theta less than one, and this dependence increases as theta decreases. Also the following equations support these findings.

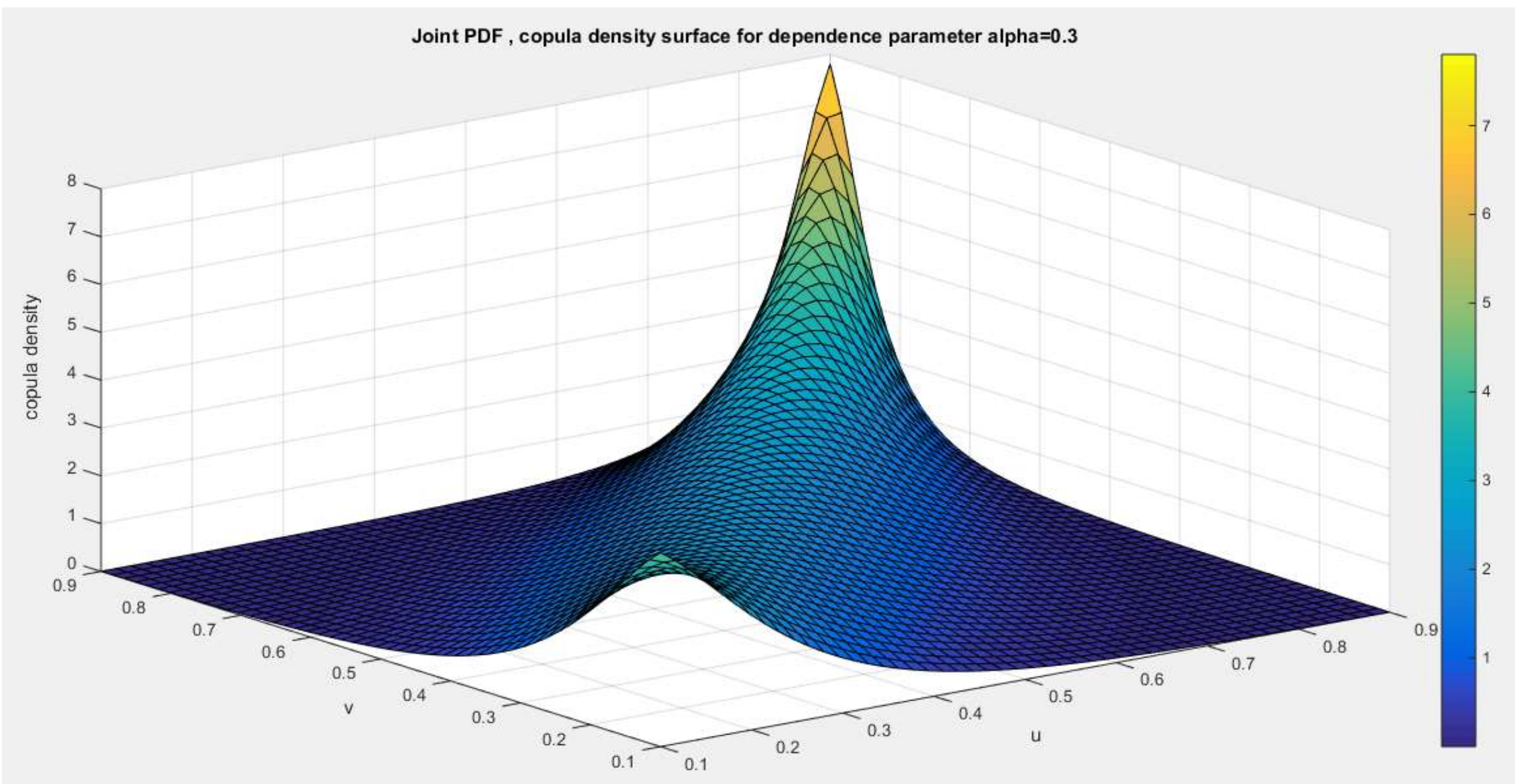

Fig.1 shows the first copula density surface with dependency parameter 0.3

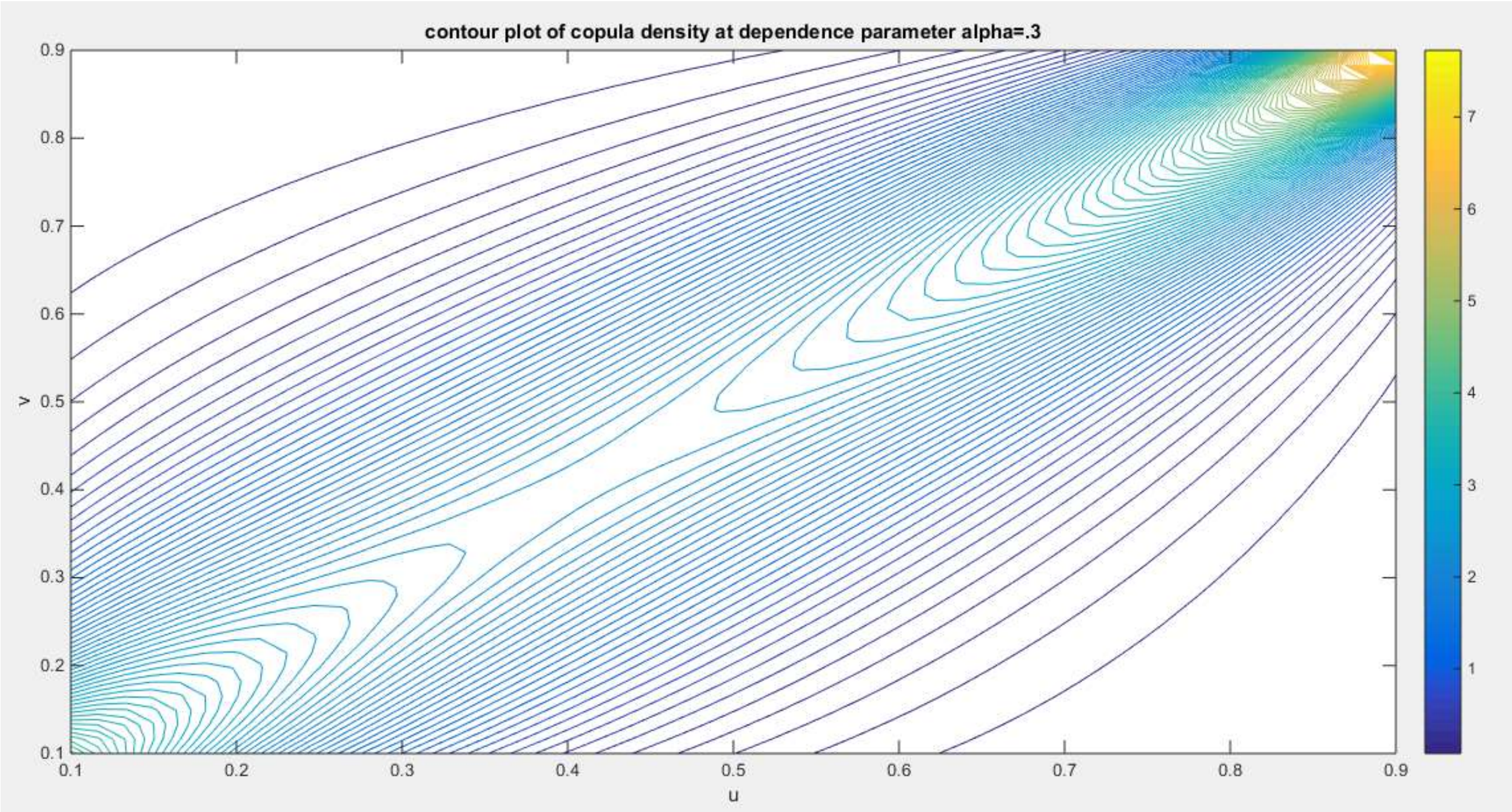

Fig.2 shows the contour plot of the first copula density with dependency parameter 0.3

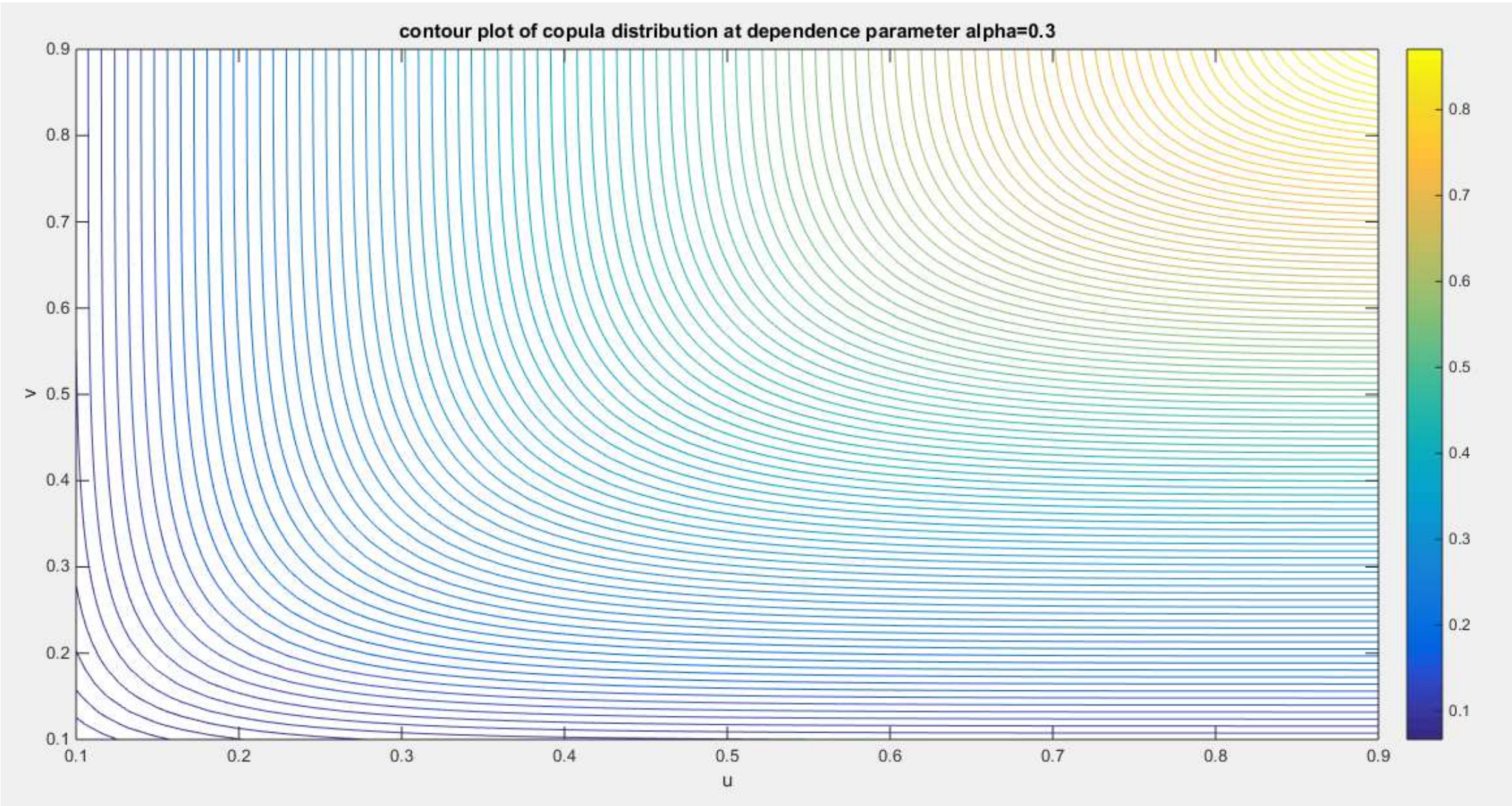

Fig. 3 shows the contour plot of the first copula distribution with dependency parameter 0.3

**Upper tail dependency:** $\lambda_U = 2 - 2.\lim_{u\to 0} \frac{\varphi'(2u)}{\varphi'(u)} = 1$

$$\lim_{u\to 0} \frac{\varphi'(2u)}{\varphi'(u)} = \lim_{u\to 0} \frac{\frac{1}{\alpha}(-\alpha \ln 2u)^{\frac{1}{\alpha}-1} \left(\frac{-\alpha}{2u}\right)}{\frac{1}{\alpha}(-\alpha \ln u)^{\frac{1}{\alpha}-1} \left(\frac{-\alpha}{u}\right)} = \lim_{u\to 0} \frac{(-\alpha \ln 2u)^{\frac{1}{\alpha}-1} \left(\frac{1}{2}\right)}{(-\alpha \ln u)^{\frac{1}{\alpha}-1}} = \lim_{u\to 0} \left(\frac{1}{2}\right)\left[\frac{\frac{2}{2u}}{\frac{1}{u}}\right]^{\frac{1}{\alpha}-1} = \left(\frac{1}{2}\right)$$

$$\lambda_U = 2 - 2.\lim_{u\to 0} \frac{\varphi'(2u)}{\varphi'(u)} = 2 - 2\left(\frac{1}{2}\right) = 1 \tag{4}$$

Equation (4) is complementary to equation (2), as in equation (2) the upper limit for the upper tail coefficient cannot be attained, equation (4) explicitly satisfies this limit when the dependency parameter is approaching 0.

**Lower tail dependency**: $\lambda_L = 2 - 2.\lim_{u\to 0} \frac{\varphi'(u)}{\varphi'(2u)} = 0$

$$\lim_{u\to 0} \frac{\varphi'(u)}{\varphi'(2u)} = \lim_{u\to 0} \frac{\frac{1}{\alpha}(-\alpha \ln u)^{\frac{1}{\alpha}-1}\left(\frac{-\alpha}{u}\right)}{\frac{1}{\alpha}(-\alpha \ln 2u)^{\frac{1}{\alpha}-1}\left(\frac{-\alpha}{2u}\right)} = \lim_{u\to 0} \frac{(-\alpha \ln u)^{\frac{1}{\alpha}-1}\left(\frac{-1}{u}\right)}{(-\alpha \ln 2u)^{\frac{1}{\alpha}-1}\left(\frac{-1}{2u}\right)} = \lim_{u\to 0} (2)\left[\frac{-\alpha \ln u}{-\alpha \ln 2u}\right]^{\frac{1}{\alpha}-1} = \lim_{u\to 0} (2)\left[\frac{\ln u}{\ln 2u}\right]^{\frac{1}{\alpha}-1}$$

Using L'Hopital rule : $\lim_{u\to 0} \ (2)\left[\frac{\frac{1}{u}}{\frac{2}{2u}}\right]^{\frac{1}{\alpha}-1} = \lim_{u\to 0} \ (2)\,[1]^{\frac{1}{\alpha}-1} = 2$

$\lambda_L = 2 - 2.\lim_{u\to 0} \frac{\varphi'(u)}{\varphi'(2u)} = 2 - 2(2) = -2$ but the lower dependency tail coefficient should be between 0 and 1 so lower tail dependency coefficient is 0, and this is consistent with equation (3).

## Section 2: (Second Copula)

Let's say the inverse generator is $C(\varphi(u) + \varphi(v)) = \Phi^{-1}(t) = z = exp(-t^{\alpha^2})$ , so the generator is $\ln z = \ln\{exp(-t^{\alpha^2})\} = -t^{\alpha^2}$ and hence, the $-\ln z = t^{\alpha^2}$. The generator $\varphi(z) = (-\ln z)^{\alpha^{-2}} = t$ , $z \in (0,1)$. And this generator should fulfill the sufficient conditions: $\varphi(0) = (-\ln 0)^{\alpha^{-2}} = \infty$ , $\varphi(1) = (-\ln 1)^{\alpha^{-2}} = 0$ and the first derivative $\varphi'(z) = \alpha^{-2}(-\ln z)^{\alpha^{-2}-1}\left(\frac{-1}{z}\right) < 0$ , ensures that the generator is a decreasing function. While the positive second derivative ensures that the generator is a convex function at $-1 \le \alpha \le 1$ excluding 0 i.e. $\alpha \in [-1,1] - \{0\}$, $\varphi''(z) = \frac{1}{\alpha^2}\left(\frac{1}{\alpha^2} - 1\right)(-\ln z)^{\alpha^{-2}-2}\left(\frac{-1}{z}\right)\left(\frac{-1}{z}\right) + \alpha^{-2}(-\ln z)^{\alpha^{-2}-1}\left(\frac{1}{z^2}\right) > 0$ .

For bivariate distribution: $\varphi(u) = (-\ln u)^{\frac{1}{\alpha^2}}$ and $\varphi(v) = (-\ln v)^{\frac{1}{\alpha^2}}$

$$\varphi(w) = \varphi(u) + \varphi(v) = (-\ln u)^{\frac{1}{\alpha^2}} + (-\ln v)^{\frac{1}{\alpha^2}}$$

$$C(\varphi(u) + \varphi(v)) = exp\left(-\left[(-\ln u)^{\frac{1}{\alpha^2}} + (-\ln v)^{\frac{1}{\alpha^2}}\right]^{\alpha^2}\right)$$

$$C(u, v) = C(F_X(x), F_Y(y))$$

For a copula to be a valid copula it should fulfill the boundary condition, the marginal uniformity and 2-increasing condition. And this is equivalent for necessary conditions of the inverse generator to be fulfilled which are the following:

**Proposition 7**: The boundary conditions are $C(0, v) = C(u, 0) = 0$ , and the marginal uniformity conditions are $C(1, v) = v$ *and* $C(u, 1) = u$ .

**Proof**:

$$C(1, v) = exp\left(-\left[(-\ln 1)^{\frac{1}{\alpha^2}} + (-\ln v)^{\frac{1}{\alpha^2}}\right]^{\alpha^2}\right) = exp\left(-\left[0 + (-\ln v)^{\frac{1}{\alpha^2}}\right]^{\alpha^2}\right) = exp\left(-\left[(-\ln v)^{\frac{1}{\alpha^2}}\right]^{\alpha^2}\right) = exp(-[-\ln v])$$

$C(1, v) = v$ and $C(u, 1) = exp\left(-\left[(-\ln u)^{\frac{1}{\alpha^2}}\right]^{\alpha^2}\right) = exp(-[-\ln u]) = u$

**Proposition 8**: for copula to be valid it should be 2-increasing, in other words; $\frac{\partial^2 C(u,v)}{\partial u \partial v} \geq 0$

**Proof**

$$\frac{\partial C(u,v)}{\partial u} = exp\left(-\left[(-\ln u)^{\frac{1}{\alpha^2}} + (-\ln v)^{\frac{1}{\alpha^2}}\right]^{\alpha^2}\right)(-\alpha^2)\left[(-\ln u)^{\frac{1}{\alpha^2}} + (-\ln v)^{\frac{1}{\alpha^2}}\right]^{\alpha^2-1}\left(\frac{1}{\alpha^2}(-\ln(u))^{\frac{1}{\alpha^2}-1}\left(\frac{-1}{u}\right)\right)$$

$$\frac{\partial C(u,v)}{\partial u} = exp\left(-\left[(-\ln u)^{\frac{1}{\alpha^2}} + (-\ln v)^{\frac{1}{\alpha^2}}\right]^{\alpha^2}\right)\left[(-\ln u)^{\frac{1}{\alpha^2}} + (-\ln v)^{\frac{1}{\alpha^2}}\right]^{\alpha^2-1}\left((-\ln(u))^{\frac{1}{\alpha^2}-1}\left(\frac{1}{u}\right)\right)$$

$$\frac{\partial}{\partial v}\left(\frac{\partial C(u,v)}{\partial u}\right) = exp\left(-\left[(-\ln u)^{\frac{1}{\alpha^2}} + (-\ln v)^{\frac{1}{\alpha^2}}\right]^{\alpha^2}\right) \times (-\alpha^2)\left[(-\ln u)^{\frac{1}{\alpha^2}} + (-\ln v)^{\frac{1}{\alpha^2}}\right]^{\alpha^2-1}$$

$$\left(\frac{1}{\alpha^2}(-\ln(v))^{\frac{1}{\alpha^2}-1}\left(\frac{-1}{v}\right)\right) \times \left[(-\ln u)^{\frac{1}{\alpha^2}} + (-\ln v)^{\frac{1}{\alpha^2}}\right]^{\alpha^2-1}\left((-\ln(u))^{\frac{1}{\alpha^2}-1}\left(\frac{1}{u}\right)\right) +$$

$$exp\left(-\left[(-\ln u)^{\frac{1}{\alpha^2}} + (-\ln v)^{\frac{1}{\alpha^2}}\right]^{\alpha^2}\right)\left((-\ln(u))^{\frac{1}{\alpha^2}-1}\left(\frac{1}{u}\right)\right) \times (\alpha^2 - 1)\left[(-\ln u)^{\frac{1}{\alpha^2}} + (-\ln v)^{\frac{1}{\alpha^2}}\right]^{\alpha^2-2}\left(\frac{1}{\alpha^2}(-\ln(v))^{\frac{1}{\alpha^2}-1}\left(\frac{-1}{v}\right)\right)$$

$$= exp\left(-\left[(-\ln u)^{\frac{1}{\alpha^2}} + (-\ln v)^{\frac{1}{\alpha^2}}\right]^{\alpha^2}\right)\left[(-\ln u)^{\frac{1}{\alpha^2}} + (-\ln v)^{\frac{1}{\alpha^2}}\right]^{\alpha^2-1}\left((-\ln(v))^{\frac{1}{\alpha^2}-1}\left(\frac{1}{v}\right)\right)\left[(-\ln u)^{\frac{1}{\alpha^2}} + (-\ln v)^{\frac{1}{\alpha^2}}\right]^{\alpha^2-1}$$

$$\left((-\ln(u))^{\frac{1}{\alpha^2}-1}\left(\frac{1}{u}\right)\right) +$$

$$exp\left(-\left[(-\ln u)^{\frac{1}{\alpha^2}} + (-\ln v)^{\frac{1}{\alpha^2}}\right]^{\alpha^2}\right)\left((-\ln(u))^{\frac{1}{\alpha^2}-1}\left(\frac{1}{u}\right)\right) \times (\alpha^2 - 1)\left[(-\ln u)^{\frac{1}{\alpha^2}} + (-\ln v)^{\frac{1}{\alpha^2}}\right]^{\alpha^2-2}\left(\frac{1}{\alpha^2}(-\ln(v))^{\frac{1}{\alpha^2}-1}\left(\frac{-1}{v}\right)\right)$$

$$= exp\left(-\left[(-\ln u)^{\frac{1}{\alpha^2}} + (-\ln v)^{\frac{1}{\alpha^2}}\right]^{\alpha^2}\right)\left[(-\ln u)^{\frac{1}{\alpha^2}} + (-\ln v)^{\frac{1}{\alpha^2}}\right]^{2\alpha^2-2}\left[\frac{(\ln(v)\ln(u))^{\frac{1}{\alpha^2}-1}}{uv}\right]$$

$$+ (\alpha^2 - 1)\left(\frac{-1}{\alpha^2}\right)\left[\frac{(\ln(v)\ln(u))^{\frac{1}{\alpha^2}-1}}{uv}\right] exp\left(-\left[(-\ln u)^{\frac{1}{\alpha^2}} + (-\ln v)^{\frac{1}{\alpha^2}}\right]^{\alpha^2}\right) \times \left[(-\ln u)^{\frac{1}{\alpha^2}} + (-\ln v)^{\frac{1}{\alpha^2}}\right]^{\alpha^2-2}$$

$$= exp\left(-\left[(-\ln u)^{\frac{1}{\alpha^2}} + (-\ln v)^{\frac{1}{\alpha^2}}\right]^{\alpha^2}\right)\left[(-\ln u)^{\frac{1}{\alpha^2}} + (-\ln v)^{\frac{1}{\alpha^2}}\right]^{\alpha^2-2}\left[\frac{(\ln(v)\ln(u))^{\frac{1}{\alpha^2}-1}}{uv}\right]\left\{\left[(-\ln u)^{\frac{1}{\alpha^2}} + (-\ln v)^{\frac{1}{\alpha^2}}\right]^{\alpha^2} - \left(\frac{1}{\alpha^2}\right)(\alpha^2 - 1)\right\}$$

The second derivative is positive at $\alpha \in [-1,1] - \{0\}$ . Integration of this copula density using MCMC simulation from zero to one is one and so it is a valid copula density. Figures (4-6) show the surface and the contour plot of the copula density and distribution at different dependency parameters.

**Preposition 9**: this copula is absolutely continuous copula and it has no singular part.

**Proof** : to test for singularity: $\frac{\varphi(0)}{\varphi'(0)} = 0$

$$\frac{\varphi(u)}{\varphi'(u)} = \frac{(-\ln u)^{\frac{1}{\alpha^2}}}{\frac{1}{\alpha^2}(-\ln u)^{\frac{1}{\alpha^2}-1}\left(\frac{-1}{u}\right)} = \lim_{u\to 0}\frac{(-\ln u)^{\frac{1}{\alpha^2}}}{\frac{1}{\alpha^2}(-\ln u)^{\frac{1}{\alpha^2}-1}\left(\frac{-1}{u}\right)} = \lim_{u\to 0}\frac{\alpha^2(-\ln u)}{\frac{-1}{u}} = \lim_{u\to 0} -\alpha^2 u(-\ln u)$$

Using L'Hopital rule: $\lim_{u\to 0}\frac{\varphi(u)}{\varphi'(u)} = \frac{\alpha^2\left(\frac{-1}{u}\right)}{\frac{1}{u^2}} = -\alpha^2 u = -\alpha^2(0) = 0$

As long as this limit is zero at u=0 so the copula has no singular part and it is absolutely continuous copula.

**Preposition 10**: Kendall tau for this copula is $\tau = 1 - \alpha^2$ and $\alpha \in [-1,1] - \{0\}$

**Proof**: using integration by parts for the following integrant

$$\tau = 4\int_0^1 \frac{\varphi(u)}{\varphi'(u)}du \ +1 = \ 4\int_0^1 (-\alpha^2 u(-\ln u))\,du \ +1 \quad and \quad \int_0^1 \frac{\varphi(u)}{\varphi'(u)}\,du = \int_0^1 -\alpha^2 u(-\ln u)\ \ du \ = \frac{-\alpha^2}{4}$$

$$\tau = 4\int_0^1 \frac{\varphi(u)}{\varphi'(u)}du \ +1 = \ 4\left(\frac{-\alpha^2}{4}\right) \ +1 = 1-\alpha^2$$

If $\alpha$ *is approaching* 0 $\tau = 1$ indicating perfect positive dependency.

If $\alpha \in \{1,-1\}$ *so* $\tau = 0$ denoting independency with product copula $C(\varphi(u)+\varphi(v)) = uv$.

**Preposition 11**: The product copula is attained at dependency parameter $\alpha \in \{1,-1\}$

**Proof**: when $\alpha \in \{1,-1\}$

$$C(\varphi(u)+\varphi(v)) \ = exp\left(-\left[(-\ln u)^{\frac{1}{\alpha^2}} + (-\ln v)^{\frac{1}{\alpha^2}}\right]^{\alpha^2}\right) = exp(-[(-\ln u) + (-\ln v)]) = uv$$

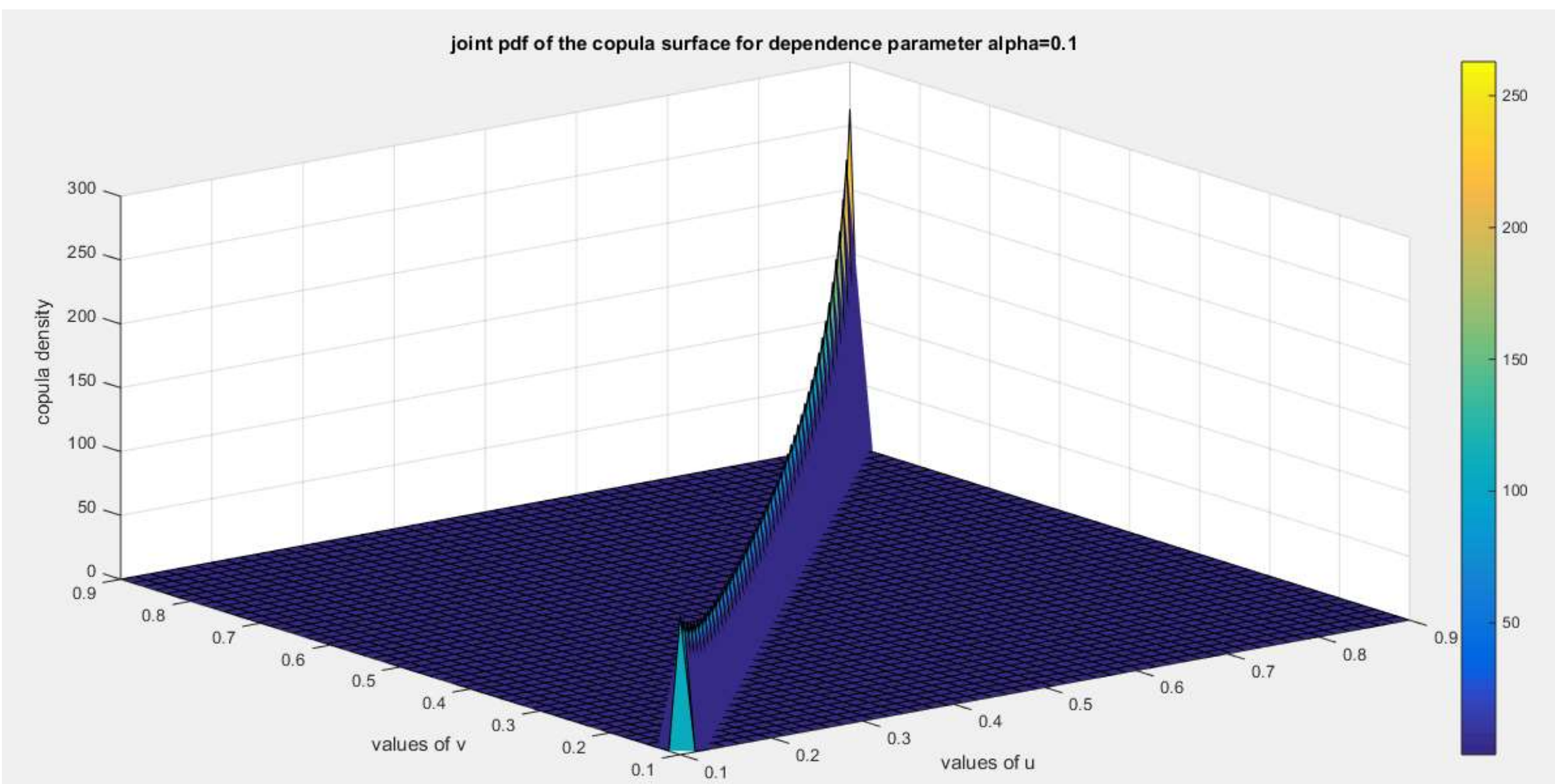

Fig. 4 shows second copula density surface with dependency parameter 0.1

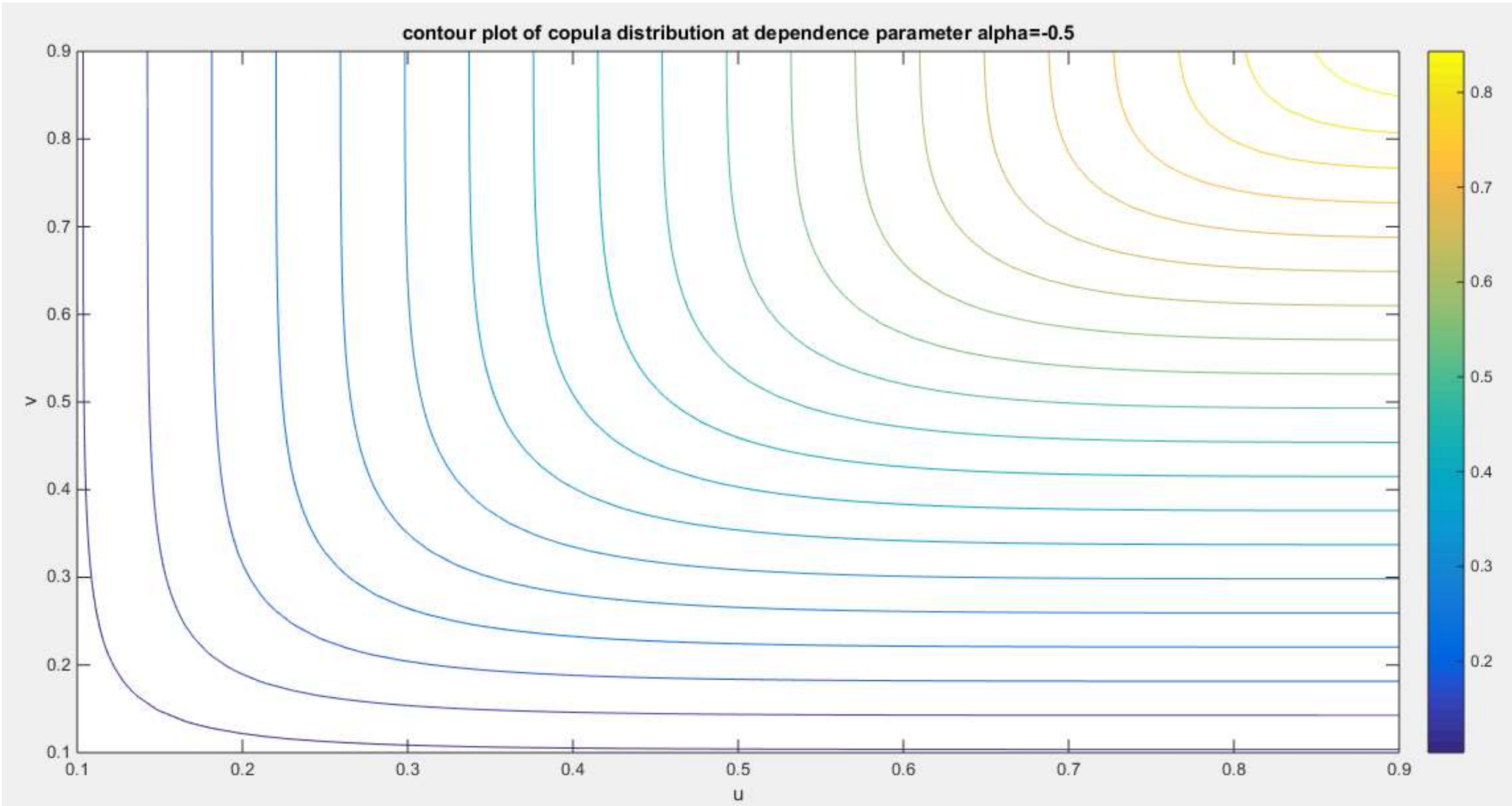

Fig. 5 shows the contour plot of the second copula distribution with dependency parameter -0.5

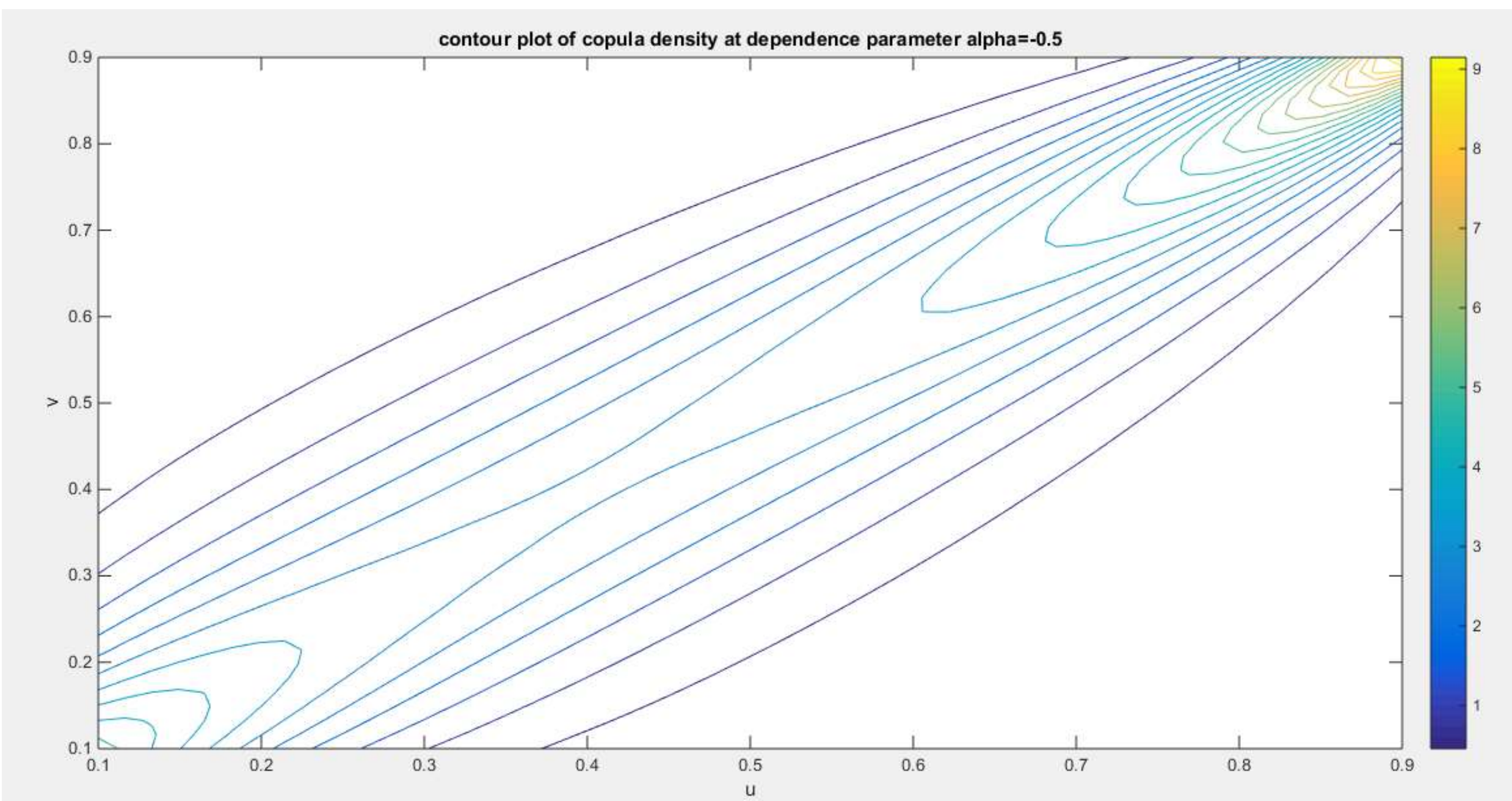

Fig. 6 shows the contour plot of the second copula distribution with dependency parameter -0.5

The limitation of this copula is that it models only positive dependency.

**Preposition 12:** The copula has no lower tail dependency and the maximum upper tail dependency coefficient is 1.

**Proof:**

$$\lambda_U = \lim_{u\to 1}\frac{1-2u+C(u,u)}{1-u} = \frac{1-2u+exp\left(-\left[2(-\ln u)^{\frac{1}{\alpha^2}}\right]^{\alpha^2}\right)}{1-u} = \frac{1-2u+exp\left(-\ 2^{\alpha^2}(-\alpha\ln u)\right)}{1-u} = \frac{1-2u+u^{2^{\alpha^2}}}{1-u}$$

Using L'Hopital rule

$$\lambda_U = \lim_{u\to 1}\frac{1-2u+u^{2^{\alpha^2}}}{1-u} = \lim_{u\to 1}\frac{-2+2^{\alpha\alpha^2}u^{\left(2^{\alpha^2}-1\right)}}{-1} = \ 2-2^{\alpha^2} \in (0,1] \quad for\ \forall\alpha \in (-1,1)-\{0\} \qquad (5)$$

$$\lambda_L = \lim_{u\to 0}\frac{C(u,u)}{u} = \lim_{u\to 0}\frac{exp\left(-\ 2^{\alpha^2}(-\ln u)\right)}{u} = \lim_{u\to 0}\frac{u^{2^{\alpha^2}}}{u} = u^{\left(2^{\alpha^2}-1\right)} = 0 \quad for\ \forall\alpha \in (-1,1)-\{0\} \qquad (6)$$

So the maximum upper tail dependence coefficient for $\forall\alpha \in (-1,1)$ of this copula is 1, however the copula has no lower tail dependency as the lower tail dependence coefficient is zero for $\forall\alpha \in (-1,1)-\{0\}$ . The $\lambda_U = 1$ , when dependency parameter approaches 0. While $\lambda_U = 0$ , when this parameter is 1 or -1, at which the product copula is produced. The $\lambda_U \in (0,1)$ when the dependency parameter $\alpha \in (-1,1)-\{0\}$ . Also the following equations support these findings.

**Upper tail dependency:** $\lambda_U = 2-2.\lim_{u\to 0}\frac{\varphi'(2u)}{\varphi'(u)} = 1$

$$\lim_{u\to 0}\frac{\varphi'(2u)}{\varphi'(u)} = \lim_{u\to 0}\frac{\alpha^{-2}(-\ln 2u)^{\alpha^{-2}-1}\ \left(\frac{-1}{2u}\right)}{\alpha^{-2}(-\ln u)^{\frac{1}{\alpha^2}-1}\ \left(\frac{-1}{u}\right)} = \lim_{u\to 0}\left(\frac{1}{2}\right)\frac{(-\ln 2u)^{\alpha^{-2}-1}}{(-\ln u)^{\alpha^{-2}-1}} = \lim_{u\to 0}\left(\frac{1}{2}\right)\left[\frac{-\ln 2\,u}{-\ln u}\right]^{\frac{1}{\alpha^2}-1}$$

$$= \lim_{u\to 0}\left(\frac{1}{2}\right)\left[\frac{\ln 2\,u}{\ln u}\right]^{\frac{1}{\alpha^2}-1}, \text{ Using L'Hopital rule: } \lim_{u\to 0}\left(\frac{1}{2}\right)\left[\frac{\frac{2}{2u}}{\frac{1}{u}}\right]^{\frac{1}{\theta}-1} = \left(\frac{1}{2}\right)$$

$$\lambda_U = 2-2.\lim_{u\to 0}\frac{\varphi'(2u)}{\varphi'(u)} = 2-2\left(\frac{1}{2}\right) = 1 \qquad (7)$$

Equation (5) is complementary to equation (7), as in equation (5) the upper limit for the upper tail coefficient cannot be explicitly attained, equation (7) explicitly satisfies this limit when the dependency parameter is approaching 0.

**Lower tail dependency:** $\lambda_L = 2-2.\lim_{u\to 0}\frac{\varphi'(u)}{\varphi'(2u)} = 0$

$$\lim_{u\to 0}\frac{\varphi'(u)}{\varphi'(2u)} = \lim_{u\to 0}\frac{\alpha^{-2}(-\ln u)^{\alpha^{-2}-1}\ \left(\frac{-1}{u}\right)}{\alpha^{-2}(-\ln 2u)^{\frac{1}{\alpha^2}-1}\ \left(\frac{-1}{2u}\right)} = \lim_{u\to 0}\frac{(-\ln u)^{\alpha^{-2}-1}\ (2)}{(-\ln 2u)^{\alpha^{-2}-1}} = \lim_{u\to 0}\ (2)\left[\frac{-\ln u}{-\ln 2u}\right]^{\frac{1}{\alpha^2}-1}$$

$$= \lim_{u\to 0}\ (2)\left[\frac{\ln u}{\ln 2u}\right]^{\frac{1}{\alpha^2}-1}, \text{ Using L'Hopital rule: } \lim_{u\to 0}\ (2)\left[\frac{\frac{1}{u}}{\frac{2}{2u}}\right]^{\frac{1}{\alpha^2}-1} = \lim_{u\to 0}\ (2)\,[1]^{\frac{1}{\alpha^2}-1} = 2$$

$\lambda_L = 2-2.\lim_{u\to 0}\frac{\varphi'(u)}{\varphi'(2u)} = 2-2(2) = -2$ but the lower dependency tail should be between 0 and 1 so lower dependency tail is 0, and this is consistent with equation (6)

## Section 3: Methods of Estimation

(Weiß, 2011), (Chen, 2007), and (Tsukahara, 2005) extensively introduced various methods of estimations like the following:

### 3.1. Maximum likelihood estimation (MLE)

It is one step where all the parameters marginal and copula are estimated for the full joint density function, taking the log-likelihood of this full joint function and differentiating it with respect to each of the marginal parameters and the copula parameter, then equating the resultant functions with zero and solving numerically this system of equations.

### 3.2. Inference function of margins (IFM)

For bivariate copula-base models for which the maximum likelihood is computationally difficult to estimate all the model parameters jointly in a one-step parametric method (Al Turk et al., 2017), two-stages estimation procedure has been proposed previously by (Joe, 2005) and (Xu, 2009). The first stage involves estimating the marginal parameters separately using maximum likelihood as shown in equation (8), where s stands for the score function or gradient

$$s_i^{(1)} = \frac{\partial}{\partial\theta_1}\ln f_X\left(x_i;\theta_{marginal\ 1}\right) \quad and \quad s_i^{(2)} = \frac{\partial}{\partial\theta_2}\ln f_Y\left(y_i;\theta_{marginal\ 2}\right) \tag{8}$$

After estimating the marginal parameter, you get the marginal CDF for each variable as shown in equation (9)

$$u_i = \hat{F}_1(x_i;\hat{\theta}_1) \quad and \quad v_i = \hat{F}_2(y_i;\hat{\theta}_2) \tag{9}$$

The second stage involves maximum likelihood of the dependence parameters with the marginal parameters held fixed from the first stage, you get the log-likelihood of the copula density as shown in equation (10), then differentiate this log-likelihood with respect to the copula dependency parameter as shown in equation (11).

$$\ln L_{(IFM)} = \sum_{i=1}^{n} \ln c\left( \hat{F}_1(x_i;\hat{\theta}_1), \hat{F}_2(y_i;\hat{\theta}_2) \quad ;\theta_{copula}\right) = \sum_{i=1}^{n} \ln c\left( u_i, v_i \quad ;\theta_{copula}\right) \tag{10}$$

$$\frac{\partial \ln L_{(IFM)}}{\partial\theta_{copula}} = \frac{\partial}{\partial\theta_{copula}}\left[\sum_{i=1}^{n} \ln c\left( \hat{F}_1(x_i;\hat{\theta}_1), \hat{F}_2(y_i;\hat{\theta}_2) \quad ;\theta_{copula}\right)\right] \tag{11}$$

In this paper, the author used IFM for both the first and second copula.

For the first copula: first derivative of the log-likelihood with respect to copula parameter is in equation (12)

Let: $t_1 = (-\alpha\ln u_i)^{1/\alpha}$ , $t_2 = (-\alpha\ln v_i)^{1/\alpha}$ , $t_3 = (-\alpha\ln u_i)^{\frac{1}{\alpha}-1}$ , $t_4 = (-\alpha\ln v_i)^{\frac{1}{\alpha}-1}$

$h = \left\{\frac{1}{\alpha}\, t_3\,(-\ln u_i) + t_1\left(\frac{-1}{\alpha^2}\right)\ln(-\alpha\ln u_i) + \frac{1}{\alpha}\, t_4\,(-\ln v_i) + t_2\left(\frac{-1}{\alpha^2}\right)\ln(-\alpha\ln v_i)\right\}$,

$g = \alpha(t_1+t_2)^{\alpha-1}h + (t_1+t_2)^{\alpha}\ln(t_1+t_2)$

$A_1 = \frac{1}{\alpha^2}(t_1+t_2)^{\alpha} + \left(\frac{-1}{\alpha}\right)g$ , $A_2 = (\alpha-2)(t_1+t_2)^{-1}h + \ln(t_1+t_2)$ ,

$A_3 = \left(\frac{1}{\alpha}-1\right)(\alpha^2\ln(u_i)\ln(v_i))^{-1}(2\alpha\ln(u_i)\ln(v_i)) + \frac{-1}{\alpha^2}\ln(\alpha^2\ln(u_i)\ln(v_i))$ $\quad A_4 = \frac{g-1}{(t_1+t_2)^{\alpha}-\alpha+1}$

$$\frac{\partial \ln L_{(IFM)}}{\partial\theta_{copula}} = A_1 + A_2 + A_3 + A_4 \tag{12}$$

For second copula: first derivative of the log-likelihood with respect to copula parameter is in equation (13)

Let: $t_1 = (-\ln u_i)^{1/\alpha^2}$ , $t_2 = (-\ln v_i)^{1/\alpha^2}$ , $t_3 = (-\ln u_i)^{1/\alpha^2}\left(\frac{-2}{\alpha^3}\right)\ln(-\ln u_i)$ ,

$t_4 = (-\ln v_i)^{1/\alpha^2}\left(\frac{-2}{\alpha^3}\right)\ln(-\ln v_i)$ , $h = -(\alpha^2)(t_1+t_2)^{\alpha^2-1}(t_3+t_4) + \{-(t_1+t_2)^{\alpha^2}(2\alpha)\ln(t_1+t_2)\}$

$A_1 = h$ , $A_2 = (\alpha^2-2)(t_1+t_2)^{-1}(t_3+t_4) + (2\alpha)\ln(t_1+t_2)$ , $A_3 = \left(\frac{-2}{\alpha^3}\right)\ln(\ln(u_i)\ln(v_i))$

$A_4 = \frac{\left(-h-\frac{2}{\alpha^3}\right)}{\left((t_1+t_2)^{\alpha^2}+\frac{1}{\alpha^2}-1\right)}$

$$\frac{\partial \ln L_{(IFM)}}{\partial \theta_{copula}} = A_1 + A_2 + A_3 + A_4 \tag{13}$$

Applying Quasi-Newton method

$\alpha^{(a+1)} = \alpha^{(a)} + M^{-1}\left(\alpha^{(a)}\right)S(\alpha^a)$ where $S(\alpha) = \frac{\partial \ln L_{(IFM)}}{\partial \theta_{copula}}$ , $M = S^T S$ and $M^{-1}(\alpha)$ can approximate the variance of the dependency parameter estimator.

Each marginal parameter has its own variance from stage one, they are treated as fixed when moving to stage two in IFM approach, so they do not contribute to the variance of the estimated dependency parameter directly. The variance of the copula parameter is obtained from the second stage and it is valid under IFM but it underestimates the total uncertainty (because it ignores the variability in marginal parameter estimation). Therefore, in theory, the propagation of uncertainty from the first stage to the second requires the Godambe information matrix (sandwich estimator)(Godambe, 1991), or simultaneous estimator to obtain the joint variance-covariance matrix of all parameters (marginals and copula). Hence, the sandwich estimator also known as robust variance estimator is used to obtain the asymptotic variance of parameter estimates when the likelihood function is estimated in multiple stages or is mis-specified in cases where not using the full joint likelihood. The variance is defined in the following equation (14):

$$var(\hat{\alpha}) = H^{-1}SH^{-1} \tag{14}$$

Where S is the outer product of the gradient of parameter estimator (marginal and copula) or the score function as shown in equation (15).

$$S = \frac{1}{n}\begin{bmatrix} \frac{\partial \ln f_X\left(x_i;\theta_{marginal}\right)}{\partial \theta_{marginal}} \\ \frac{\partial \ln f_Y\left(y_i;\theta_{marginal}\right)}{\partial \theta_{marginal}} \\ \frac{\partial \ln L_{(IFM)}}{\partial \theta_{copula}} \end{bmatrix}\begin{bmatrix} \frac{\partial \ln f_X\left(x_i;\theta_{1marginal}\right)}{\partial \theta_{1marginal}} & \frac{\partial \ln f_Y\left(y_i;\theta_{2marginal}\right)}{\partial \theta_{2marginal}} & \frac{\partial \ln L_{(IFM)}}{\partial \theta_{copula}} \end{bmatrix} \tag{15}$$

$H^{-1}$ is the inverse of the negative of expectation of the information matrix, in other word, the variance of the parameter estimator (marginal and copula). It is a diagonal matrix. $H = -E\left(\frac{\partial^2 l}{\partial \alpha^2}\right) = I$ , the inverse of this H or inverse of this I is the variance of marginal parameter estimators in the first stage and variance of copula estimator in the second stage. So let's call this inverse the A matrix, so sandwich estimator can be expressed as shown in equation (16-17)

$$var\left(\hat{\theta}\right) = H^{-1}SH^{-1} = ASA \tag{16}$$

$$var - cov\left(\hat{\theta}\right) = \begin{bmatrix} var(\theta_1) & 0 & 0 \\ 0 & var(\theta_2) & 0 \\ 0 & 0 & var\left(\alpha_{copula}\right) \end{bmatrix}\begin{bmatrix} s_{11} & s_{12} & s_{13} \\ s_{21} & s_{22} & s_{23} \\ s_{31} & s_{32} & s_{33} \end{bmatrix}\begin{bmatrix} var(\theta_1) & 0 & 0 \\ 0 & var(\theta_2) & 0 \\ 0 & 0 & var\left(\alpha_{copula}\right) \end{bmatrix}$$

$$var - cov(\hat{\theta}) = \begin{bmatrix} v_1^2 s_{11} & v_1 v_2 s_{12} & v_1 v_3 s_{13} \\ v_1 v_2 s_{21} & v_2^2 s_{22} & v_2 v_3 s_{23} \\ v_1 v_3 s_{31} & v_2 v_3 s_{32} & v_3^2 s_{33} \end{bmatrix} \quad (17)$$

This variance covariance matrix captures the covariance between the marginal and the copula parameter estimators which are not captured in the naive block separated IFM variances.

### 3.3. Semiparametric method (SP)

The third method introduced by (Kim et al., 2007) is the estimation by semiparametric method (SP). It is carried out in two stages like IFM. The first stage involves transforming the variables into their marginal empirical distribution (eCDF), in other words, transforming the observations into pseudo-observations by using the empirical distribution function for each marginal distribution as shown in equation (18)

$$\tilde{F}_i\,(x_i) = \frac{\sum_{i=1}^{n} I(X_i \leq x_i)}{n+1} \qquad and \qquad \tilde{F}_i\,(y_i) = \frac{\sum_{i=1}^{n} I(Y_i \leq y_i)}{n+1} \quad (18)$$

The second stage involves estimating the copula dependency parameter using the pseudo likelihood function as shown in equation (19-20) and maximizing this function.

$$\ln L_{(SP)} = \sum_{i=1}^{n} \ln c\left(\tilde{F}_i\,(x_i), \tilde{F}_i\,(y_i) \quad ; \theta_{copula}\right) \quad (19)$$

$$\frac{\partial \ln L_{(SP)}}{\partial \theta_{copula}} = \frac{\partial}{\partial \theta_{copula}}\left[\sum_{i=1}^{n} \ln c\left(\tilde{F}_i\,(x_i), \tilde{F}_i\,(y_i) \quad ; \theta_{copula}\right)\right] \quad (20)$$

## Section 4: real data analysis

## 4.1 data description

The OECD data platform delivered by the OECD (organization for economic Cooperation and development) collects data on a wide range of economic, social, and environmental sectors including data about economy, health, trade, education, labor, innovation, and development. It contains many indicators describing these data. The data can be found on the following site: In this analysis, 3 indicators are used to find the dependency between and how the previously new copulas can model this dependency. The indicators are the quality of support network, water quality, and feeling safe walking alone at night. The data are collected from 41 countries as shown in Table (1) . https://stats.oecd.org/index.aspx?DataSetCode=BLI .

The indicators ( let's call it the y variable) in Table (1) and in Table (5) are distributed as Median Based Unit Rayleigh (MBUR) previously discussed by Iman Attia (Attia, 2024) with the following PDF in equation (21), CDF in equation (22), and, quantile function in equation (23):

$$f(y) = \frac{6}{\alpha^2}\left[1 - y^{\frac{1}{\theta^2}}\right] y^{\left(\frac{2}{\theta^2}-1\right)} \;,\;\; 0 < y < 1\,, \theta > 0 \quad (21)$$

$$F(y) = 3y^{\frac{2}{\theta^2}} - 2y^{\frac{3}{\theta^2}} \;,\;\; 0 < y < 1\,, \theta > 0 \quad (22)$$

$$y = \left\{-.5\left(cos\left[\frac{cos^{-1}(1-2u)}{3}\right] - \sqrt{3}\,sin\left[\frac{cos^{-1}(1-2u)}{3}\right]\right) + .5\right\}^{\theta^2} ,0 < u < 1\,, \theta > 0 \quad (23)$$

Table (1): The three indicators of the OECD data

| Country | Water quality | Feeling safe walking alone | Quality of support network |
|---|---|---|---|
| Australia | 0.92 | 0.67 | 0.93 |
| Austria | 0.92 | 0.86 | 0.92 |
| Belgium | 0.79 | 0.56 | 0.90 |
| Canada | 0.90 | 0.78 | 0.93 |
| Chile | 0.62 | 0.41 | 0.88 |
| Colombia | 0.82 | 0.50 | 0.80 |
| Costa Rica | 0.87 | 0.47 | 0.82 |
| Czechia | 0.89 | 0.77 | 0.96 |
| Denmark | 0.93 | 0.85 | 0.95 |
| Estonia | 0.86 | 0.79 | 0.95 |
| Finland | 0.97 | 0.88 | 0.96 |
| France | 0.78 | 0.74 | 0.94 |
| Germany | 0.91 | 0.76 | 0.90 |
| Greece | 0.67 | 0.69 | 0.78 |
| Hungary | 0.81 | 0.74 | 0.94 |
| Iceland | 0.97 | 0.85 | 0.98 |
| Ireland | 0.80 | 0.76 | 0.96 |
| Israel | 0.77 | 0.80 | 0.95 |
| Italy | 0.77 | 0.73 | 0.89 |
| Japan | 0.87 | 0.77 | 0.89 |
| Korea | 0.82 | 0.82 | 0.80 |
| Latvia | 0.83 | 0.72 | 0.92 |
| Lithuania | 0.83 | 0.62 | 0.89 |
| Luxembourg | 0.85 | 0.87 | 0.91 |
| Mexico | 0.75 | 0.42 | 0.77 |
| Netherlands | 0.91 | 0.83 | 0.94 |
| New Zealand | 0.85 | 0.66 | 0.95 |
| Norway | 0.98 | 0.93 | 0.96 |
| Poland | 0.82 | 0.71 | 0.94 |
| Portugal | 0.89 | 0.83 | 0.87 |
| Slovak Republic | 0.81 | 0.76 | 0.95 |
| Slovenia | 0.93 | 0.91 | 0.95 |
| Spain | 0.76 | 0.80 | 0.93 |
| Sweden | 0.97 | 0.79 | 0.94 |
| Switzerland | 0.96 | 0.86 | 0.94 |
| Turkey | 0.62 | 0.59 | 0.85 |
| United kingdom | 0.82 | 0.78 | 0.93 |
| United states | 0.88 | 0.78 | 0.94 |
| Brazil | 0.70 | 0.45 | 0.83 |
| Russia | 0.62 | 0.64 | 0.89 |
| South Africa | 0.72 | 0.40 | 0.89 |

The 'water quality' indicator measures the self-reported satisfaction with water quality as regard drinking clean water and extent of water pollution. It is expressed as the percentage of people who report satisfaction. The indicator of 'feeling safe walking alone' measures the percentage of people who report feeling safe when walking alone in their local area after dark. The 'quality of support network' indicator measures the accessibility and trustworthiness of social support for individuals and it is explicated as the percentage of population with social support who described having someone they can depend on for support when they need it. Table (2) shows the descriptive statistics of the 3 indicators. Table (3) shows the empirical Kendall tau coefficient and its associated p-value. Figure (7) shows the boxplot for each indicator. Figure (8-10) show the scatter plot of the dependent indicators. Each indicator fits the Median Based Unit Rayleigh Distribution (MBUR). Table (4) shows the statistical validity indices for each indicator. Figure (11) shows the histogram and the fitted MBUR curve for each indicator.

Table (2): descriptive statistics for the 3 indicators:

| indicator | min | mean | Standard deviation | skewness | kurtosis | 25percentile | 50percentile | 75percentile | max |
|---|---|---|---|---|---|---|---|---|---|
| Water quality | 0.62 | 0.8332 | 0.0972 | -0.6059 | 2.9144 | 0.7775 | 0.83 | 0.91 | 0.98 |
| Feeling safe walking alone | 0.4 | 0.7207 | 0.143 | -0.9486 | 3.0353 | 0.655 | 0.76 | 0.8225 | 0.93 |
| Quality of support network | 0.77 | 0.9078 | 0.0538 | -1.176 | 3.5406 | 0.89 | 0.93 | 0.95 | 0.98 |

Table (3): empirical Kendall tau coefficient of the three indicators:

| | Water quality | Feeling safe walking alone | Quality of support network |
|---|---|---|---|
| Water quality | 1 | 0.5206 (0.000) | 0.3929 (0.0006) |
| Feeling safe walking alone | 0.5206 (0.000) | 1 | 0.4344 (0.0001) |
| Quality of support network | 0.3929 (0.0006) | 0.4344 (0.0001) | 1 |

Table (4): statistical validity indices of the distributional fit (MBUR) of the three indicators:

| indicator | Estimated theta | variance | AIC | CAIC | BIC | HQIC | KS-test | Ho | p-value of KS-test |
|---|---|---|---|---|---|---|---|---|---|
| Water quality | 0.4776 | 0.00072157 | -78.9952 | -78.8926 | -77.2817 | -78.3712 | 0.0991 | Fail to reject | 0.7789 |
| Feeling safe walking alone | 0.6494 | 0.0013 | -47.1347 | -47.0321 | -45.4211 | -46.5107 | 0.1206 | Fail to reject | 0.5492 |
| Quality pf support network | 0.3444 | 0.00037494 | -131.6505 | -131.6505 | -131.5479 | -131.0265 | 0.1806 | Fail to reject | 0.1217 |

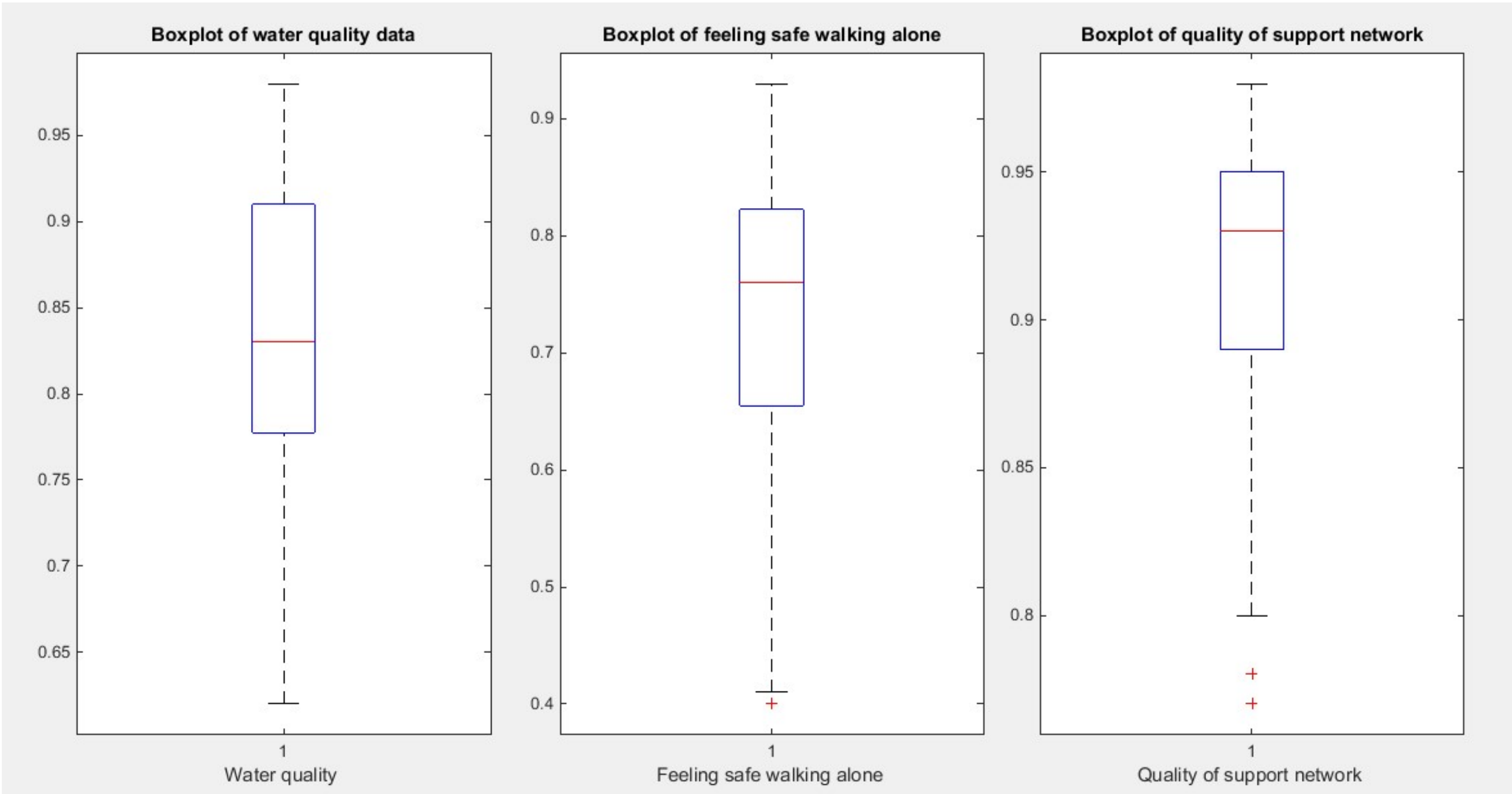

Fig. 7 shows the Boxplot for each indicator. The first two indicators show similar pattern of skewness (left sides skewness) while the third indicator exhibits more left sided skewness than the previous two indicators.

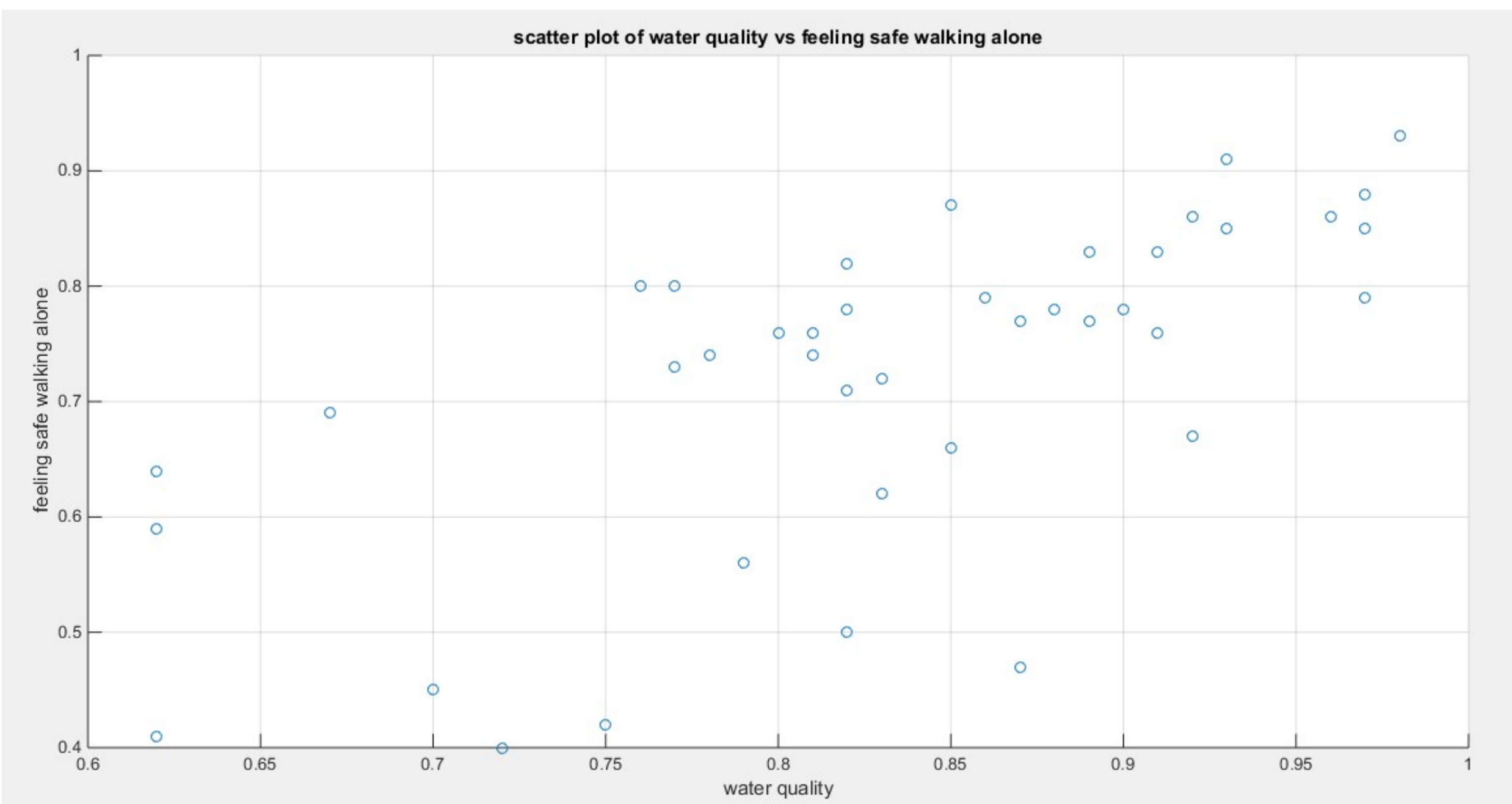

Fig. 8 shows the scatter plot of the water quality variable vs the variable describing feeling safe walking alone. The data are mainly concentrated on the upper left side of the graph most probably reflecting upper tail dependency.

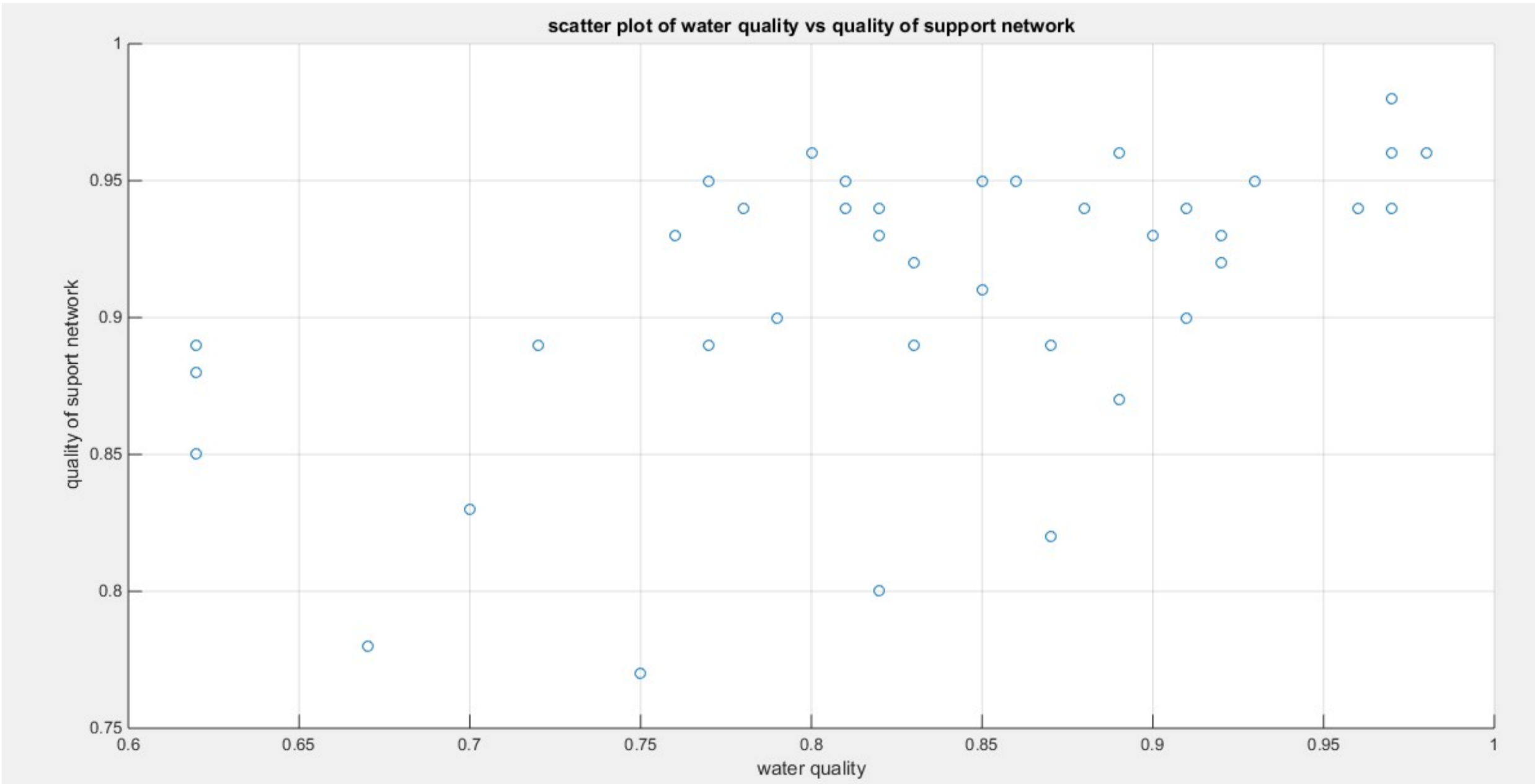

Fig. 9 shows the scatter plot of the variable water quality vs the variable quality of support network, also the data are mainly concentrated in the upper left corner of the scatter presumably reflecting upper tail dependency.

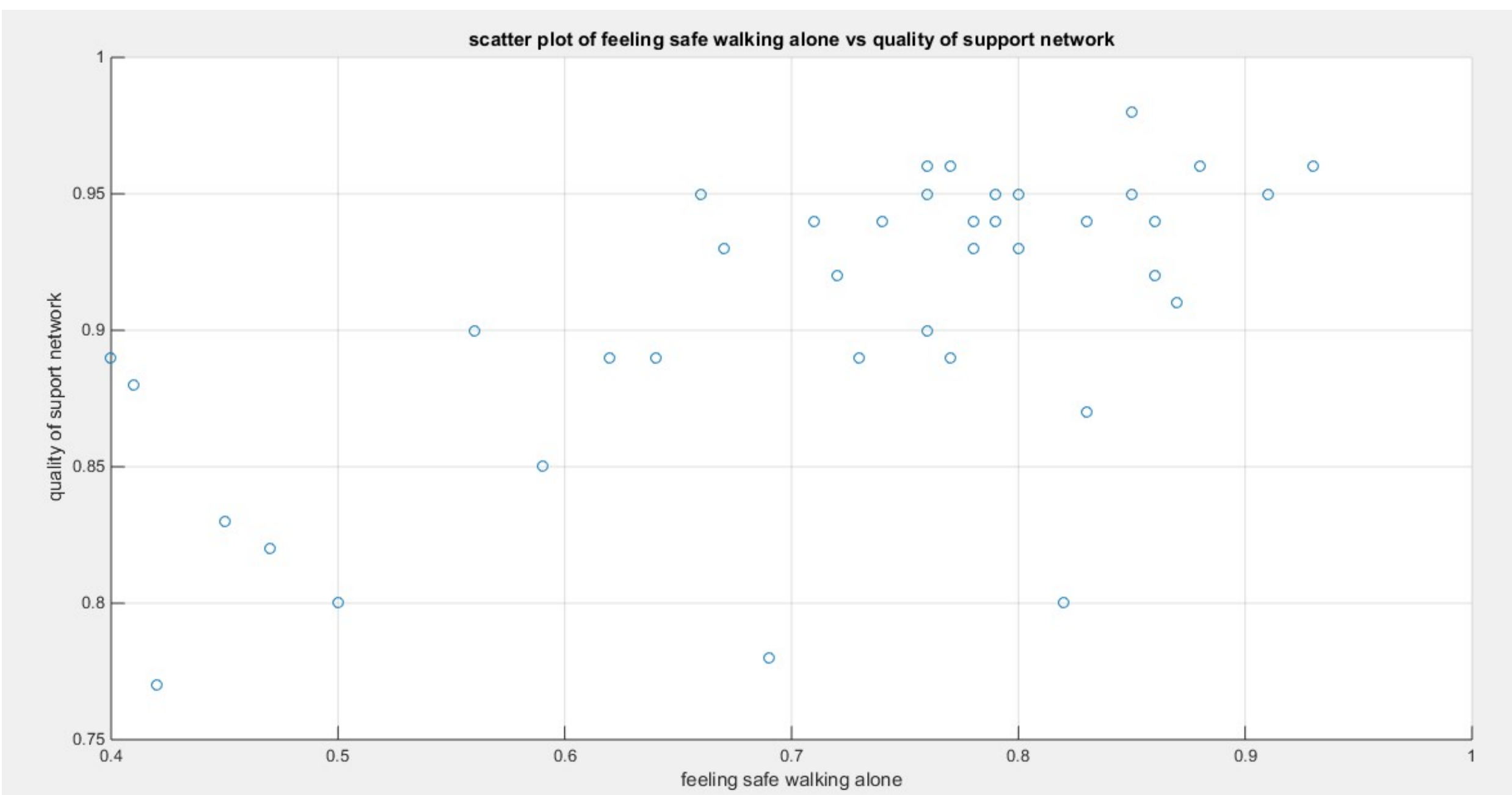

Fig. 10 shows the scatter plot of the variable feeling safe walking alone vs the variable of the quality of support network. The data are also concentrated on the right upper corner of the graph likely denoting the presence of upper tail dependency.

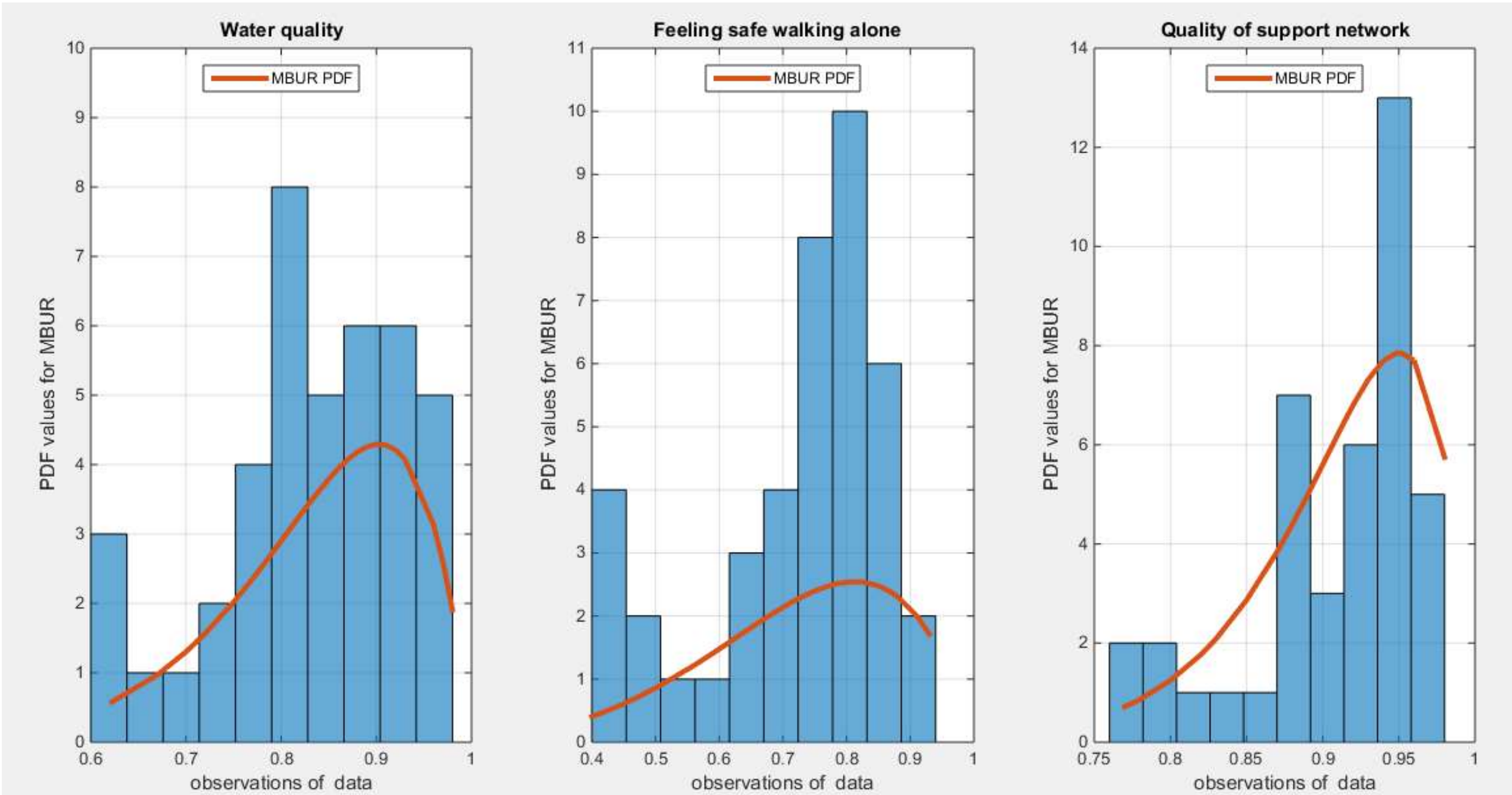

Fig. 11 shows the histogram of the three indicators which exhibit left sided skewness. The graph also shows the fitted MBUR curve for each indicator. The three indicators more or less have near similar kurtosis and skewness. Copula that can fit symmetrical dependency can model dependency between such variables.

### 4.2. procedure of analysis (IFM)

After running the inference function for margins (IFM) and estimating the marginal parameter for each pair of variables then fitting each copula to estimate its parameter. Goodness of Fit tests are conducted to assess the dependence relationship for each pair of the variables.

The steps of analysis are as follows:

1- Estimate the marginal parameters for each variable.
2- Use the IFM procedure to estimate the dependency parameter of the proposed copula model.
3- Obtain the theoretical tau from the specific relation between dependency parameter and the Kendall tau dor each copula.
4- Obtain the Cramer Von Mises test from this estimation process and call it $CVM_{data}$ that will be compared to $CVM_{samples}$
5- Run sampling technique like Metropolis Hasting (MH) procedure to test the null hypothesis using the estimated dependency parameter and both the estimated marginal parameters.
6- Construct sampling distribution for the estimated-$theta_{samples}$ , theoretical-$tau_{samples}$ , and $CVM_{samples}$ .

For each of the sample generated by MH using the estimated dependency parameter and the estimated marginal parameters, they will be transformed back into variables using the specific quantiles. These variables will be subjected to the IFM procedure; first estimate the marginal parameters then using the estimated marginal CDF for each variable to estimate the dependency parameter, in other words, repeat steps (1 to 4) for each sample so you can be able to construct the sampling distribution, hence, constructing the confidence interval for each of the statistical indices. $theta_{data}$ denotes theta estimated from the data.

The null hypothesis for the dependency parameter is: $H_0$: $theta_{samples} = theta_{data}$
The alternative hypothesis for the dependency parameter is: $H_1$: $theta_{samples} \neq theta_{data}$

The null hypothesis for the Kendall tau coefficient is: $H_0$: $tau_{samples} = tau_{data}$
The alternative hypothesis for the Kendall tau coefficient is: $H_1$: $tau_{samples} \neq tau_{data}$

The null hypothesis for the CVM is: $H_0$: $CVM_{samples} = CVM_{data}$
The alternative hypothesis for the CVM is: $H_1$: $CVM_{samples} \neq CVM_{data}$

**First Copula:**

***Model the dependency between the 'water quality' and 'feeling safe walking alone'***: The empirical tau equals 0.5206 while the theoretical tau is 0.5202 with difference value of 3.6968e-04. The estimated dependency parameter, theta $\hat{\theta}$ , is 0.4798 with estimated variance of 8.3533e-06. The Cramer Von Mises test (CVM) is 0.1388. The null hypothesis for dependency was tested by conducting resampling using the Metropolis Hasting algorithm of MCMC procedure. The proposed null hypothesis was that the dependency parameter equals the estimated theta against the alternative hypothesis that the population dependency parameter does not equal the estimated theta. Also the null hypothesis for the population Kendall tau coefficient being equal to the theoretical tau against the alternative hypothesis of not being equal to it was tested. The null hypothesis for the CVM proposing that the sampling CVM equals observed CVM obtained from the estimation procedure was investigated. The sampling distribution of each of the previously mentioned indices, thetas, tau and CVM is shown in Figures (12-14). For each figure, the descriptive statistics indices are shown. From the figures the indices are placed within the acceptance zone between the 2.5$^{th}$ and 97.5$^{th}$ quantiles. So the null hypotheses fail to be rejected. The copula fits the data well. And the dependence parameter models, within the context of the copula, the relation between the two variables. The confidence interval (CI) for the estimated dependence parameter is (0.1822, 0.8084), for the theoretical tau; it is (0.1916, 0.8178), and for CVM; it is (0.0316, 0.2937).

According to the sandwich variance the variance covariance matrix for this model is:

$$var - cov = \begin{bmatrix} 0.0081 & 0.0144 & 0.00002 \\ 0.0144 & 0.0257 & 0.00003 \\ 0.00002 & 0.00003 & 2.0301e-07 \end{bmatrix}$$

The validity indices for the second stage are:

$AIC = -21.4395$ , $CAIC = -21.337$ $BIC = -19.726$ $HQIC = -20.8156$

The validity indices for the whole model is the summation of these indices for marginal (first stage) and copula (second stage):

$AIC = -147.5694$ , $CAIC = -147.2617$ $BIC = -142.4287$ $HQIC = -145.6974$

Theoretical upper tail coefficient is 0.6054. At quantile 0.85; the empirical upper tail coefficient in the direction of the 'water quality' variable is 0.2 (one point showing joint match out of 5 points in the upper tail at this threshold, 0.85) and in the direction of the 'feeling safe walk at night' variable is 0.5 (one point showing joint match out of 2 points in the upper tail at this threshold, 0.85) while the confidence interval of the empirical distribution obtained from bootstrap sampling under the null hypothesis is [0,1]. This is expected because of the small sample size and left skewness of data.

Figure 15 shows the P-P plot of the empirical vs the theoretical copula. Figure 16-18 show the surface of the copula density, the contour plot of the copula density, and the contour plot of the copula distribution at the estimated parameter.

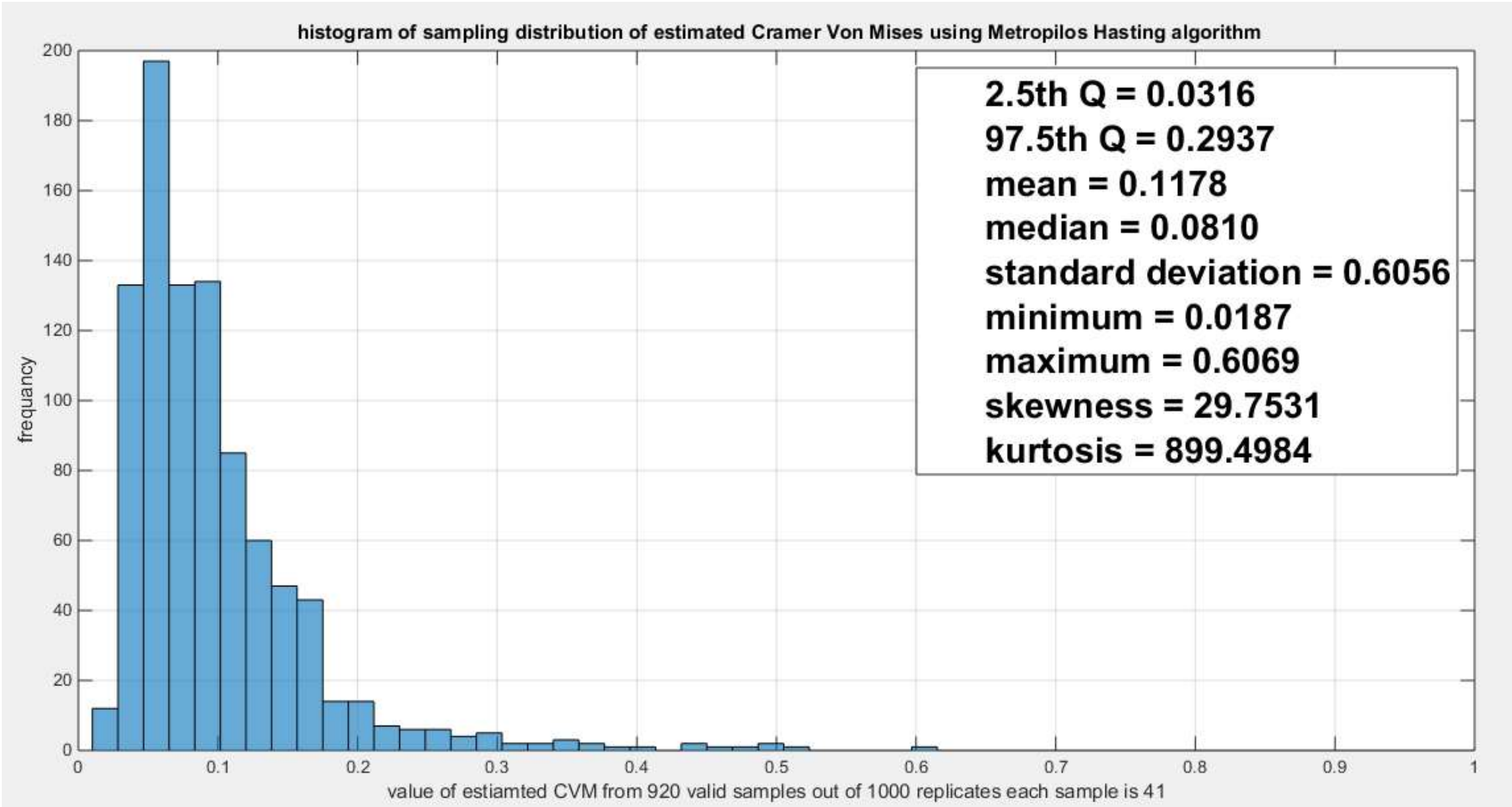

Fig. 12 shows the sampling distribution of the CVM. The CVM data is 0.1388 and it lays between the 2.5th and the 97.5th quantiles so null hypothesis is failed to be rejected and the copula fits the data well. The p-value (probability of values less than or equal to 0.1388) is 0.8207, in other words, this 0.1388 does not lay in either of the tail region of rejection.

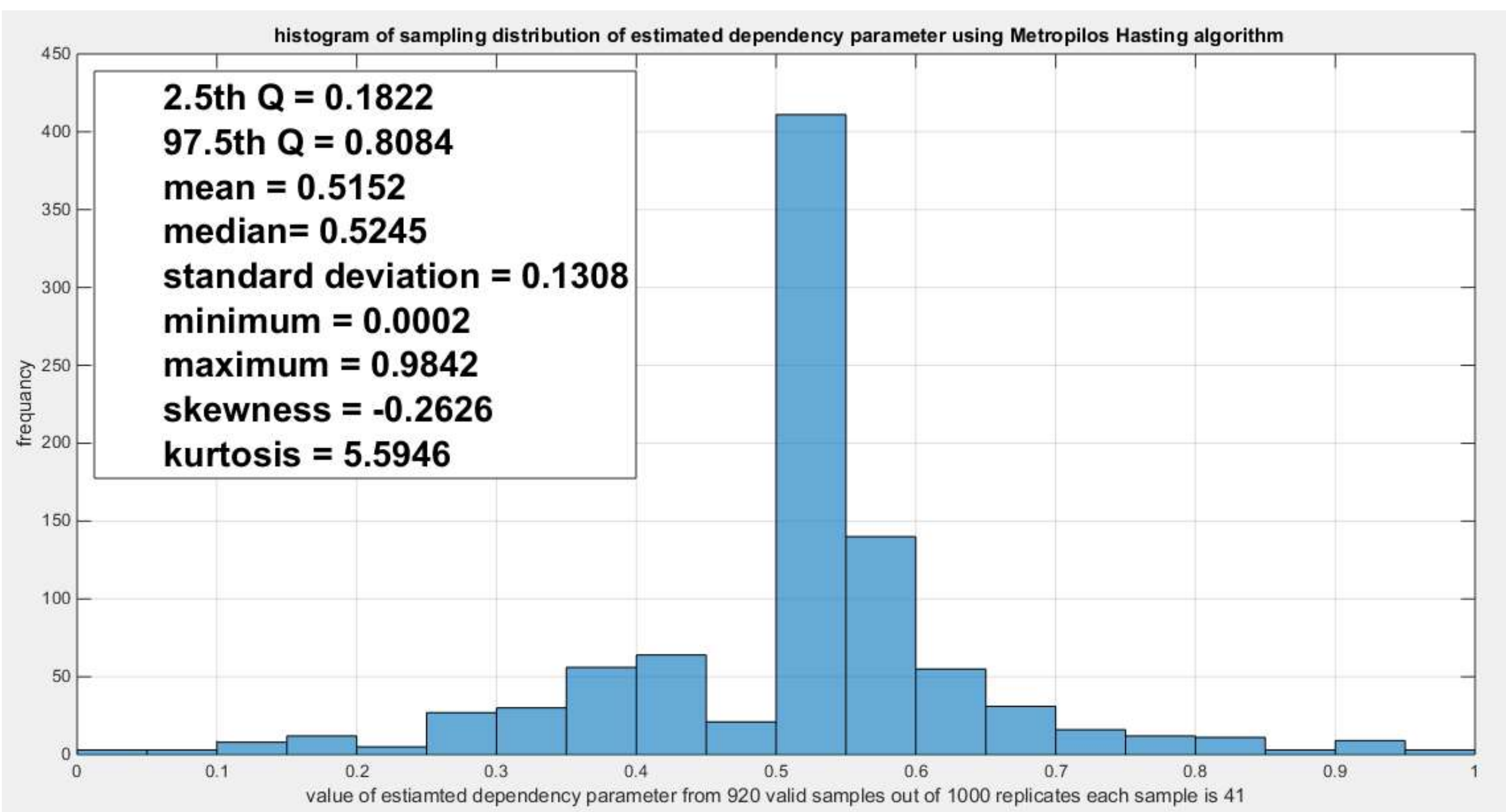

Fig. 13 shows the sampling distribution of the dependency parameter. The 0.4798 lays in the acceptance region between the two quantiles, the 2.5th and the 97.5th quantiles. So the null hypothesis fails to reject its assumption. The p-value (probability of the values being less than or equal to this 0.4798) is 0.2304 which indicates that 0.4798 is not placed in either zones of rejection.

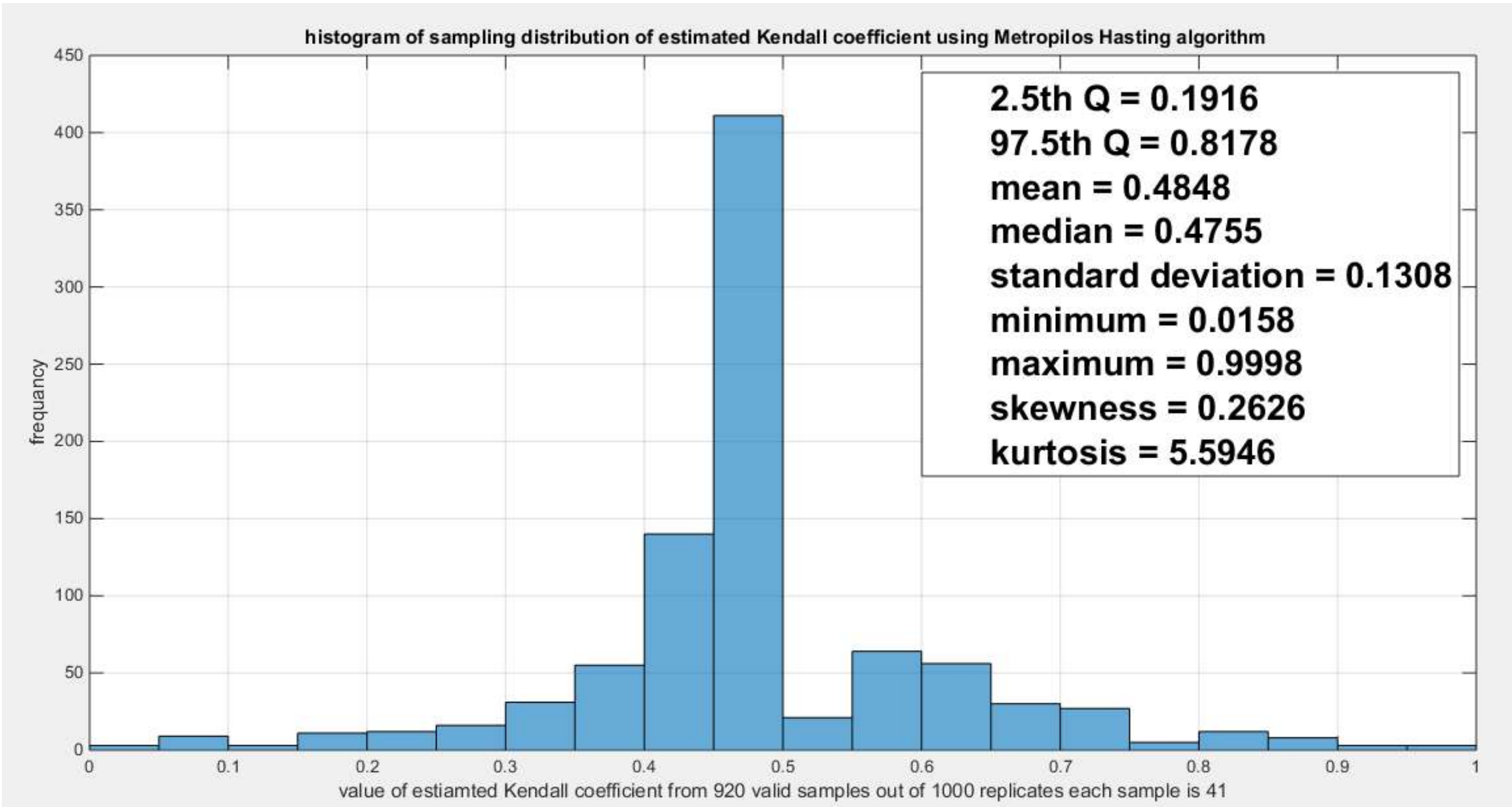

Fig. 14 shows the sampling distribution of the Kendal tau coefficient. Its value 0.5202 is between the two quantiles 2.5th and 97.5th which is the acceptance zones. So the null hypothesis fails to reject its assumption. The p-values (probability of the values less than or equal to 0.5202) is 0.7696 which denotes that this value is far away from both rejection zones.

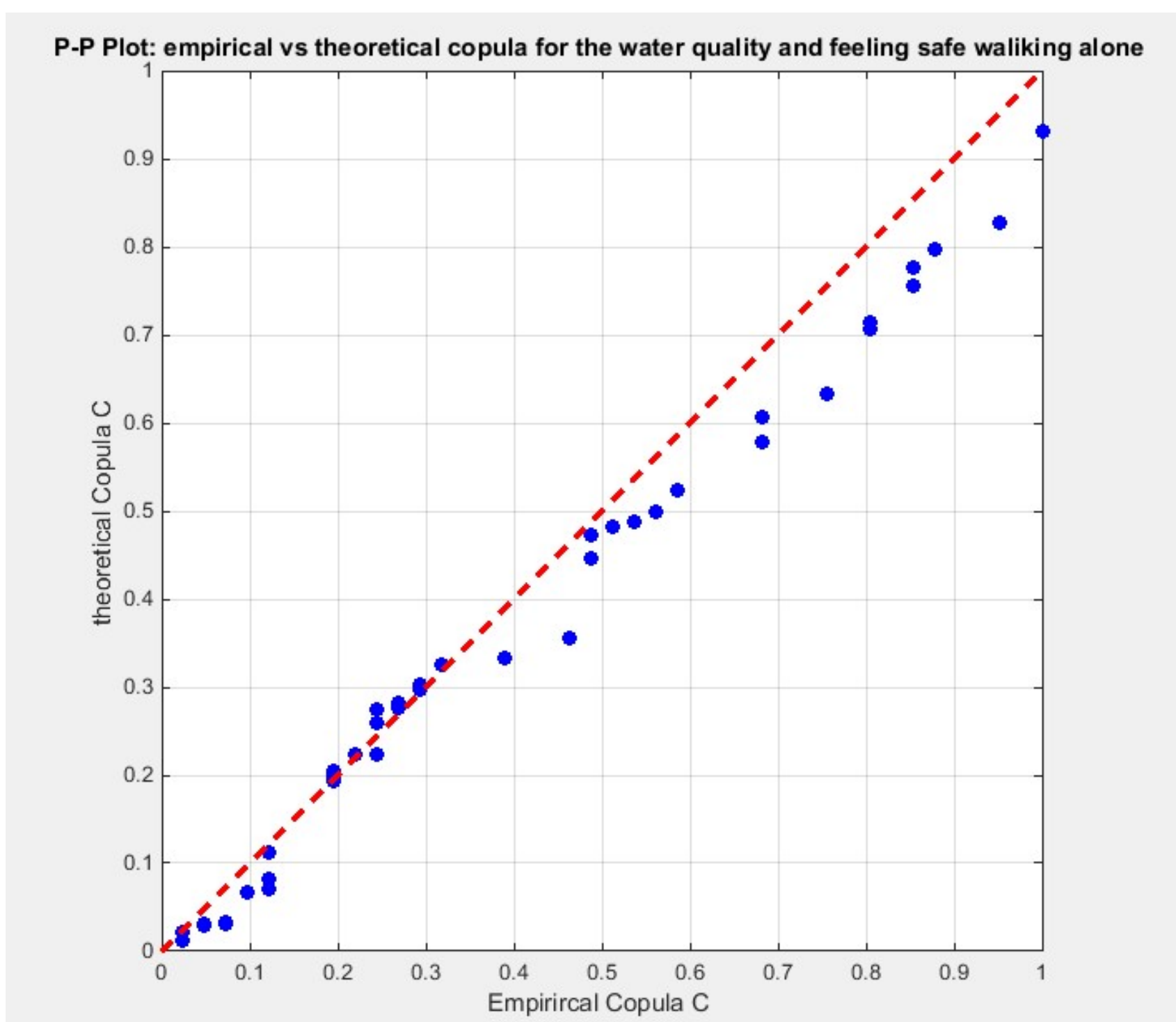

Fig. 15 shows the P-P plot of the empirical vs the theoretical copula. The theoretical copula shows near perfect alignment of its lower part with the diagonal line than its upper part

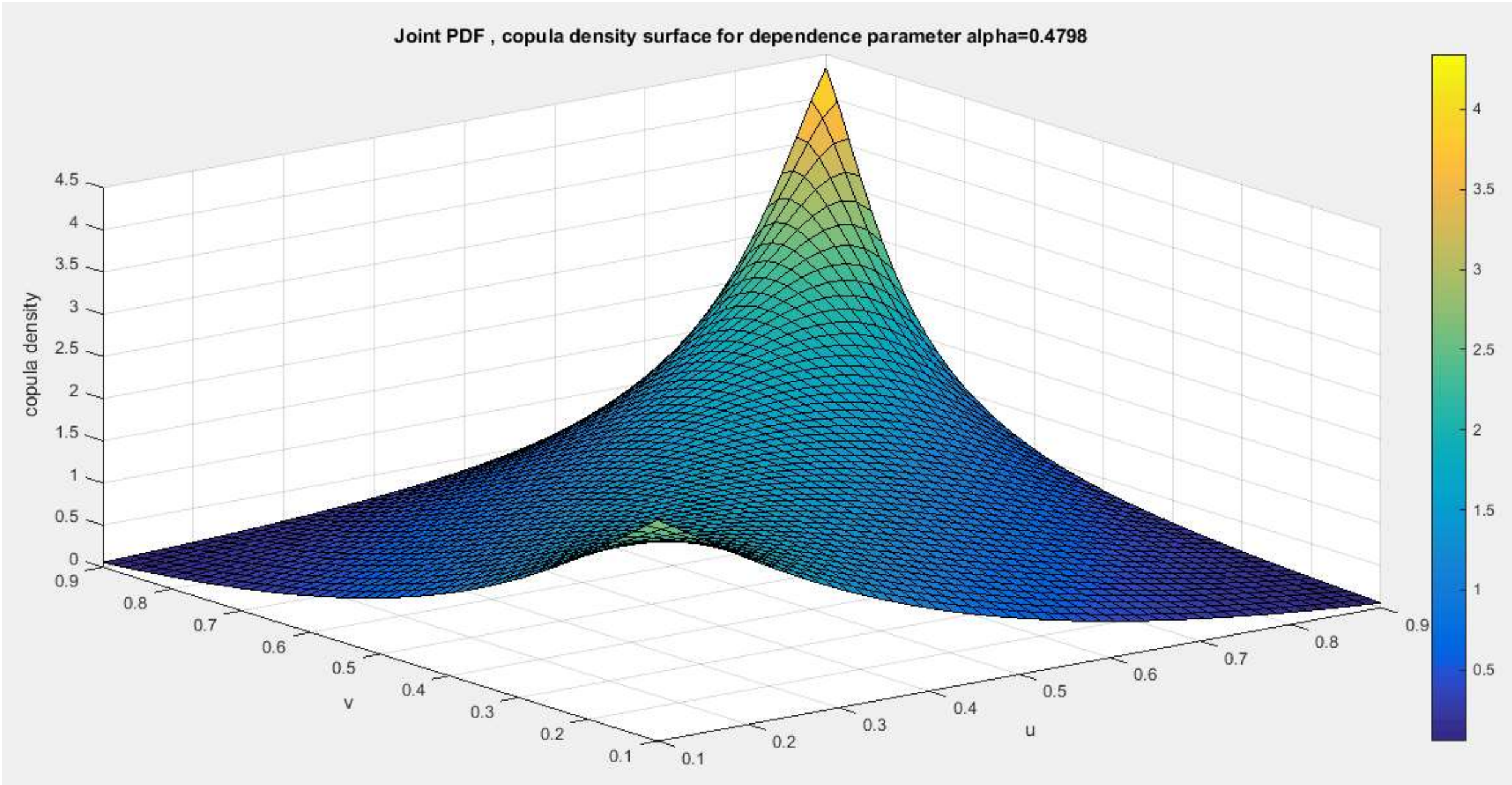

Fig.16 shows the PDF surface of copula density with parameter 0.4798

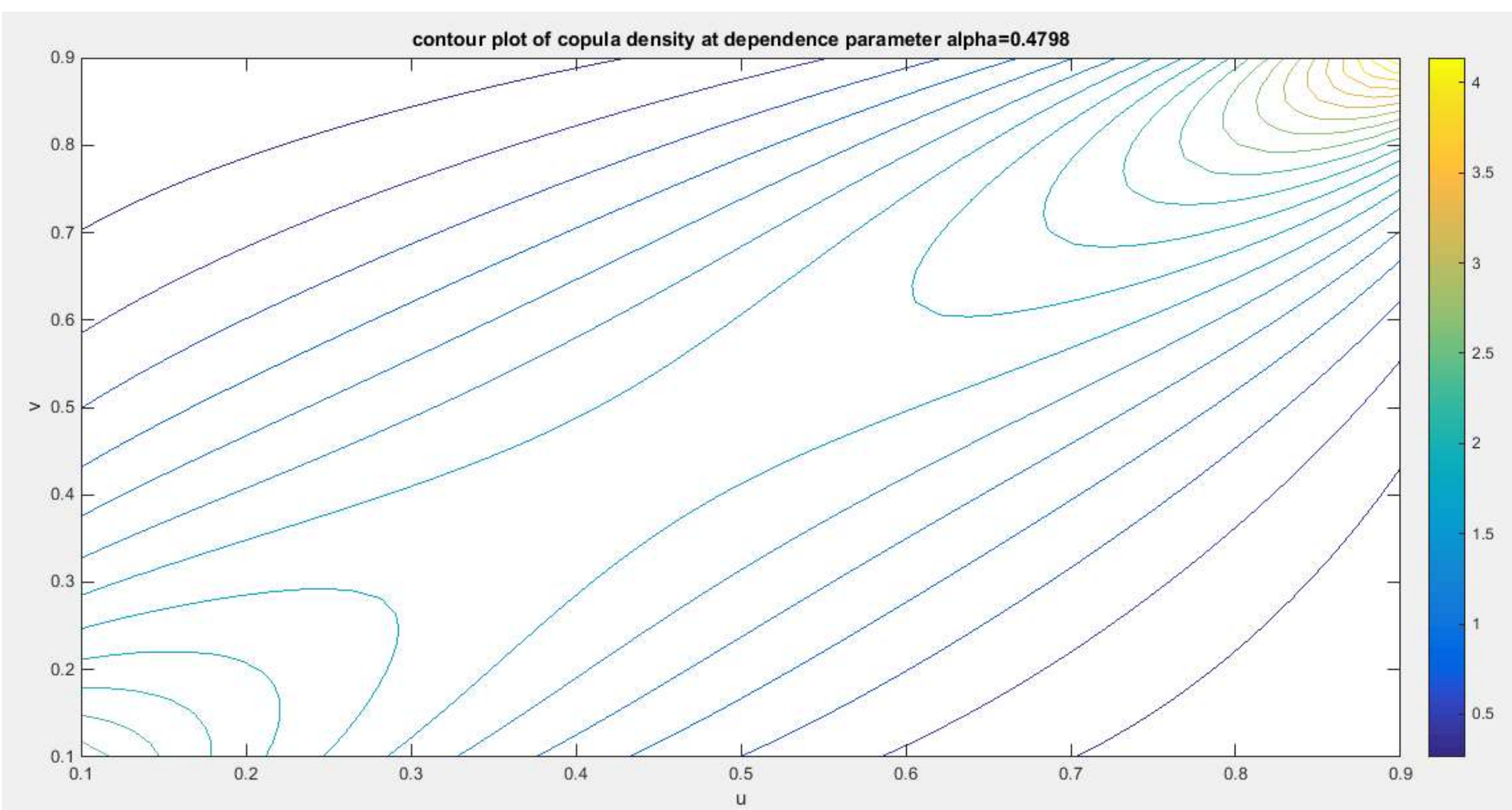

Fig.17 shows the contour plot of the copula density with parameter 0.4798

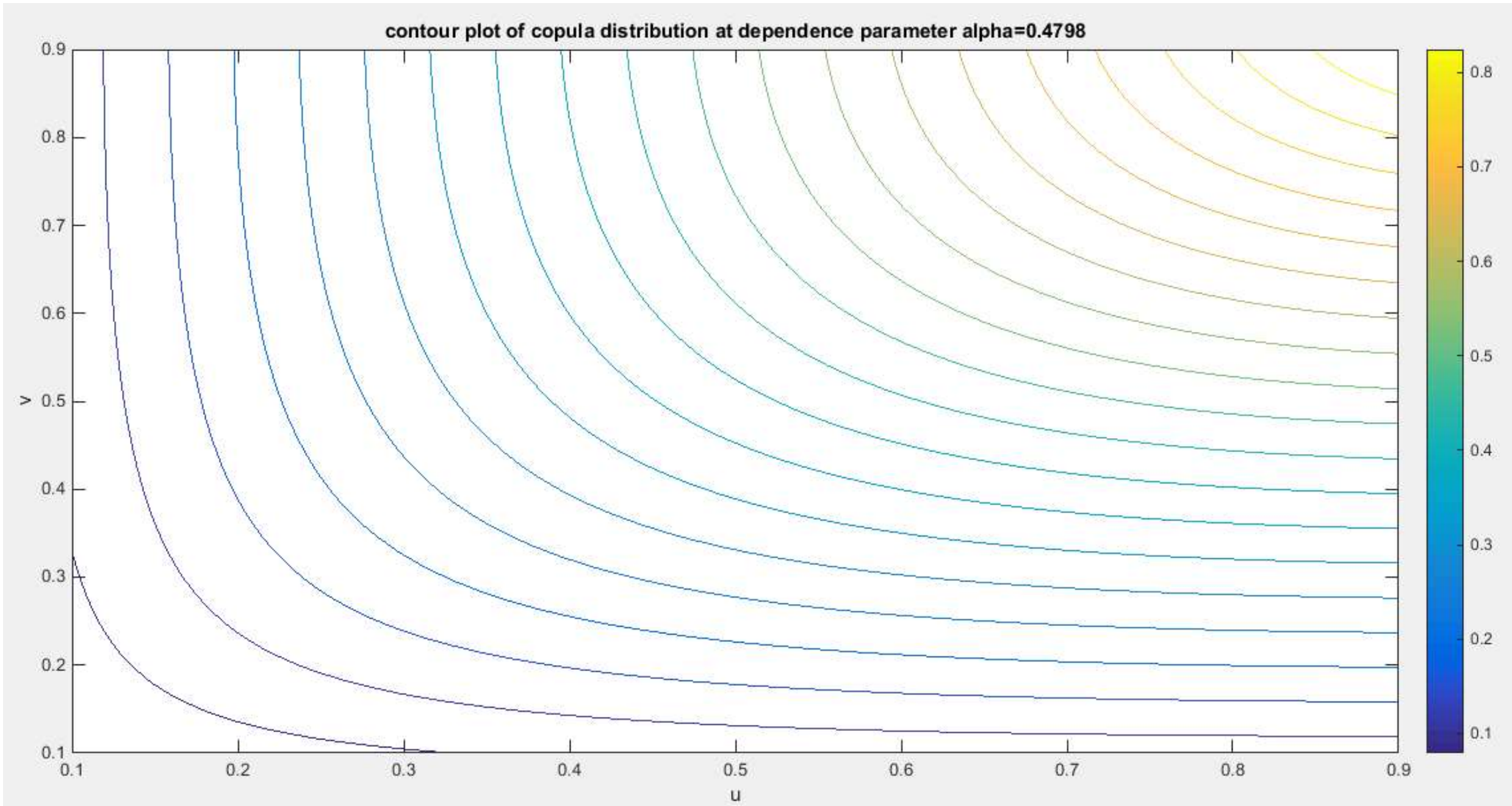

Fig. 18 shows the contour plot of the copula distribution with dependency parameter 0.4798

## Second copula

***Model the dependency between the 'water quality' and 'quality of support network'***: The empirical tau is 0.3929 and the theoretical tau is 0.3926 with dissimilarity value of 3.8927e-04 The estimated dependency parameter, theta $\hat{\theta}$ , is 0.7794 and estimated variance of 0.0018. The Cramer Von Mises test (CVM) is 0.2142. to test the null hypothesis for dependency, resampling using the Metropolis Hasting algorithm of MCMC procedure was performed. The propositional null hypothesis was stated that the population dependency parameter equals the estimated theta against the alternative hypothesis that the population dependency parameter does not equal the estimated theta. Moreover, the null hypothesis for the equality of the population Kendall tau with the theoretical tau is tested against the alternative hypothesis of not being equal to each other. The null hypothesis for the CVM suggesting that the sampling CVM equals the observed CVM attained from the estimation procedure was also investigated. These are all two sided tests. The sampling distribution of each of the previously mentioned indices, thetas, tau and CVM is shown in Figures (19-21). The descriptive statistical indices are shown in each figure. Analysis of the figures denoted that the indices are located within the acceptance region between the $2.5^{th}$ and $97.5^{th}$ quantiles. So the null hypotheses fail to be rejected. The copula can be said to pass some of the goodness of fit (GOF) tests for this data. And the dependence parameter mockups, within the milieu of this copula, the connection between the two variables. The confidence interval (CI) for the estimated dependence parameter is (0.2571, 0.9530), for the theoretical tau; it is (0.0918, 0.9078), and for CVM; it is (0.0485, 0.3715).

According to the sandwich variance the variance covariance matrix for this model is:

$$var - cov = \begin{bmatrix} 0.0269 & 0.0143 & -0.0004 \\ 0.0143 & 0.0076 & -0.0002 \\ -0.0004 & -0.0002 & 0.00005 \end{bmatrix}$$

The validity indices for the second stage are:

$$AIC = -11.1291 \quad , \quad CAIC = -11.0265 \quad BIC = -9.4155 \qquad HQIC = -10.5051$$

The validity indices for the whole model is the summation of these indices for marginal (first stage) and copula (second stage):

$AIC = -221.7747$ , $CAIC = -221.467$ $BIC = -216.634$ $HQIC = -219.9028$

Theoretical upper tail coefficient is 0.4764. At quantile 0.85; the empirical upper tail coefficient in the direction of the 'water quality' variable is 0.2 (one point showing joint match out of 5 points in the upper tail at this threshold, 0.85) and in the direction of the 'Quality of support network' variable is 1 (one point showing joint match out of 1 point in the upper tail at this threshold, 0.85) while the confidence interval of the empirical distribution obtained from bootstrap sampling under the null hypothesis is [0,1]. This is expected because of the small sample size and left skewness of data.

Figure 22 shows the P-P plot of the empirical vs the theoretical copula. Figure 23-25 show the contour plot of the copula density, the contour plot of the copula distribution, and the surface of the copula density at the estimated parameter.

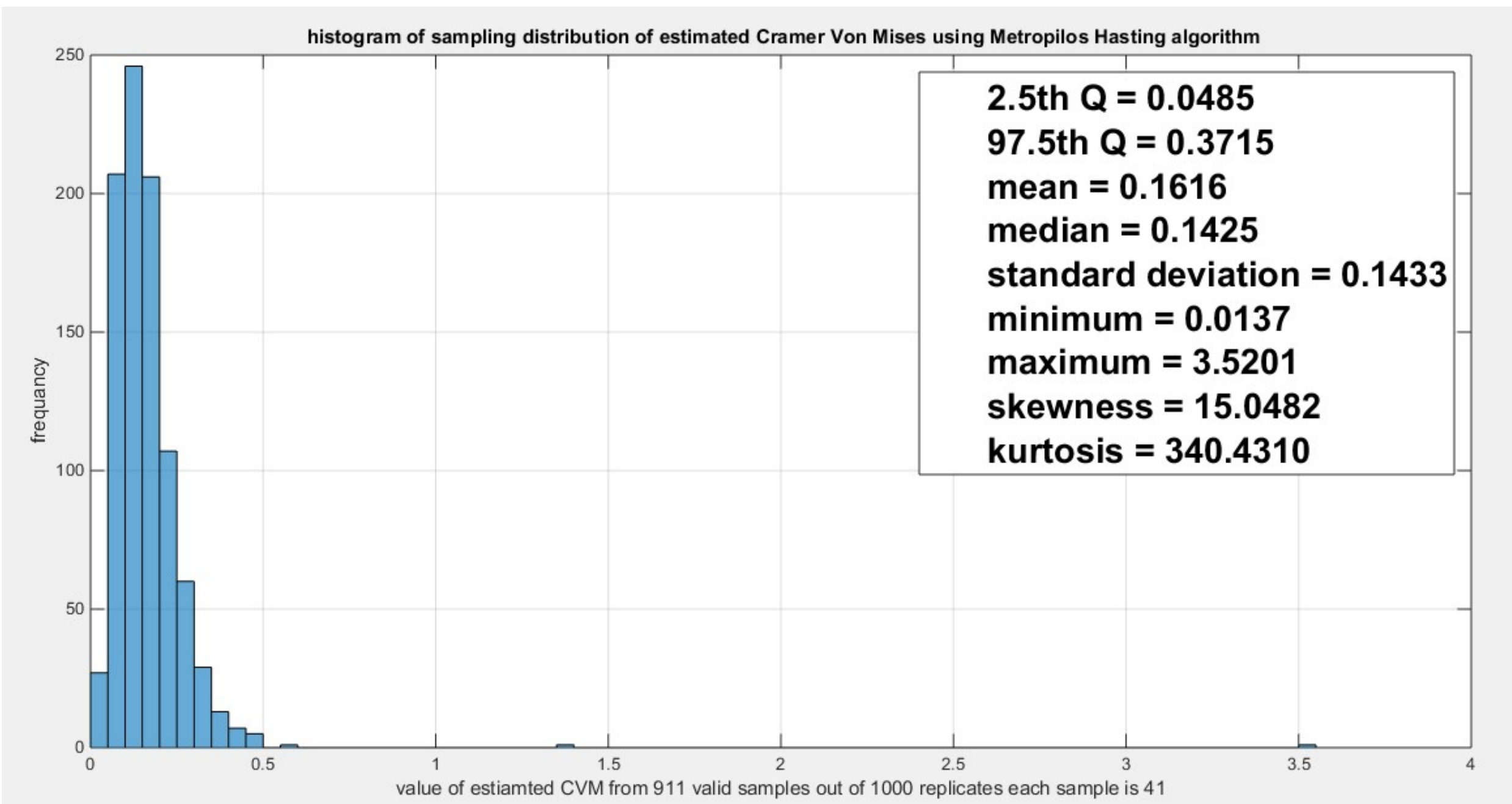

Fig. 19 shows the sampling distribution of the CVM. The data obtained CVM is 0.2142 and it lies between the 2.5th and the 97.5th quantiles so the null hypothesis is failed to be rejected and the copula fits the data well. The p-value (probability of values less than or equal to 0.2142) is 0.7846, in other words, this 0.2142 does not lay in either of the tail region of rejection.

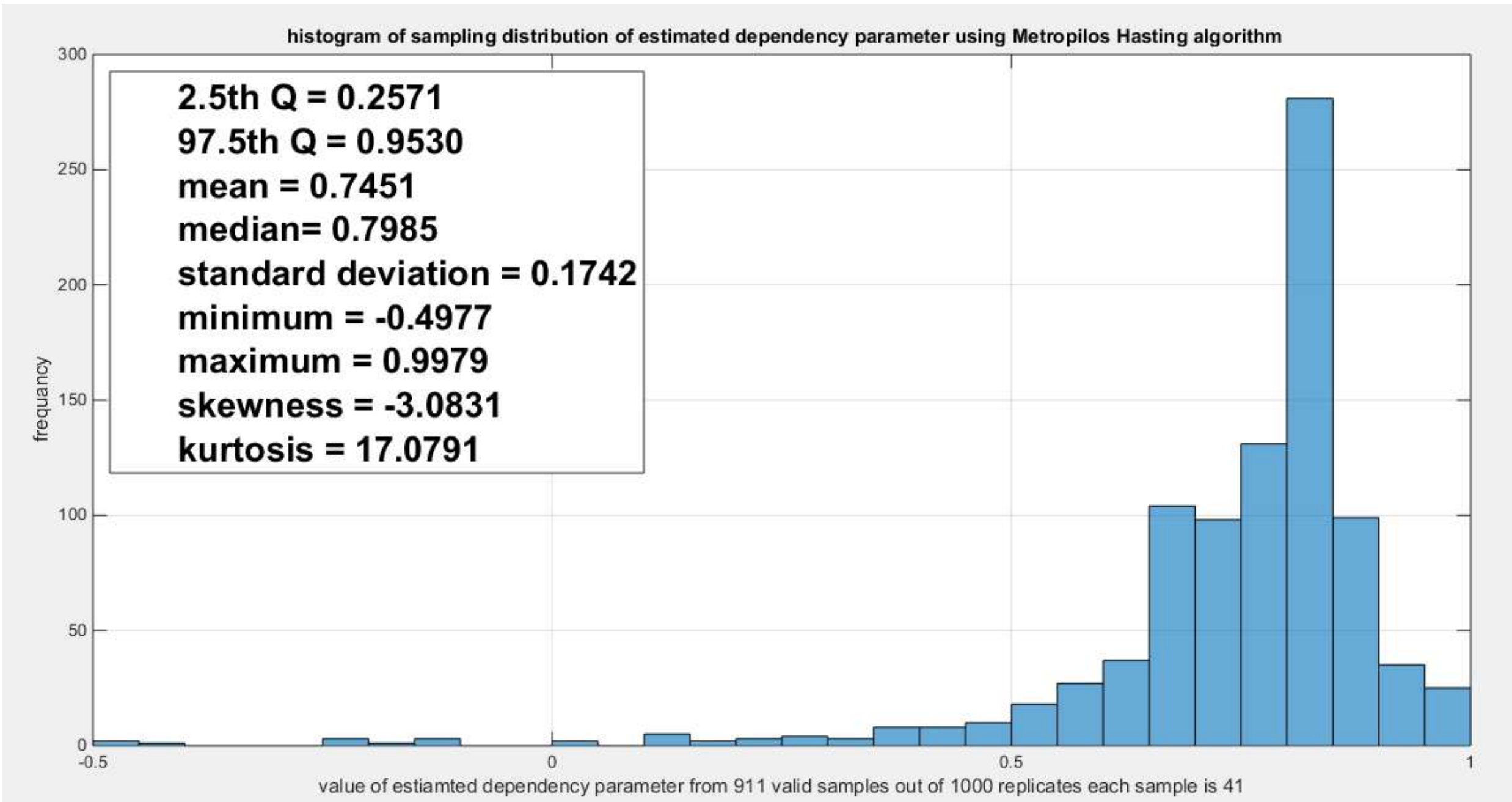

Fig. 20 shows the sampling distribution of the dependency parameter. The 0.7794 lies in the acceptance region between the two quantiles, the 2.5th and the 97.5th quantiles. So the null hypothesis fails to reject its assumption. The p-value (probability of the values being less than or equal to this 0.7794) is 0.3725which indicates that 0.7794 is not placed in either zones of rejection.

.

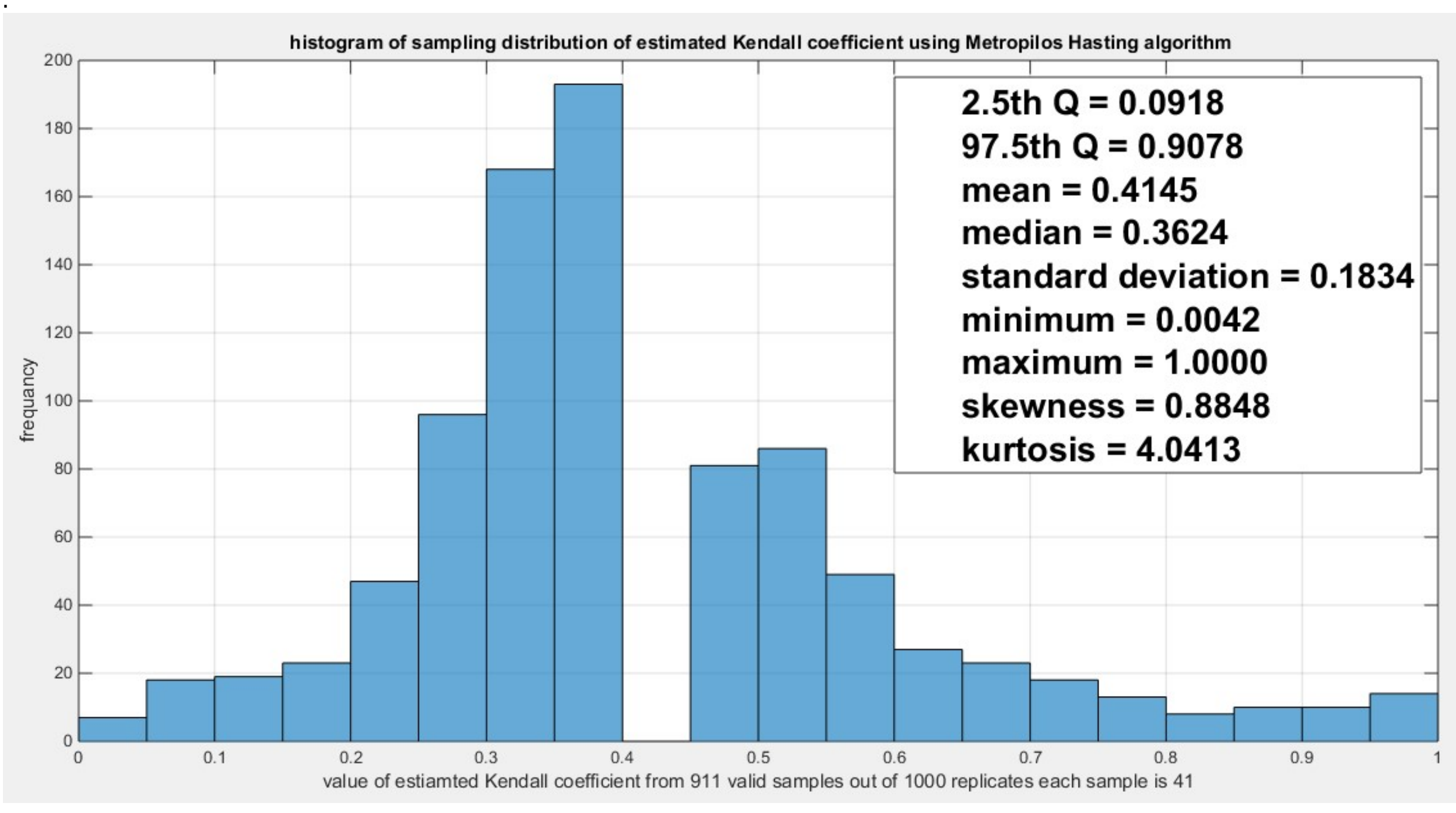

Fig. 21 shows the sampling distribution of the Kendal tau coefficient. Its value 0.3926 is between the two quantiles 2.5th and 97.5th which is the acceptance zones. So the null hypothesis fails to reject its assumption. The p-values (probability of the values less than or equal to 0.3926) is 0.6275 which denotes that this value 0.3926 is far away from both rejection zones on the right and the left.

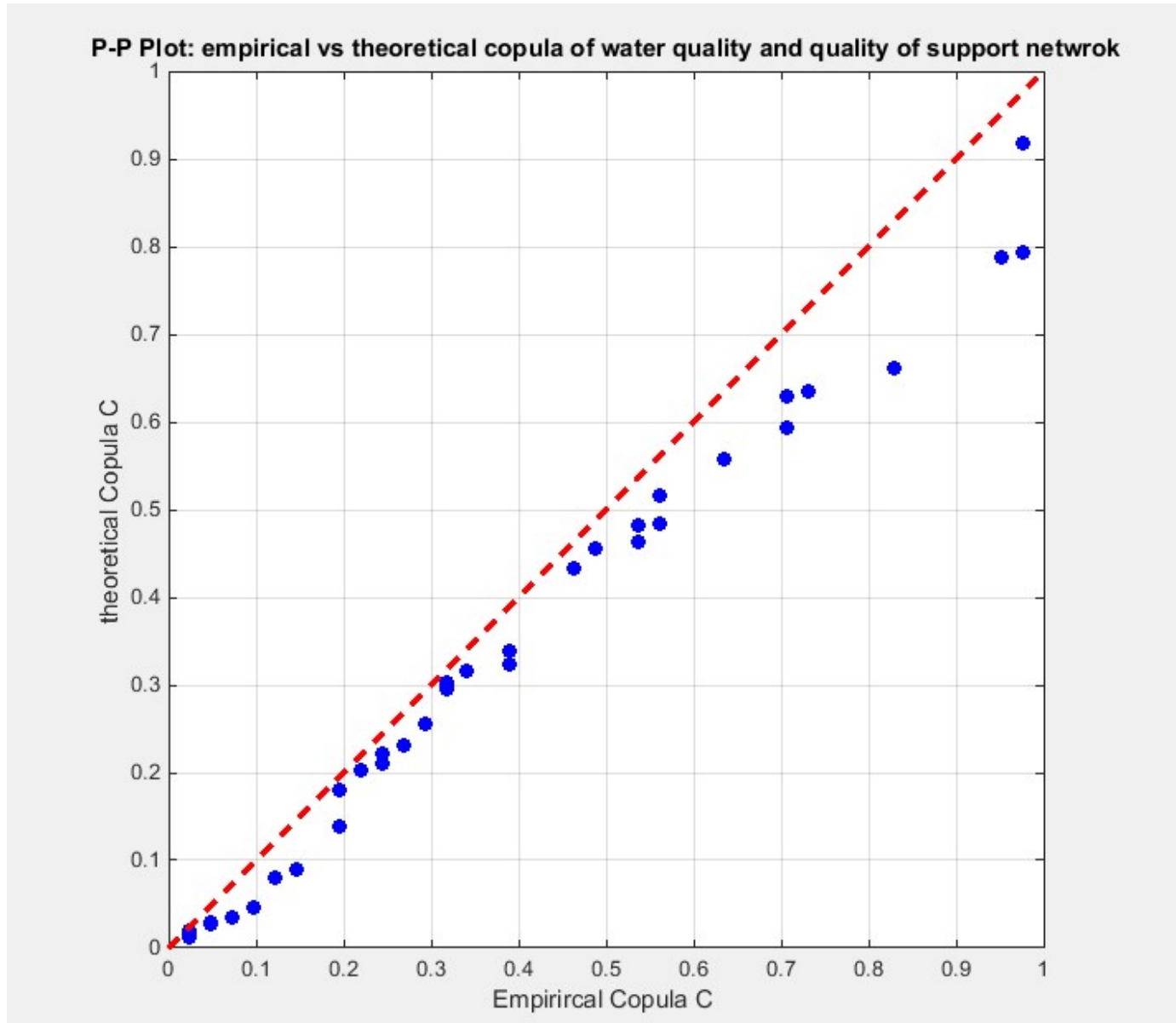

Fig. 22 shows the P-P plot of the empirical vs the theoretical copula. The theoretical copula shows near perfect alignment of its lower part and its center with the diagonal line than its upper part.

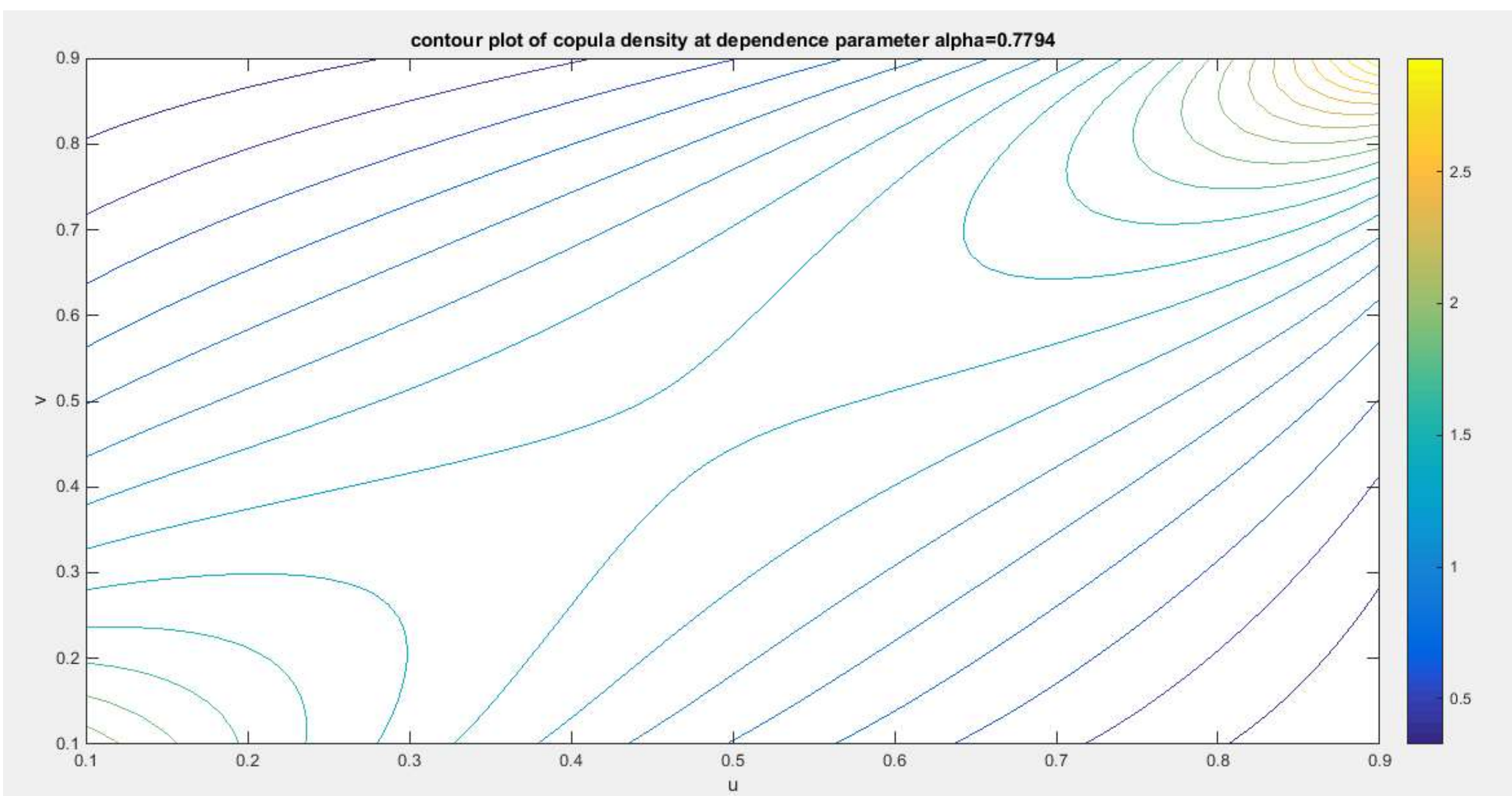

Fig. 23 shows the contour plot of the copula density with parameter 0.7794

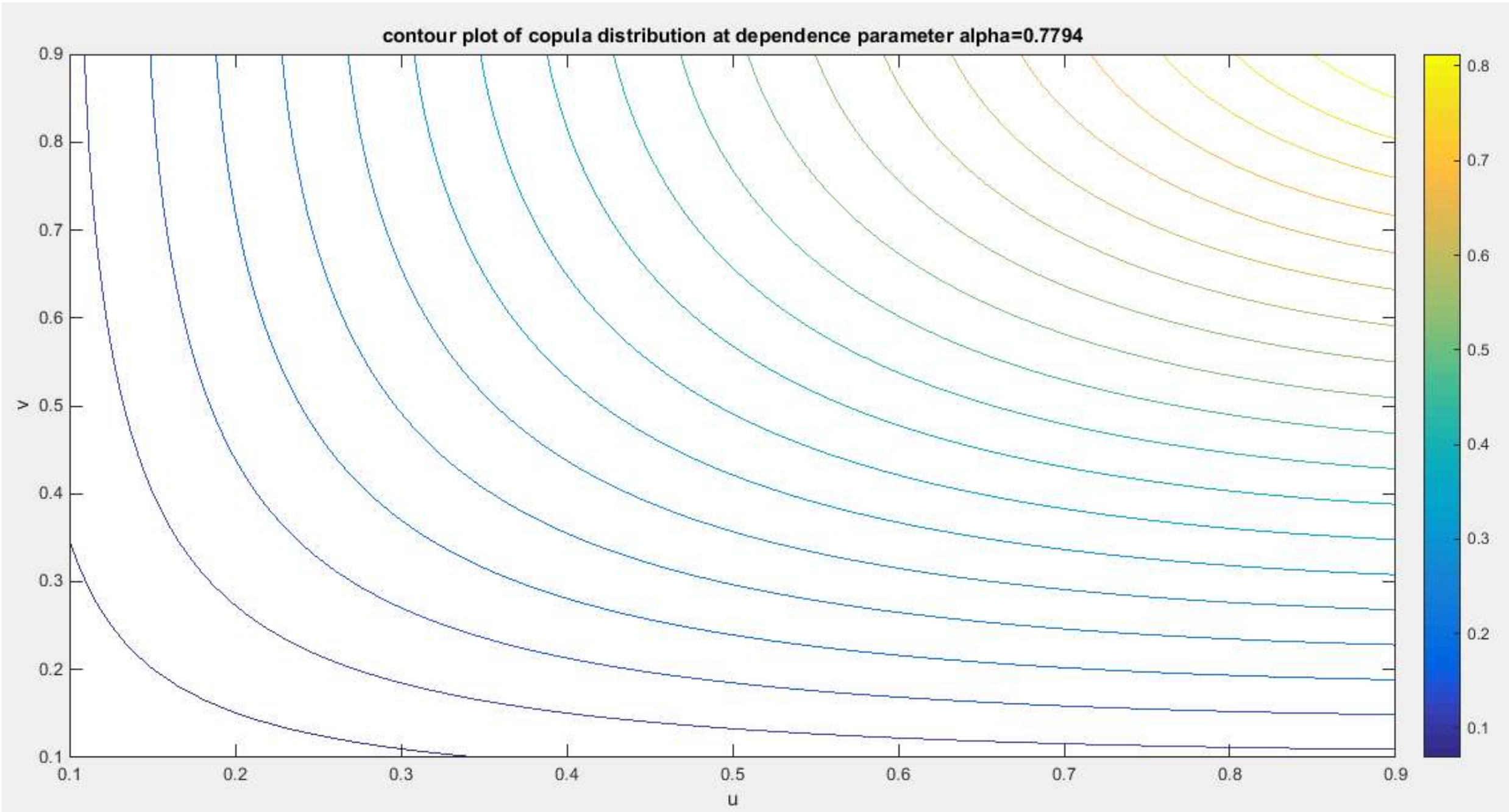

Fig. 24 shows the contour plot of the copula distribution with parameter 0.7794

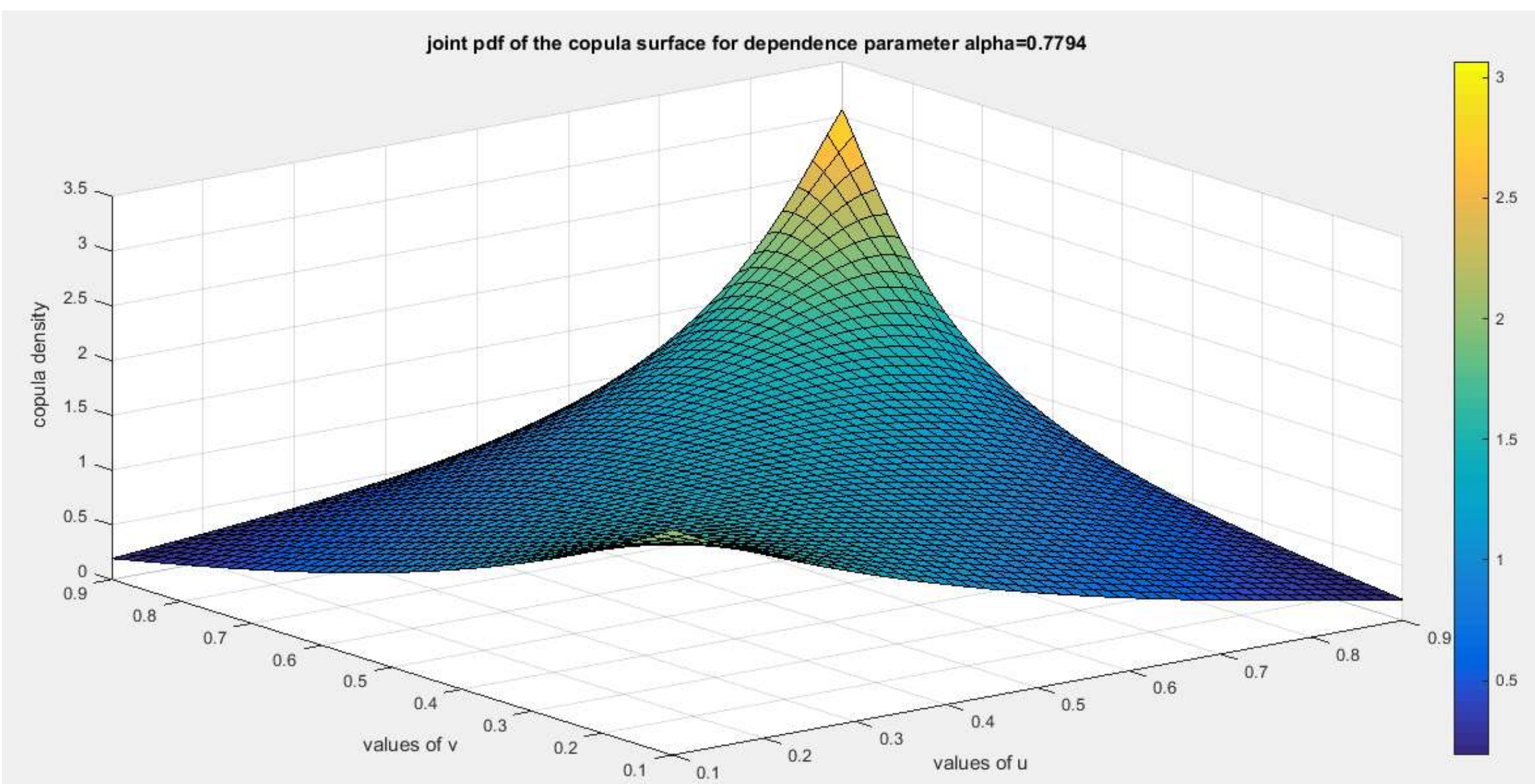

Fig. 25 shows the surface of the copula density with dependency parameter 0.7794

**Another 2 pair of variables for the second copula:**

For the OECD data, another two pair of variables are analyzed the household net wealth and the personal earnings. The household net wealth measures the total value of assets possessed by the household minus their obligations like debts, loans, and mortgages. These data are available in USD, and converted to unit interval through division by 1000000. The personal earnings measure the average income earned by persons, typically on an annual or monthly basis, reflecting the income levels and earning potential of the workforce.

These data are also expressed in USD and converted to belong to unit interval through division by 100000. These data are available for 24 countries in this same data platform as shown in Table (5).

Table (5): Household net wealth indicator versus personal earnings indicator in the OECD countries.

| | Austria | Belgium | Canada | Chile | Denmark | Estonia | Finland |
|---|---|---|---|---|---|---|---|
| Household net wealth | 309637 | 447607 | 478240 | 135787 | 149864 | 188627 | 230032 |
| Personal earnings | 53132 | 54327 | 55342 | 26729 | 58430 | 30720 | 46230 |
| | France | Germany | Greece | Hungary | Ireland | Italy | Latvia |
| Household net wealth | 298639 | 304317 | 148323 | 150296 | 370341 | 295020 | 79245 |
| Personal earnings | 45581 | 53745 | 27207 | 25409 | 49474 | 37769 | 29876 |
| | Lithuania | Luxembourg | Netherlands | Poland | Portugal | Slovak Republic | Slovenia |
| Household net wealth | 182039 | 941162 | 248599 | 233221 | 255303 | 171425 | 233286 |
| Personal earnings | 31811 | 65854 | 58828 | 32527 | 28410 | 23619 | 41445 |
| | Spain | United Kingdom | United states | | | | |
| Household net wealth | 366534 | 524422 | 684500 | | | | |
| Personal earnings | 37922 | 47147 | 69392 | | | | |

Table (6): descriptive statistics of the indicators (household net wealth and personal earnings)

| indicator | min | mean | Standard deviation | skewness | kurtosis | $25_{percentile}$ | $50_{percentile}$ | $75_{percentile}$ | max |
|---|---|---|---|---|---|---|---|---|---|
| Household net wealth | 0.0792 | 0.3094 | 0.1951 | 1.8194 | 6.9370 | 0.1767 | 0.252 | 0.3684 | 0.9412 |
| Personal earnings | 0.2362 | 0.4296 | 0.1369 | 0.2503 | 1.8851 | 0.3030 | 0.4351 | 0.5404 | 0.6939 |

Table (7): Kendall tau coefficient of the indicators

| | Household net wealth | Personal earnings |
|---|---|---|
| Household net wealth | 1 | 0.5507 (0.0001) |
| Personal earnings | 0.5507 (0.0001) | 1 |

Table (8): statistical validity indices of the two indicators

| indicator | Estimated theta | variance | AIC | CAIC | BIC | HQIC | KS-test | Ho | p-value of KS-test |
|---|---|---|---|---|---|---|---|---|---|
| **Household net wealth** | **1.2689** | **0.0087** | **-10.423** | **-10.2411** | **-9.2449** | **-10.1104** | **0.2341** | **Fail to reject** | **0.1219** |
| **Personal earnings** | **1.0418** | **0.0058** | **-12.9361** | **-12.7543** | **-11.758** | **-12.6236** | **0.2496** | **Fail to reject** | **0.0837** |

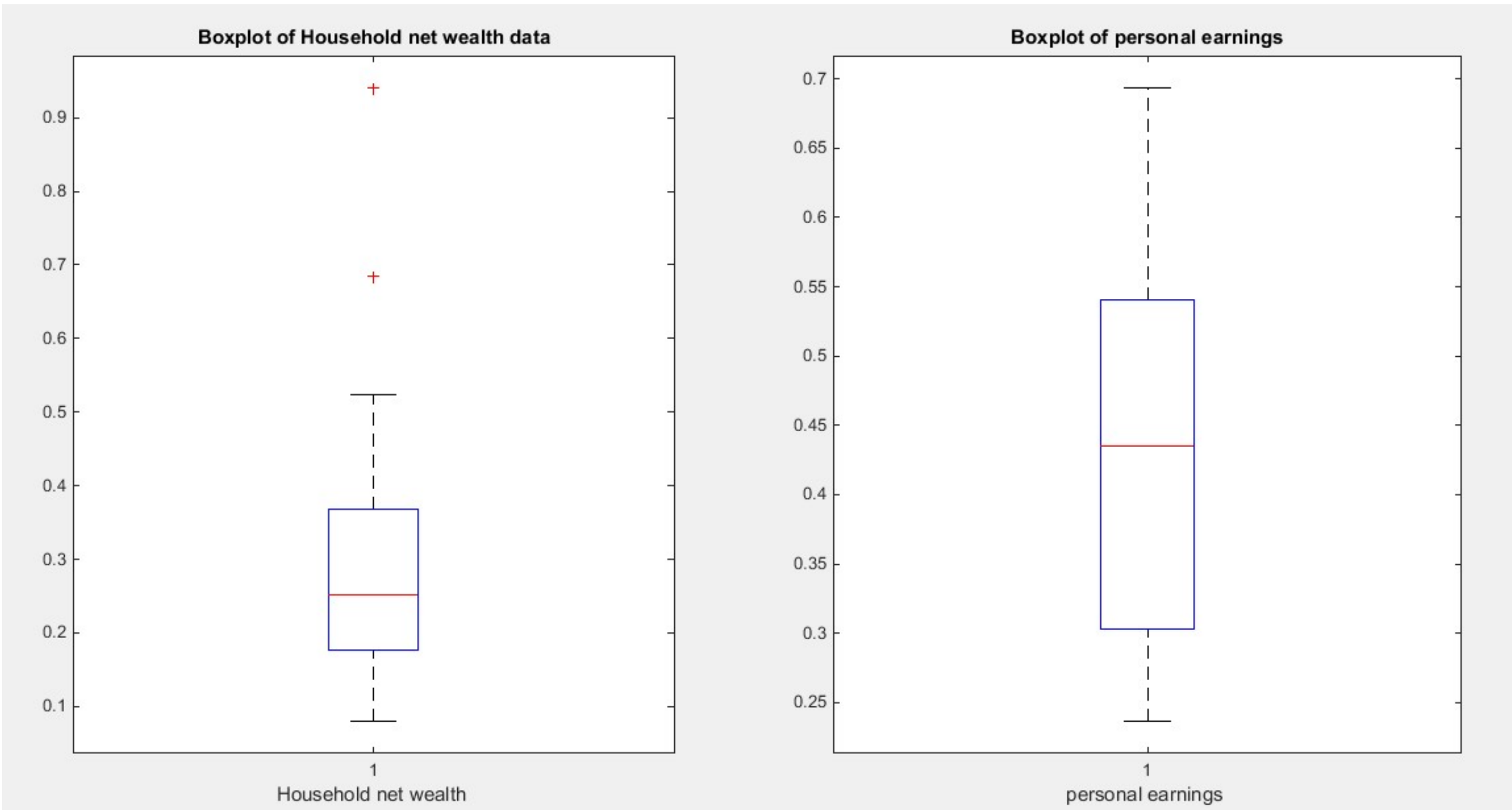

Fig. 26 shows Boxplot of household net wealth and the personal earnings denoting right side skewness.

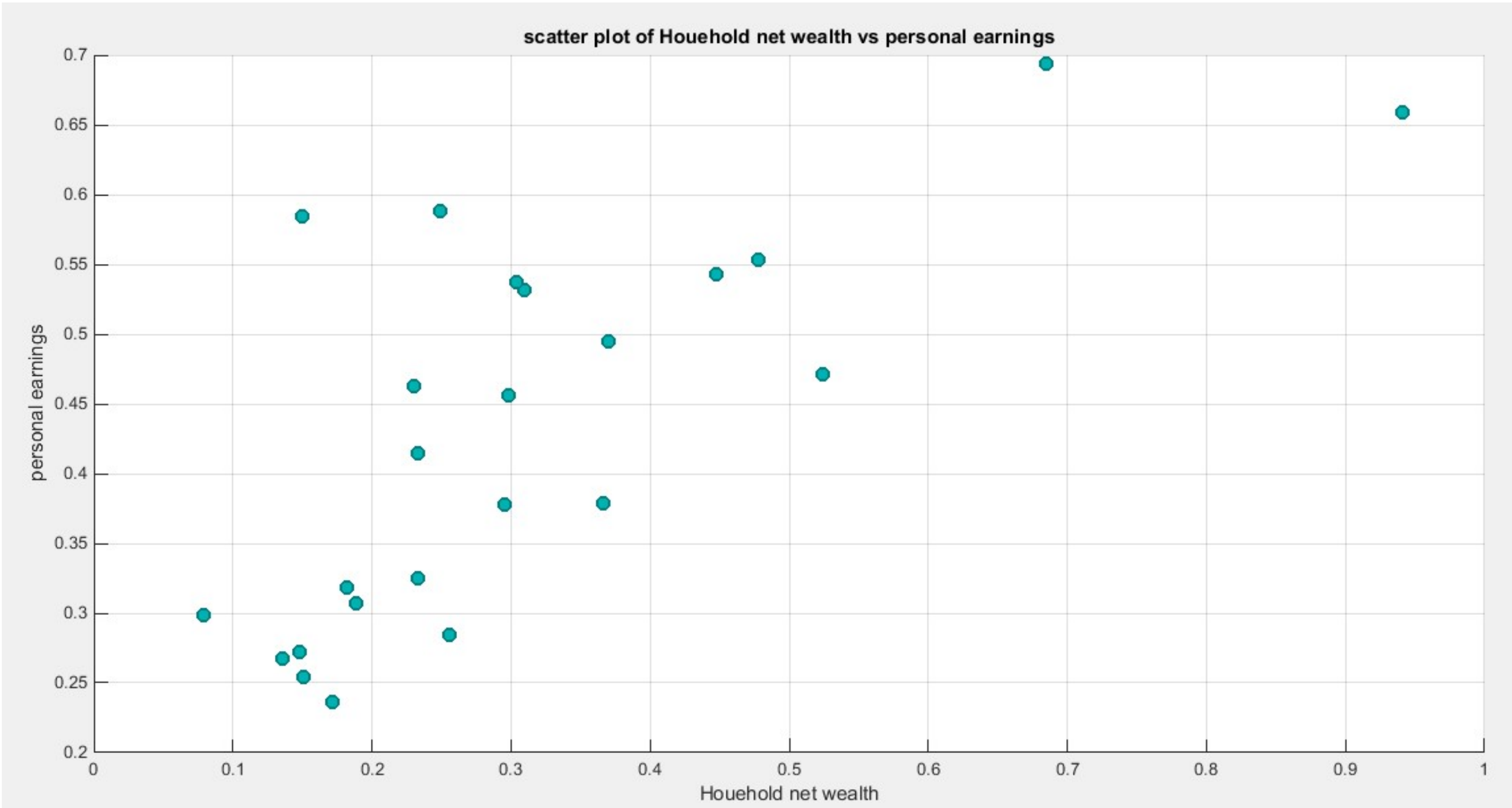

Fig. 27 shows scatter plot of the household net wealth vs personal earnings and the data are mainly concentrated on the left lower corner and the center pointing to the presence of lower tail dependency.

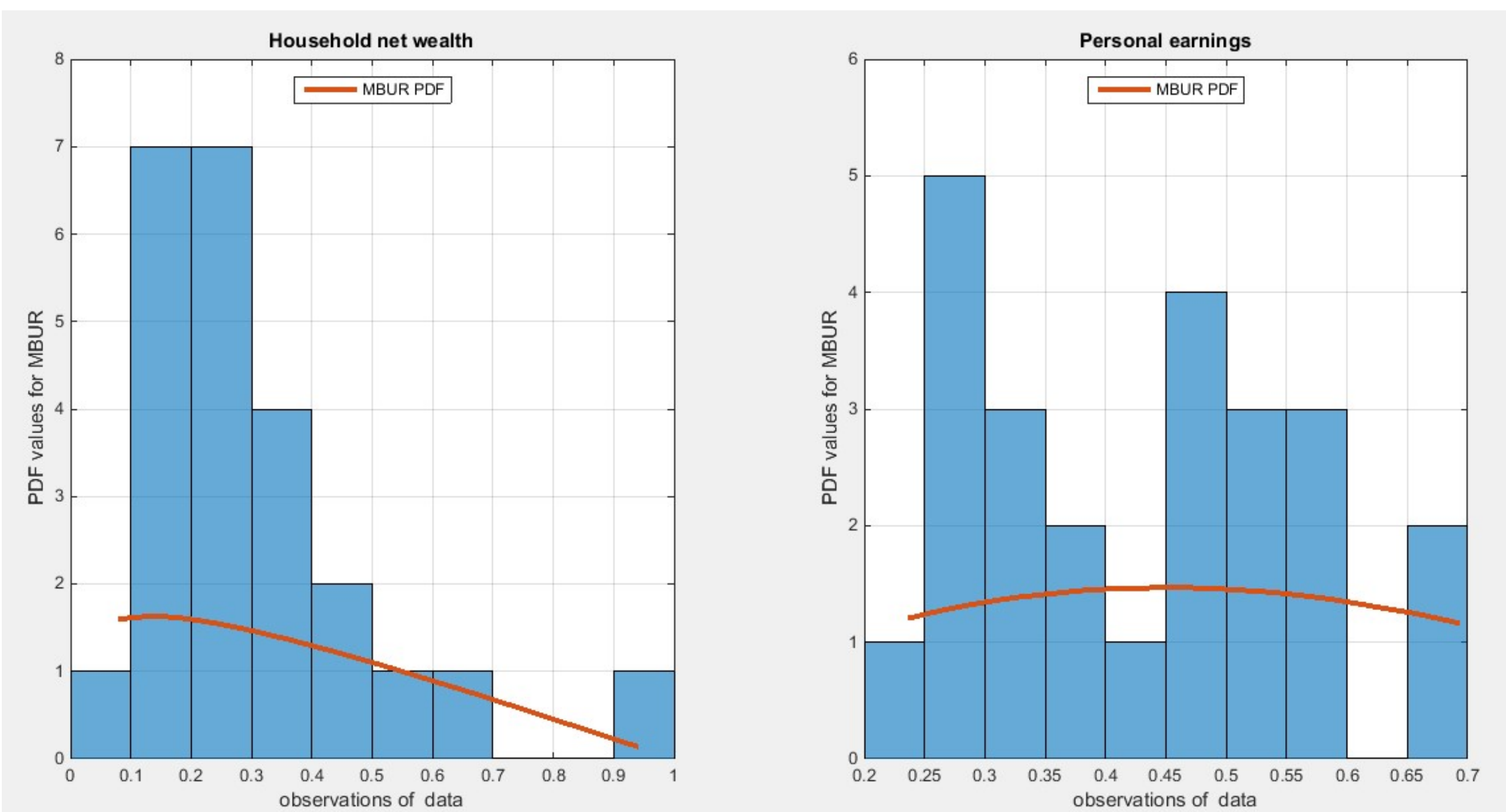

Fig. 28- shows the histogram of the two indicators which exhibit right sided skewness. The graph also shows the fitted MBUR curve for each indicator. The household net wealth indicator exhibits more kurtosis and skewness than the personal earnings indicator.

***Second copula Model the dependency between the 'household net wealth' and 'personal earnings'***: The empirical tau is 0.5507 and the theoretical tau is 0.5507 with dissimilarity value of 7.1985e-05. The estimated dependency parameter, theta $\hat{\theta}$ , is 0.6703 and estimated variance of 0.0018. The Cramer Von Mises test (CVM) is 0.4281. to test the null hypothesis for dependency, resampling using the Metropolis Hasting algorithm of MCMC procedure was performed. The propositional null hypothesis was stated that the population dependency parameter equals the estimated theta against the alternative hypothesis that the population dependency parameter does not equal the estimated theta. Moreover, the null hypothesis for the equality of the population Kendall tau with the theoretical tau is tested against the alternative hypothesis of not being equal to each other. The null hypothesis for the CVM suggesting that the sampling CVM equals the observed CVM attained from the estimation procedure was also investigated. These are all two sided tests. The sampling distribution of each of the previously mentioned indices, thetas, tau and CVM is shown in Figures (29-31). The descriptive statistical indices are shown in each figure. Analysis of the figures denoted that the indices are located within the acceptance region between the 2.5th and 97.5th quantiles. So the null hypotheses fail to be rejected. The copula can be said to pass some of the goodness of fit (GOF) tests for this data. And the dependence parameter mockups, within the milieu of this copula, the connection between the two variables. The confidence interval (CI) for the estimated dependence parameter is (-0.0408, 0.9268), for the theoretical tau; it is (0.1314, 0.9756), and for CVM; it is (0.0368, 0.4871).

According to the sandwich variance the variance covariance matrix for this model is:

$$var - cov = \begin{bmatrix} 0.1418 & 0.0964 & 0.0002 \\ 0.0964 & 0.0657 & 0.0002 \\ 0.0002 & 0.0002 & 0.0001 \end{bmatrix}$$

The validity indices for the second stage are:

$AIC = -10.833$ , $CAIC = -10.6512$ $BIC = -9.655$ $HQIC = -10.5205$

The validity indices for the whole model is the summation of these indices for marginal (first stage) and copula (second stage):

$AIC = -34.1921$ , $CAIC = -33.6466$ $BIC = -30.6579$ $HQIC = -33.2544$

Theoretical upper tail coefficient is 0.6346. At quantile 0.7; the empirical upper tail coefficient in the direction of the 'households net wealth' variable is 0.6667 (2 points showing joint match out of 3 points in the upper tail at this threshold, 0.7) and in the direction of the 'personal earnings' variable is 1 (2 points showing joint match out of 2 points in the upper tail at this threshold, 0.7) while the confidence interval of the empirical distribution obtained from bootstrap sampling under the null hypothesis is [333,1] in the ' households net wealth' and [0.2857,1] in the direction of 'personal earnings'. This is expected because of the small sample size and left skewness of data.

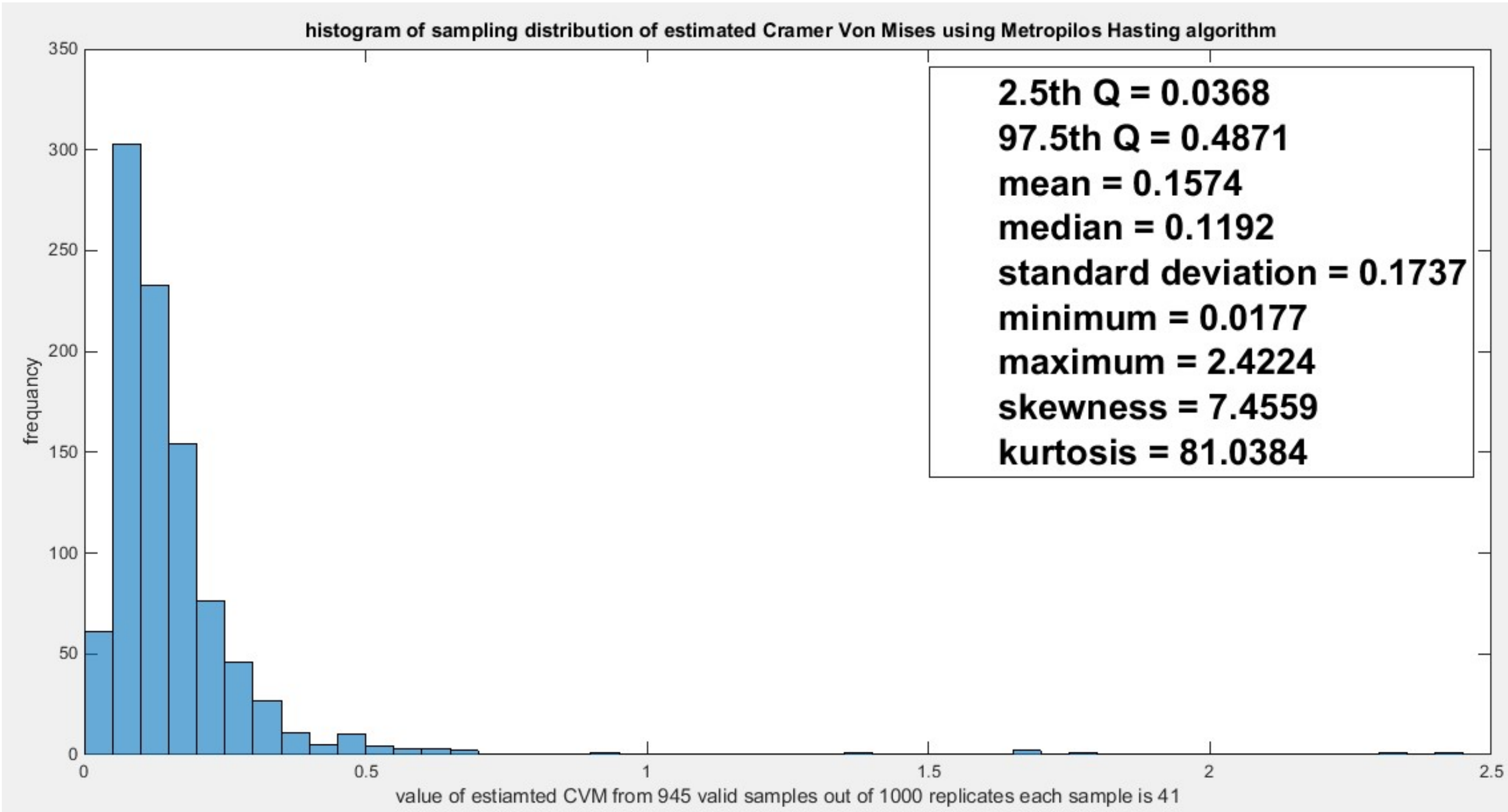

Fig. 29 shows the sampling distribution of the CVM. The data obtained CVM is 0.4281 and it lays between the 2.5th and the 97.5th quantiles so the null hypothesis is failed to be rejected and the copula fits the data well. The p-value (probability of values less than or equal to 0.4281) is 0.9683, in other words, this 0.4281 does not lay in either of the tail regions of rejection.

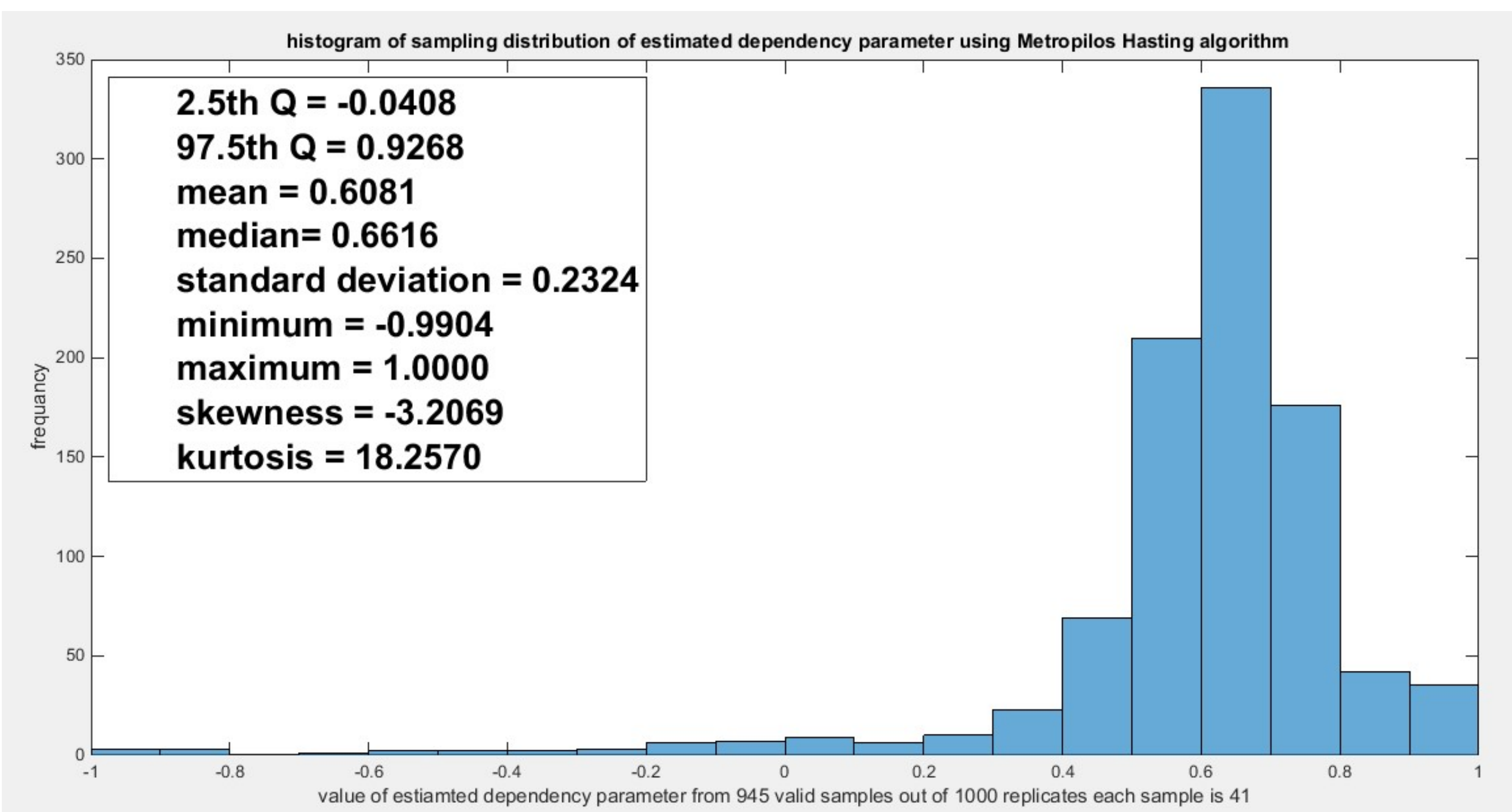

Fig. 30 shows the sampling distribution of the dependency parameter. The 0.6703 lays in the acceptance region between the two quantiles, the 2.5th and the 97.5th quantiles. So the null hypothesis fails to reject its assumption. The p-value (probability of the values being less than or equal to this 0.6703) is 0.5651, which indicates that 0.6703 is not placed in either zones of rejection.

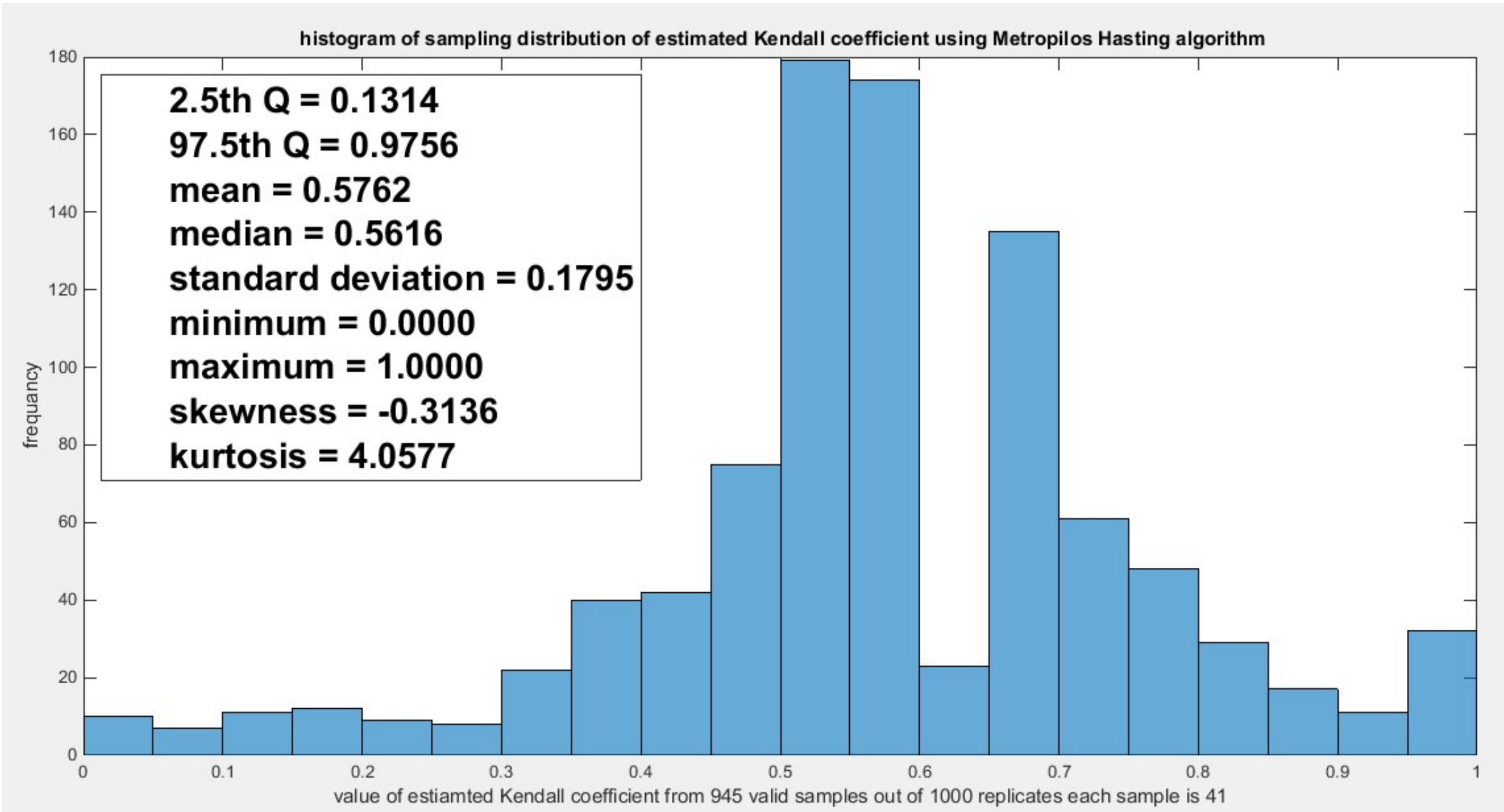

Fig. 31 shows the sampling distribution of the Kendal tau coefficient. Its value 0.5507 is between the two quantiles 2.5th and 97.5th which is the acceptance zone. So the null hypothesis fails to reject its assumption. The p-values (probability of the values less than or equal to 0.5507) is 0.4423 which denotes that this value 0.5507 is far away from both rejection zones on the right and the left.

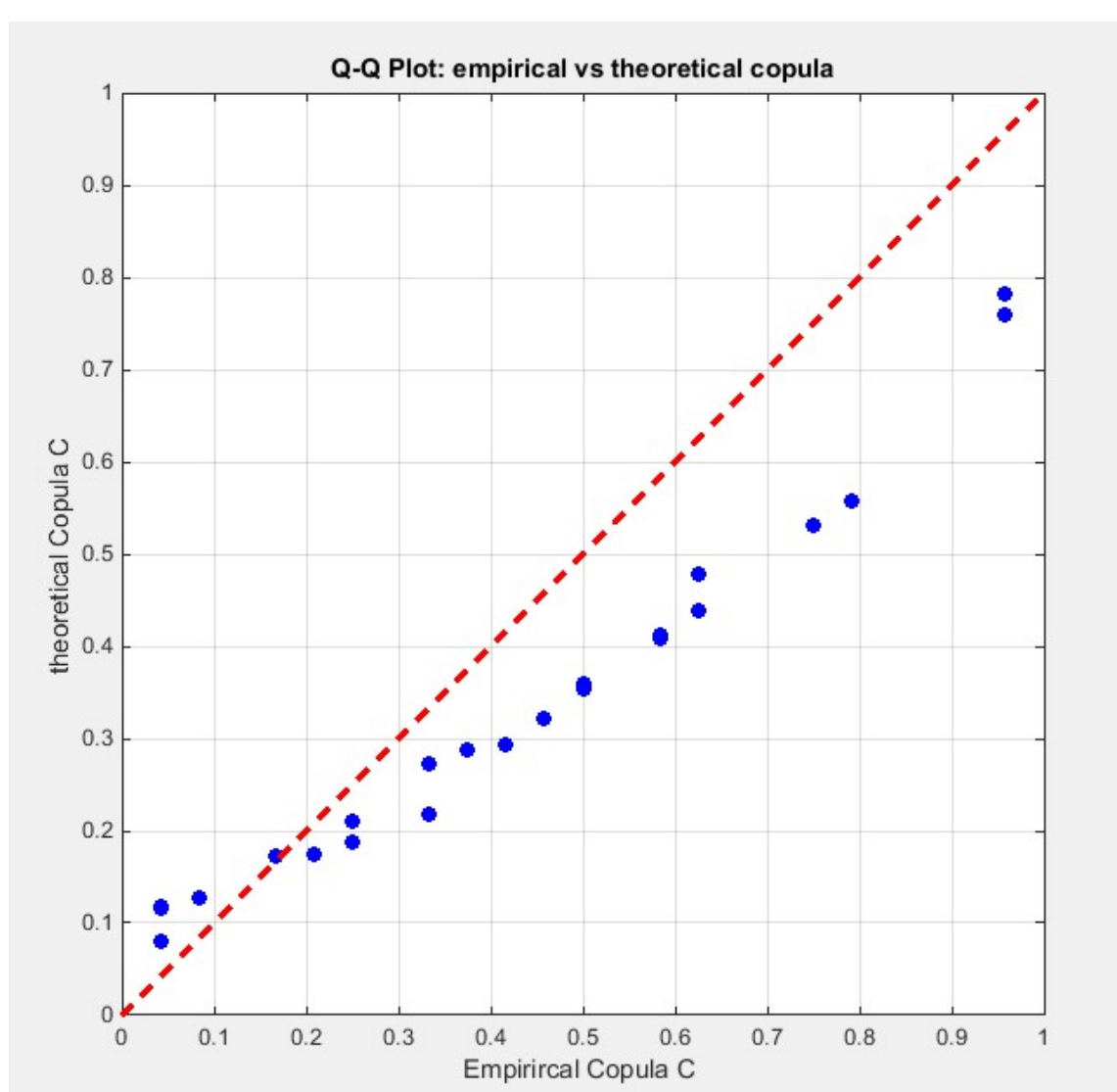

Fig. 32 shows the P-P plot of the empirical vs the theoretical copula. The theoretical copula shows near perfect alignment of its lower part with the diagonal line than its upper part.

Figure 32 shows the P-P plot of the empirical and theoretical copula. Figure 33-35 show the surface of the density copula, the contour plot of the density copula and the contour plot of the distribution copula at the estimated dependency parameter.

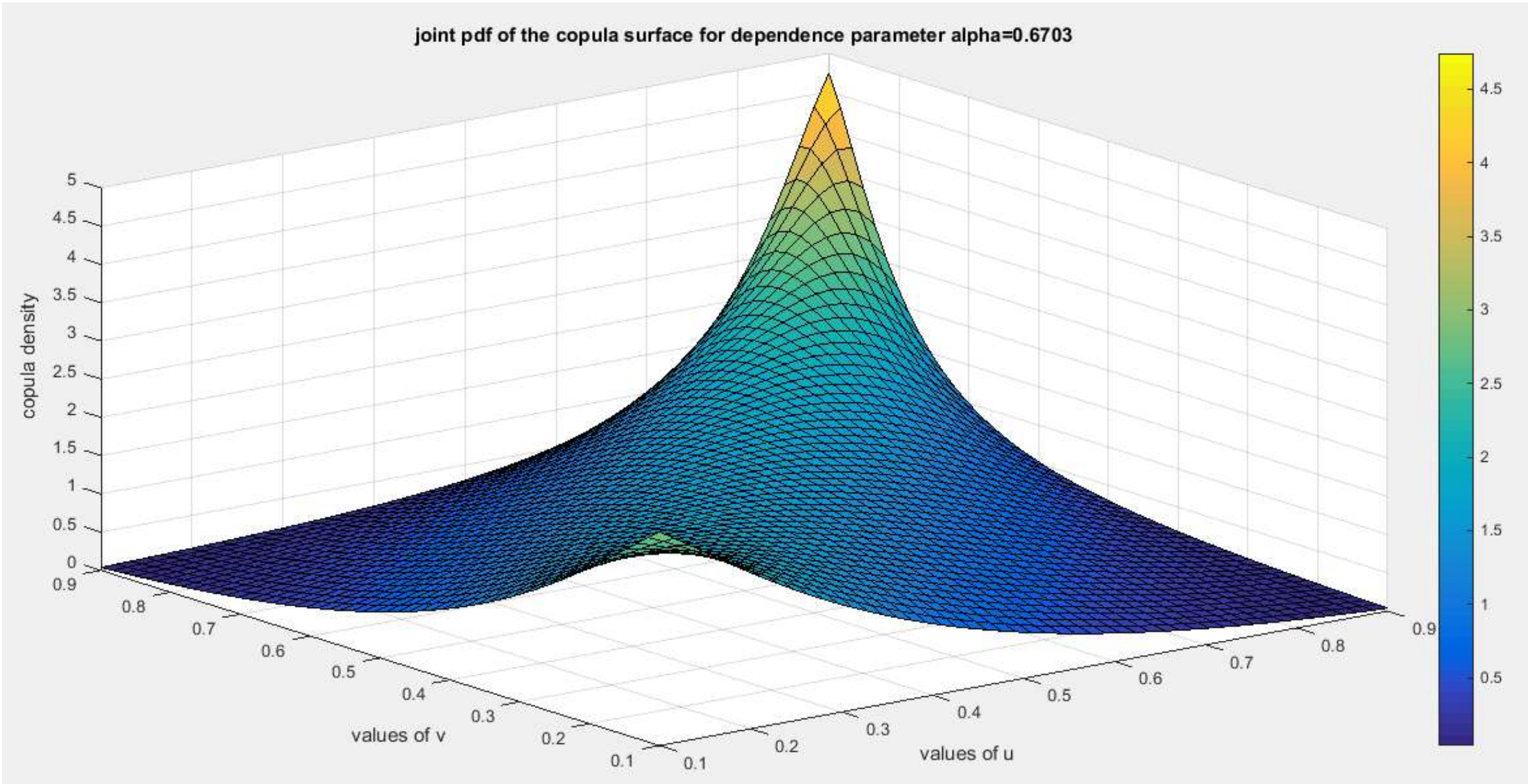

Fig. 33 shows the surface of the copula density with dependency parameter 0.6703

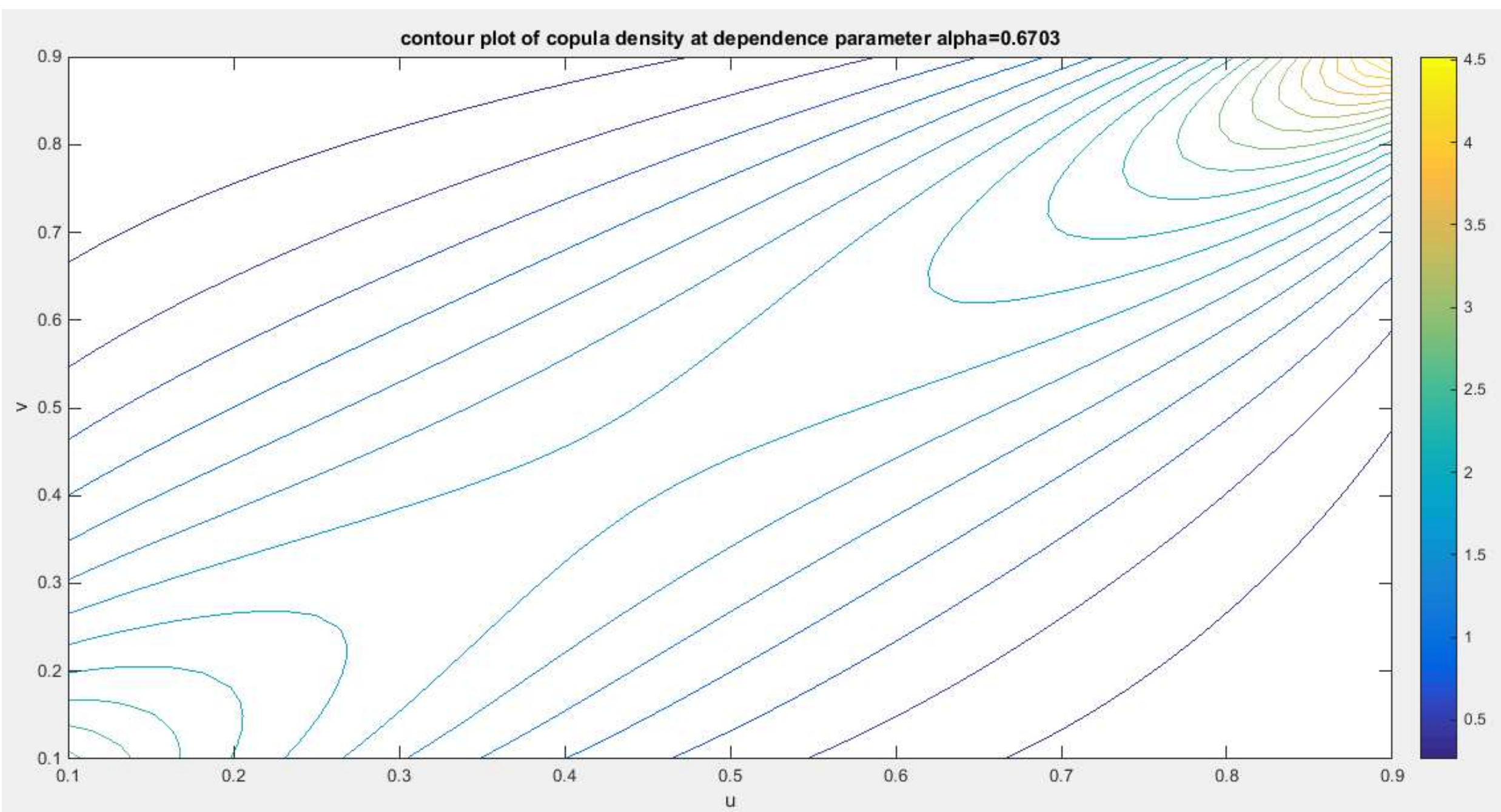

Fig. 34 shows the contour plot of the copula density with dependency parameter 0.6703

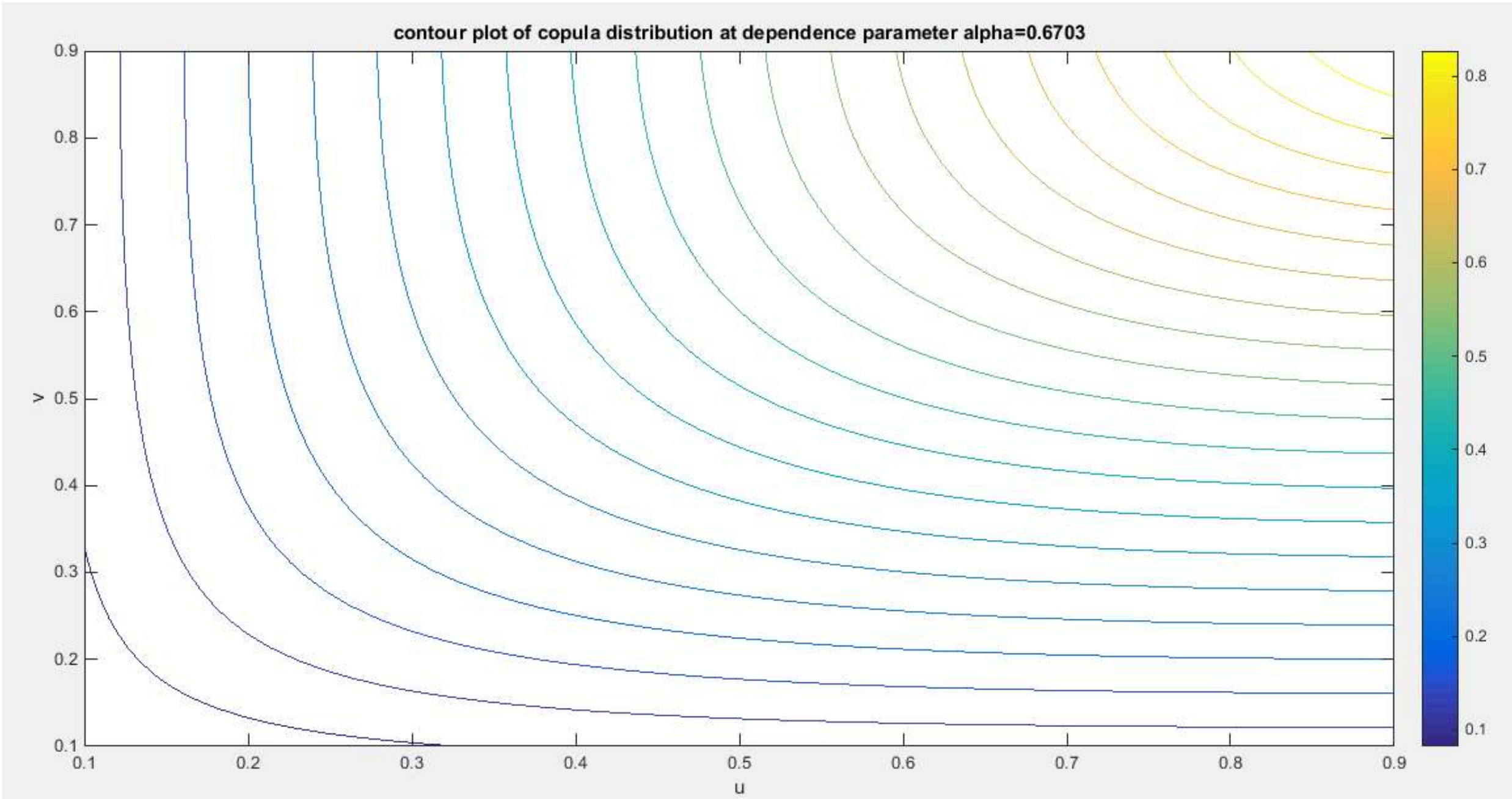

Fig. 35 shows the contour plot of the copula distribution with dependency parameter 0.6703

## Section 5: Conclusions

Both copulas model positive dependency. They can model upper tail dependency but they have no lower tail dependency. The author discussed 2 datasets exhibiting left skewness. Both copulas fit the data well. The author also discussed one dataset exhibiting right skewness and the second copula fit this set well. Both copula fits data that has more or less similar characteristics. The complexity of the the copula density function for both copulas mandates use of IFM method of estimation. This method is preferred over the MLE for its easiness. The robust variance or sandwich variance captures the covariance between the marginal and the copula parameter estimators which are not captures by estimating the variances separately in each stage. The variance of the estimators of marginal parameters increases while the variance of the estimator of copula parameter decreases. Although the sample sizes used in the analysis were small (n=24) and moderate(n=41) both copula model the upper tail dependency well and this was supported by the sampling distribution of the empirical upper tail coefficient.

## Section 6: Future Work

These copulas can be mixed with other copulas to better model data with asymmetric dependency. Also other methods of estimation can be attempted like the semiparametric estimation and MLE. Bayesian inference can be considered.

**Declarations:**
**Ethics approval and consent to participate**
Not applicable.
**Consent for publication**
Not applicable
**Availability of data and material**
Not applicable. Data sharing not applicable to this article as no datasets were generated or analyzed during the current study.
**Competing interests**
The author declares no competing interests of any type.
**Funding**
No funding resource. No funding roles in the design of the study and collection, analysis, and interpretation of data and in writing the manuscript are declared

**Authors' contribution**

AI (Attia Iman) carried the conceptualization by formulating the goals, aims of the research article, formal analysis by applying the statistical, mathematical and computational techniques to synthesize and analyze the hypothetical data, carried the methodology by creating the model, software programming and implementation, supervision, writing, drafting, editing, preparation, and creation of the presenting work.

**Acknowledgement**

Not applicable

# APPENDIX

**Although the figures of the copula density of the first and the second copula may present some resemblance in overall shape but they differ in the mathematical behaviors and properties, the values attained for data fitness, and theoretical Kendall tau. The associated color bar can explain and can depict the differences.**

# First Copula

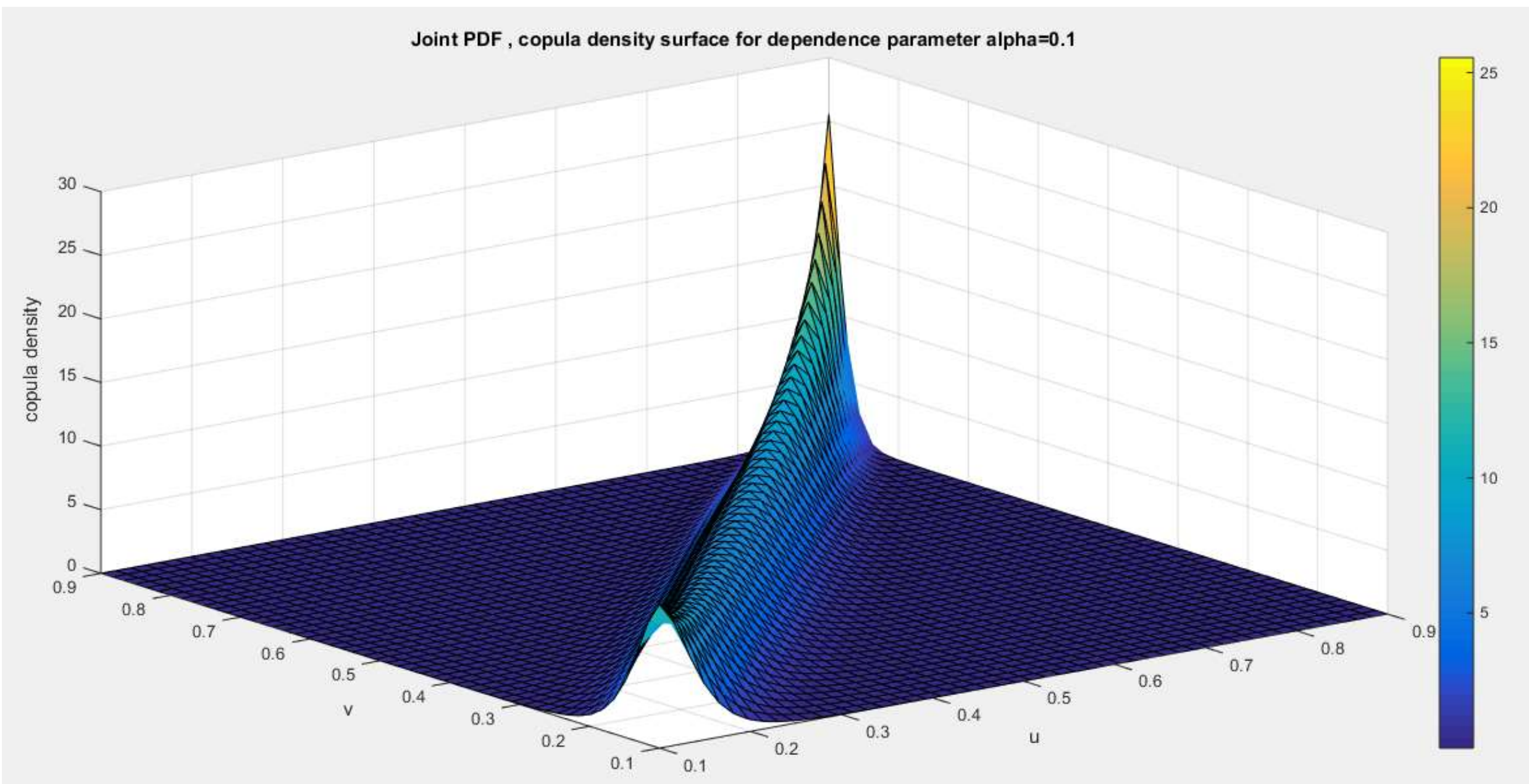

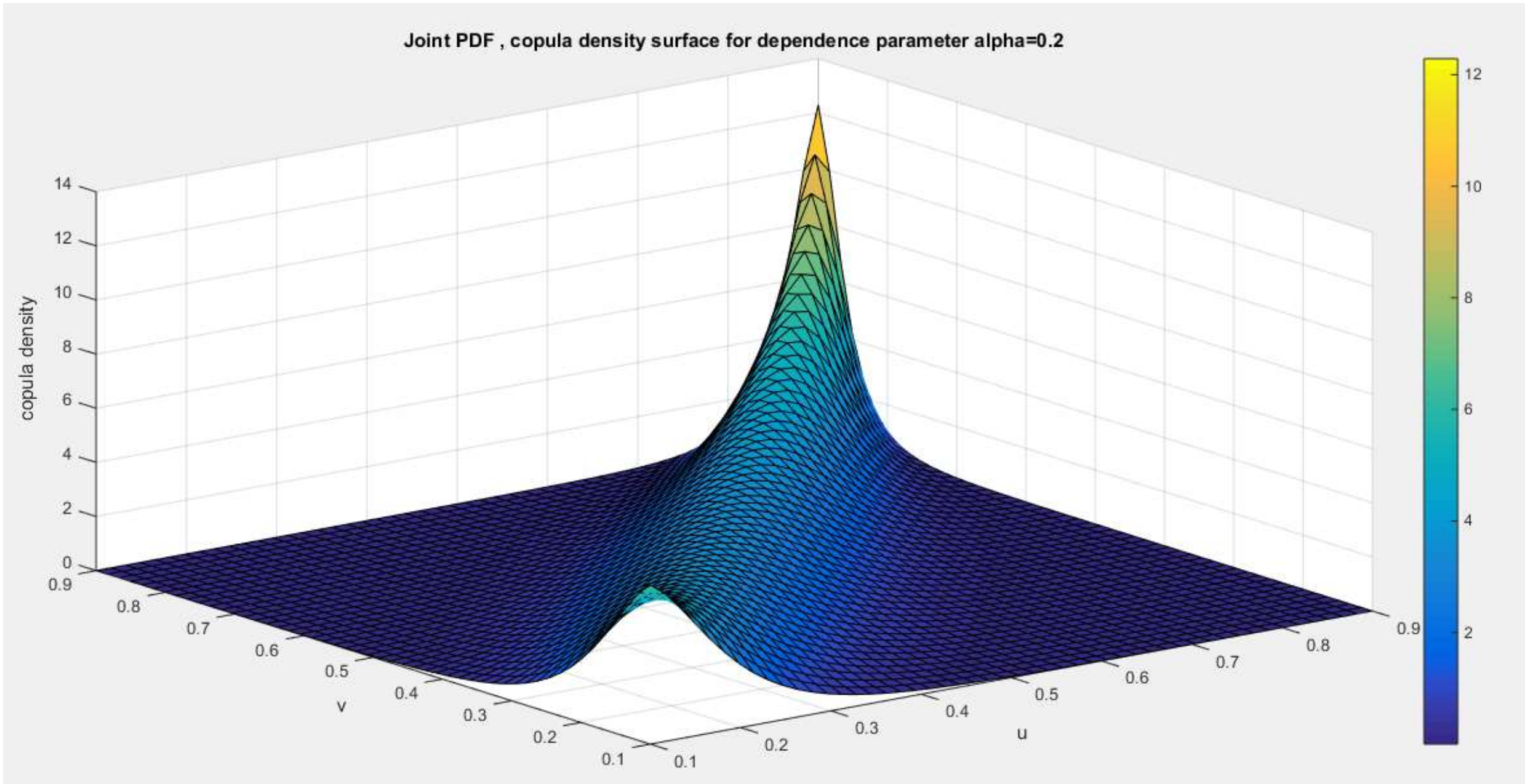
Joint PDF , copula density surface for dependence parameter alpha=0.2
copula density
v
u

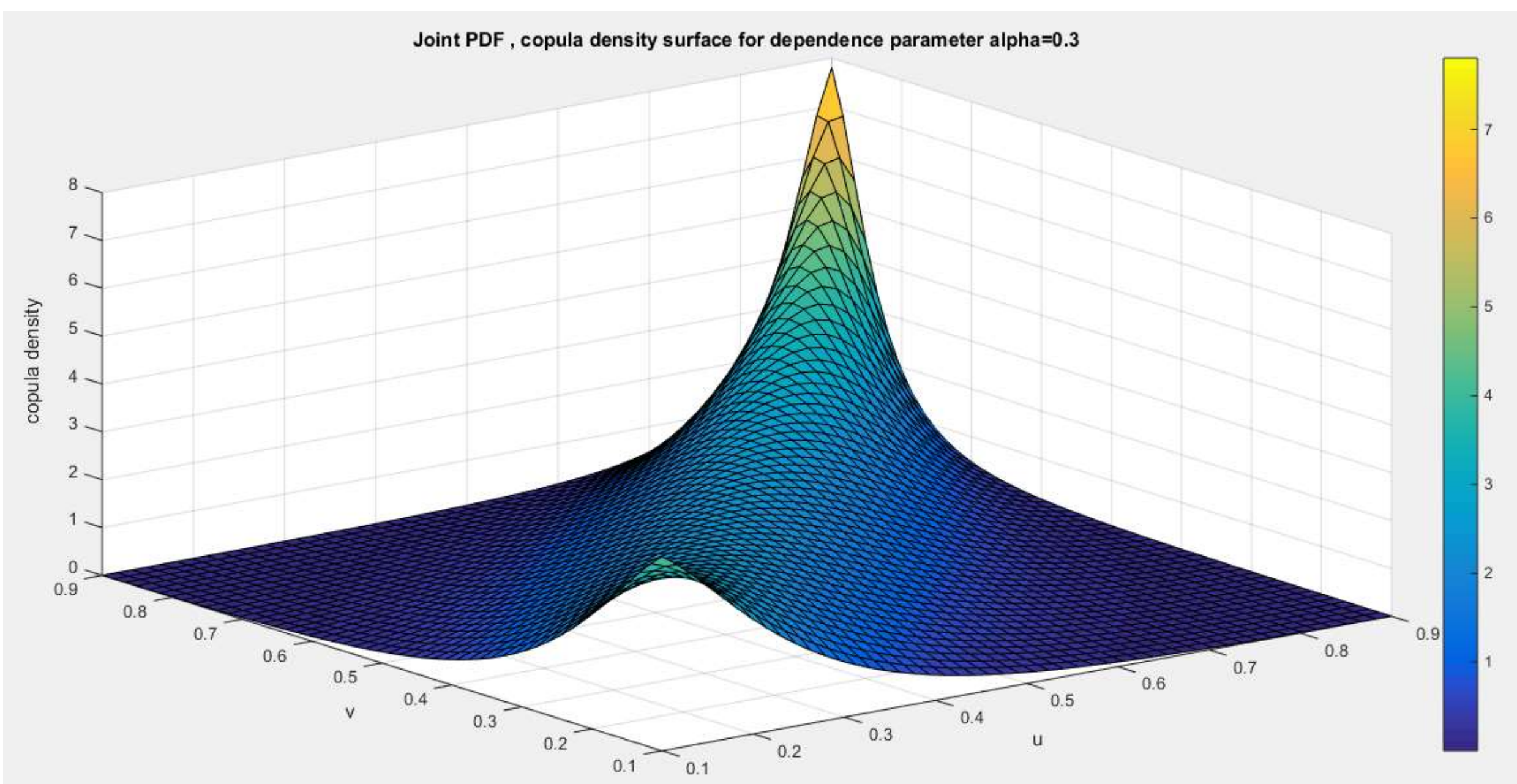
Joint PDF , copula density surface for dependence parameter alpha=0.3
copula density
v
u

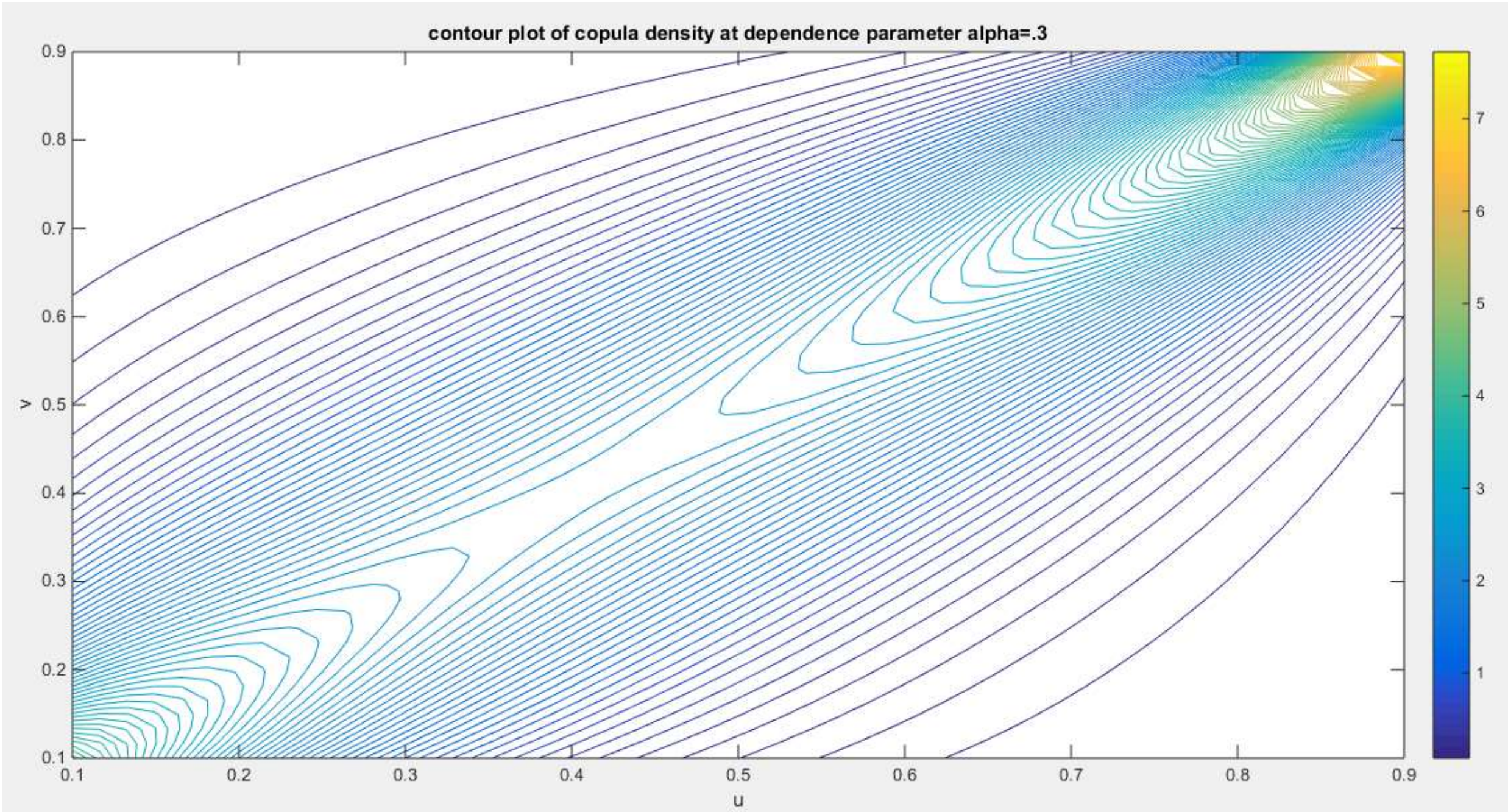
contour plot of copula density at dependence parameter alpha=.3
v
u
0.1
0.2
0.3
0.4
0.5
0.6
0.7
0.8
0.9
1
2
3
4
5
6
7

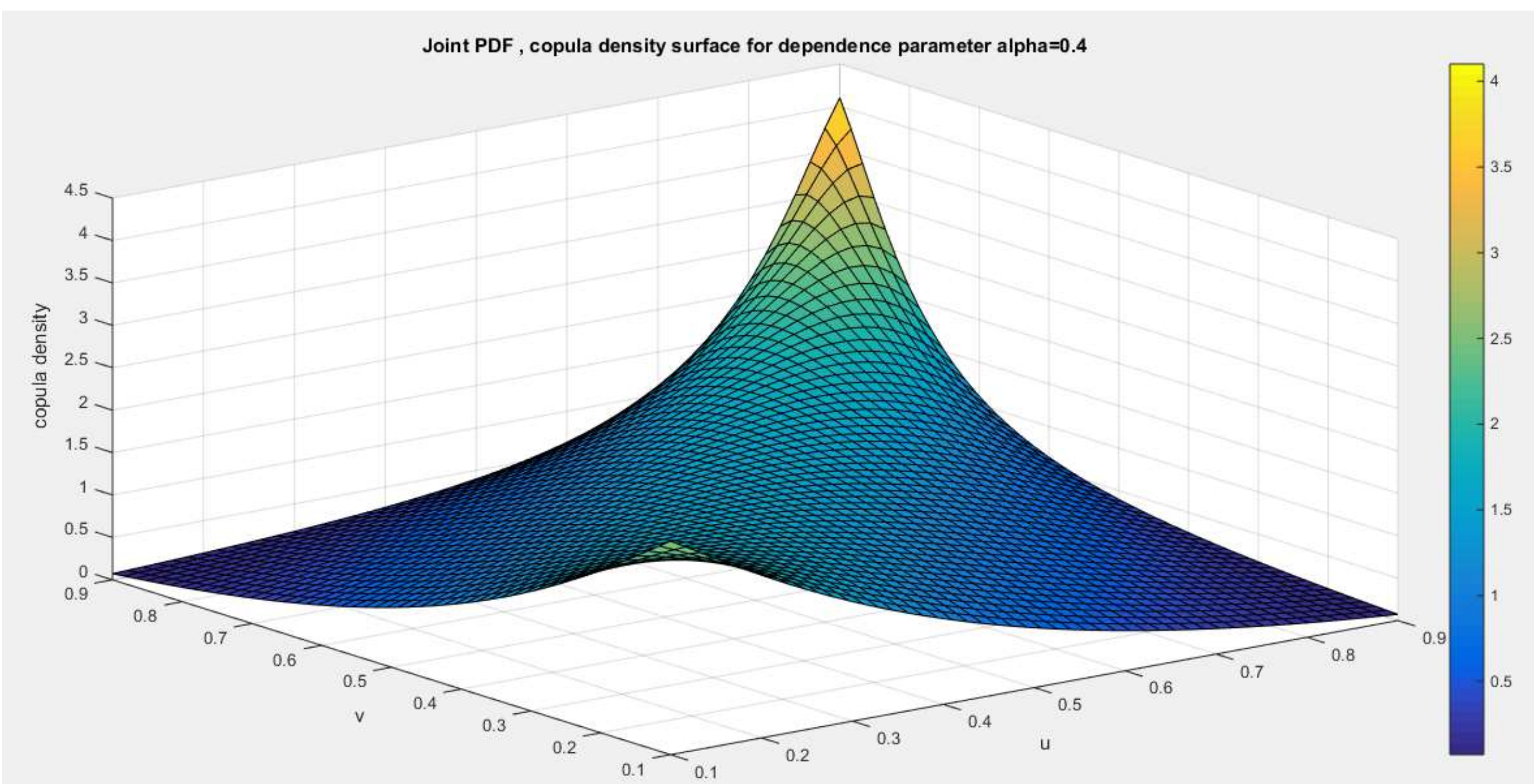
Joint PDF , copula density surface for dependence parameter alpha=0.4
copula density
v
u
0
0.5
1
1.5
2
2.5
3
3.5
4
4.5

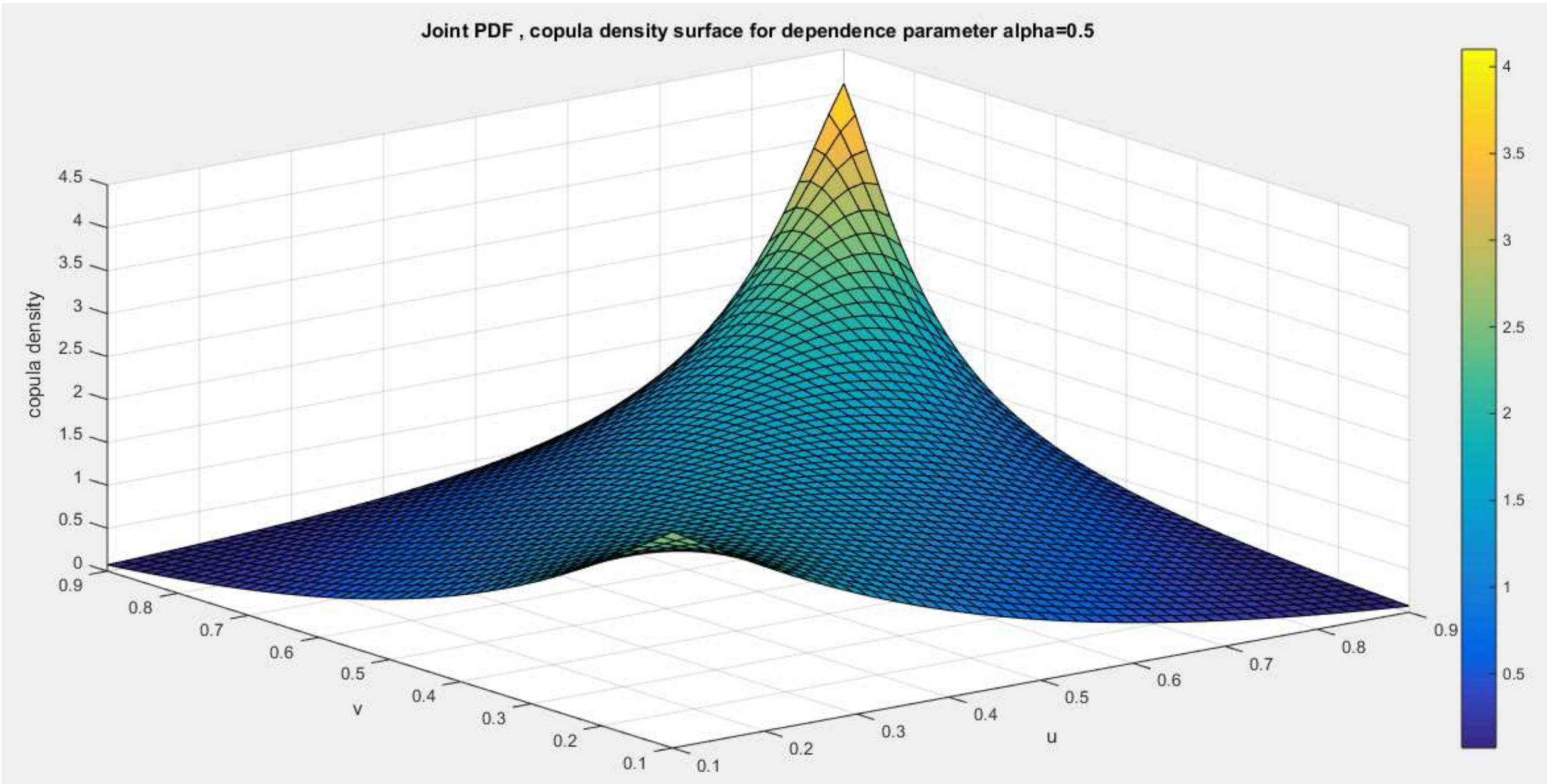
Joint PDF , copula density surface for dependence parameter alpha=0.5
copula density
v
u

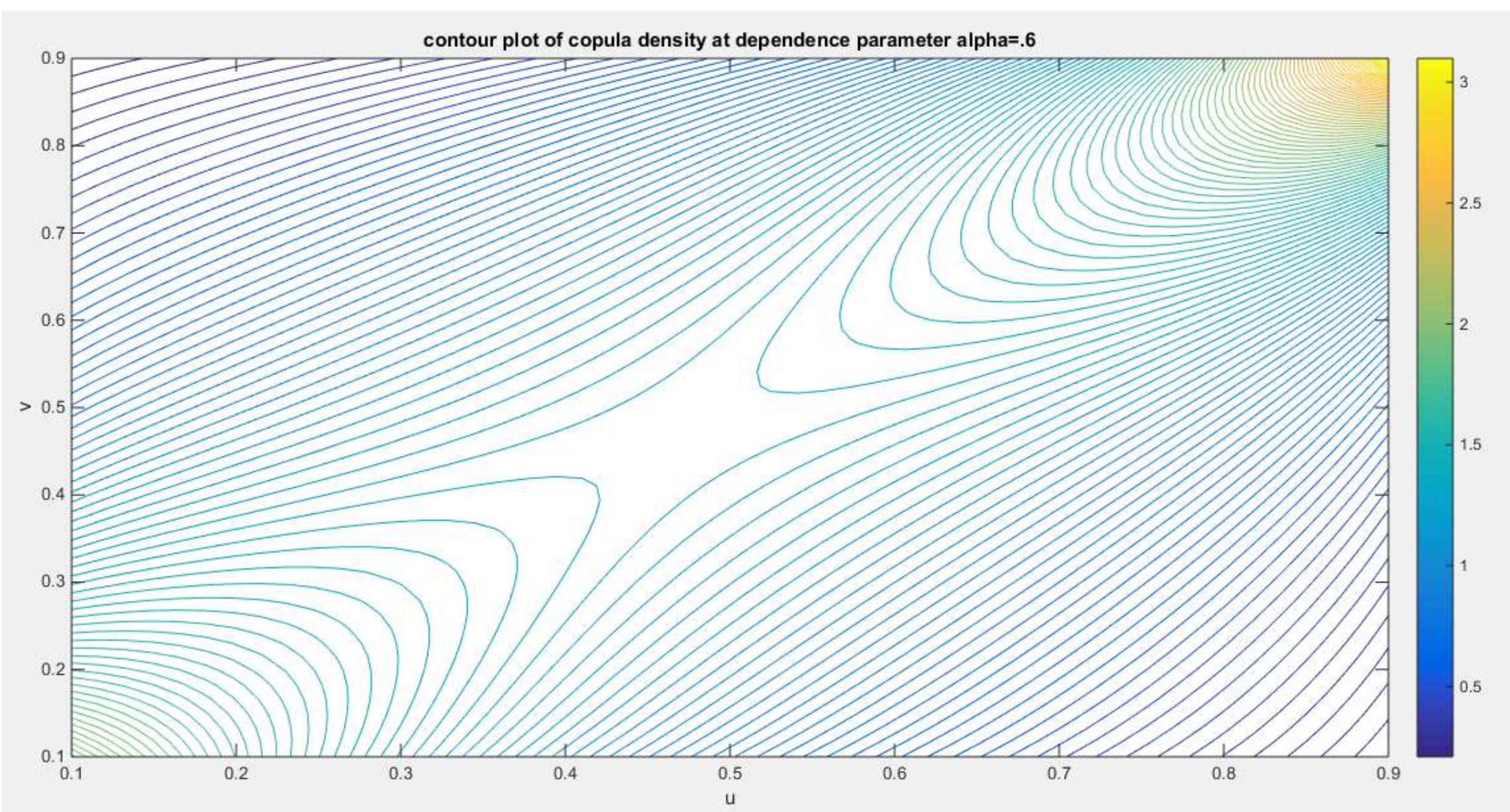
contour plot of copula density at dependence parameter alpha=.6
v
u

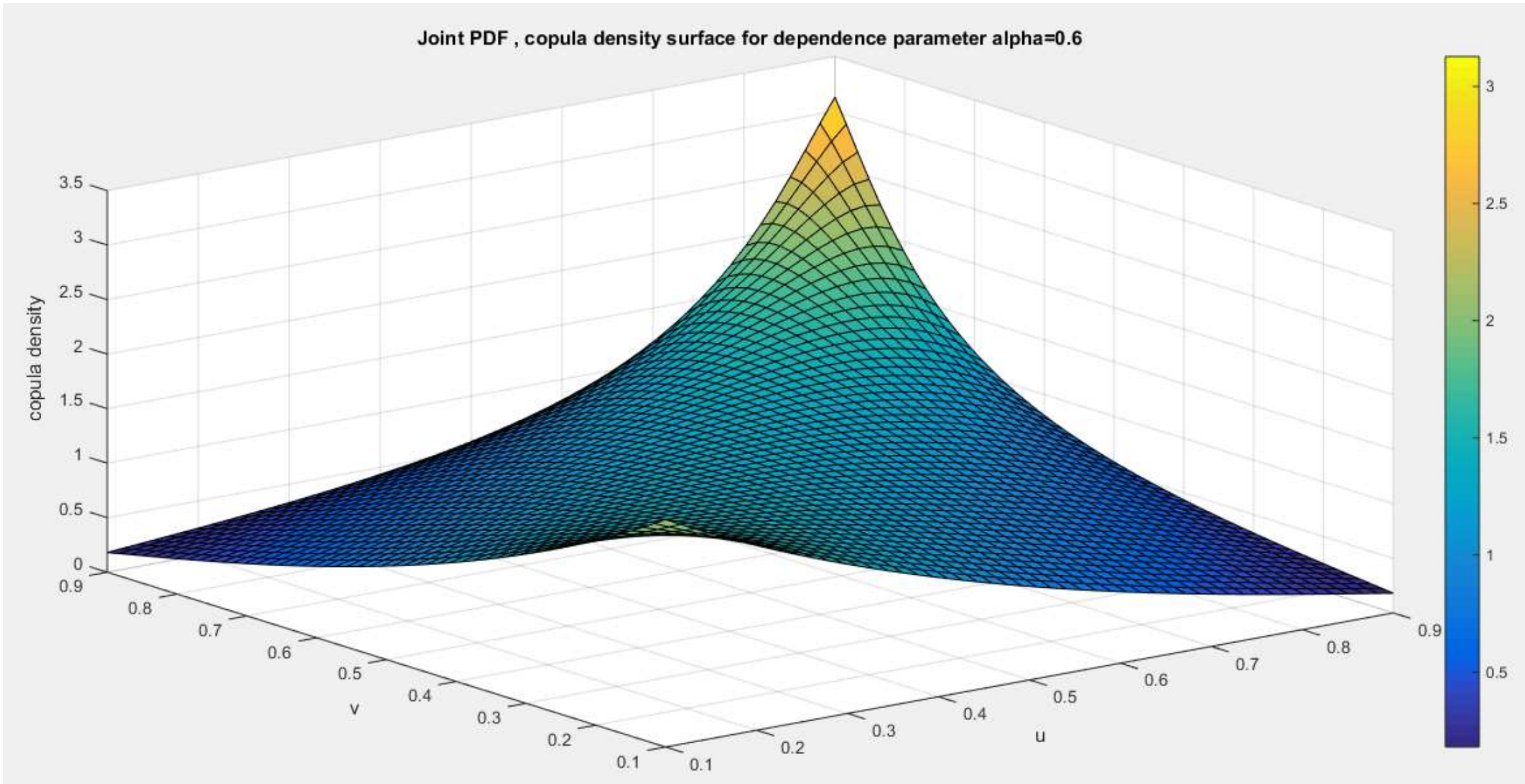
Joint PDF , copula density surface for dependence parameter alpha=0.6
copula density
v
u

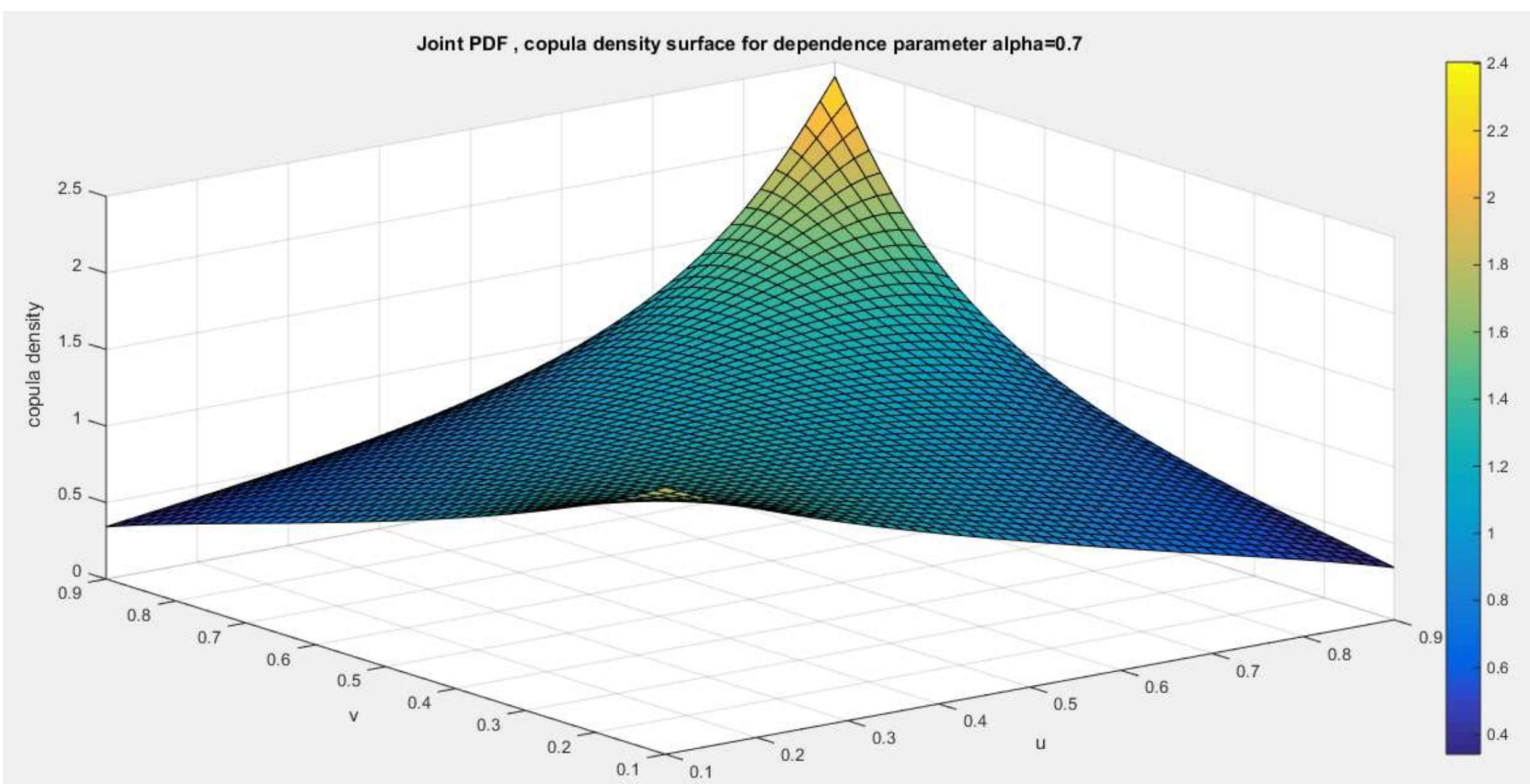
Joint PDF , copula density surface for dependence parameter alpha=0.7
copula density
v
u

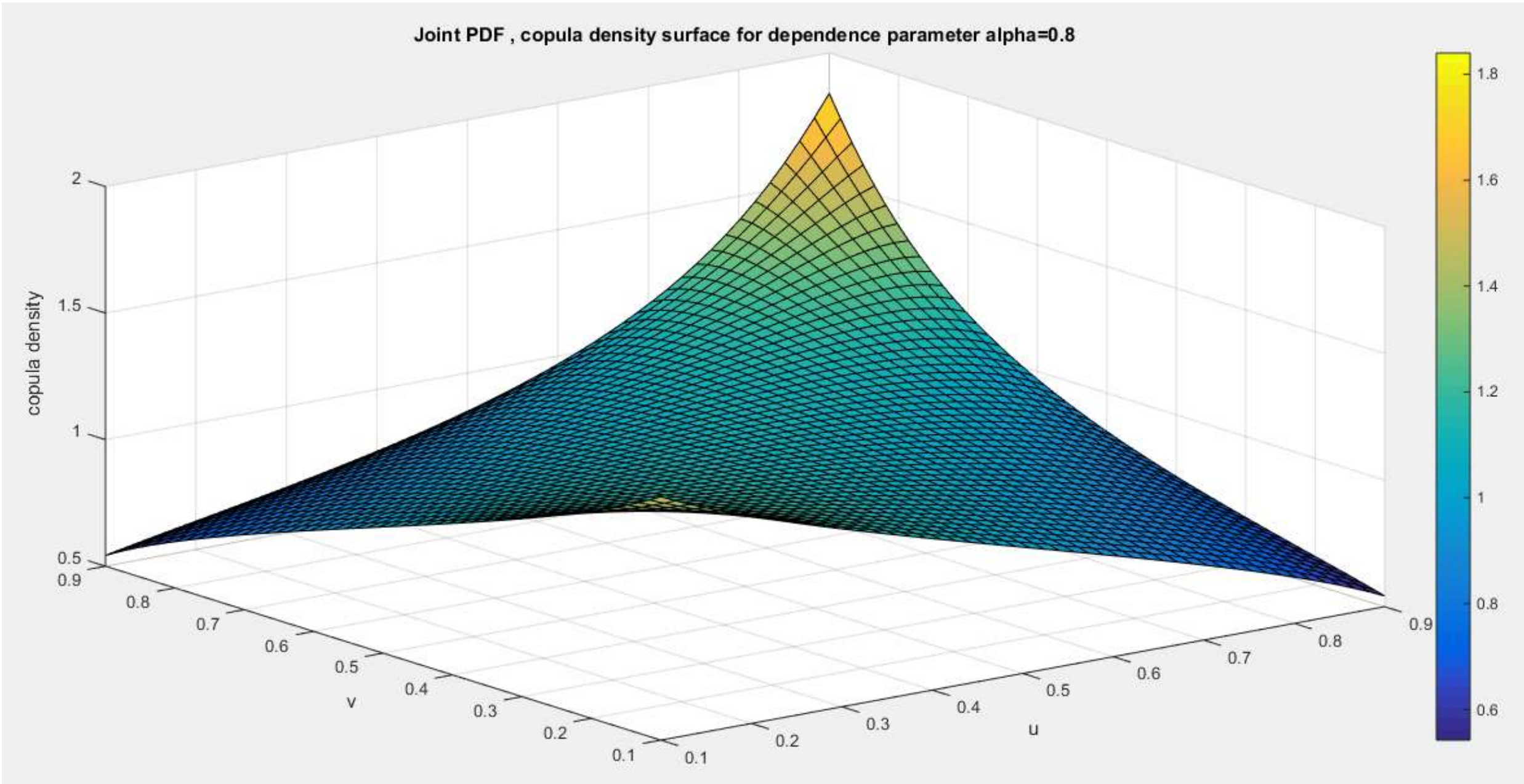
Joint PDF , copula density surface for dependence parameter alpha=0.8
copula density
v
u

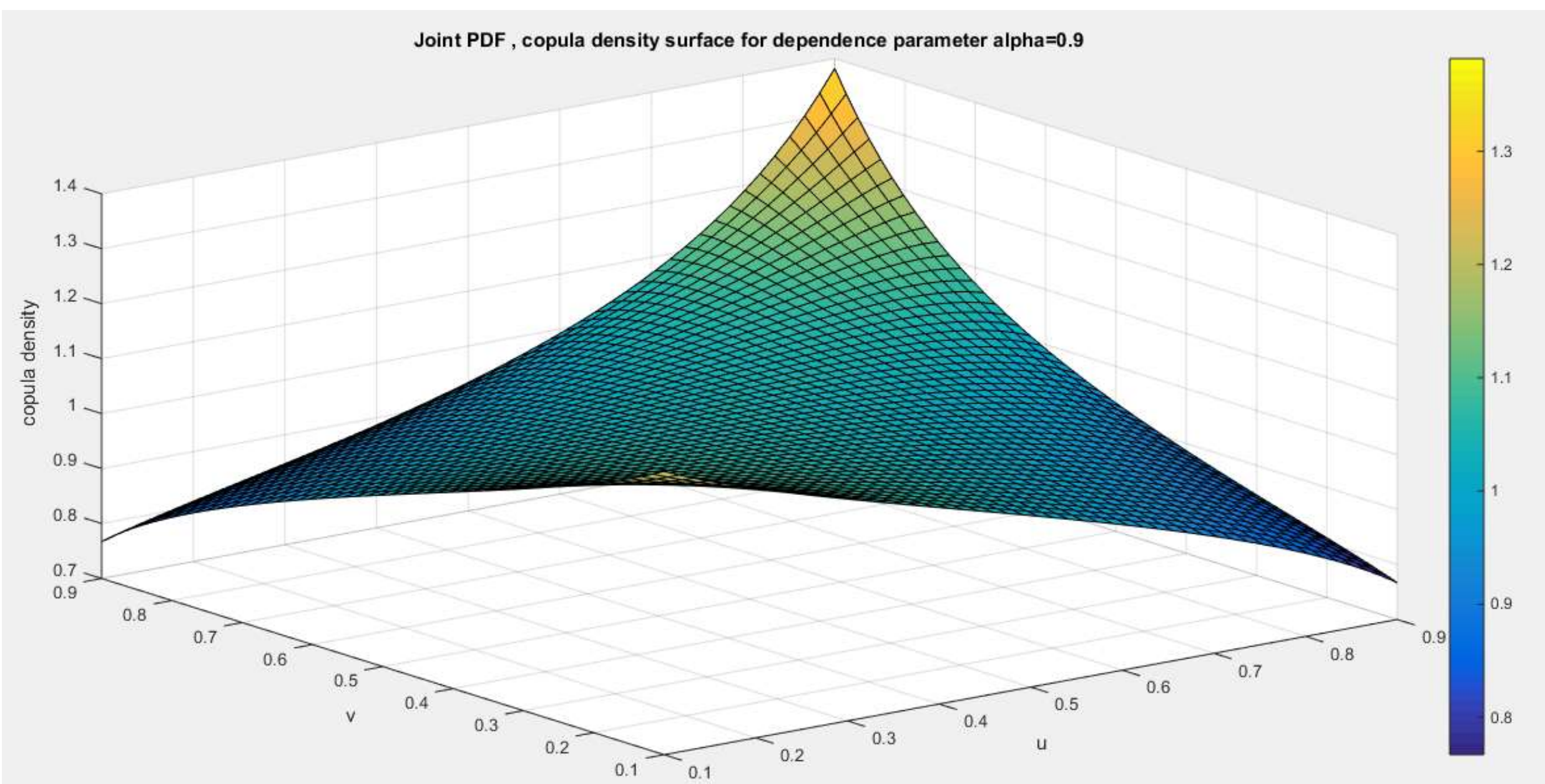
Joint PDF , copula density surface for dependence parameter alpha=0.9
copula density
v
u

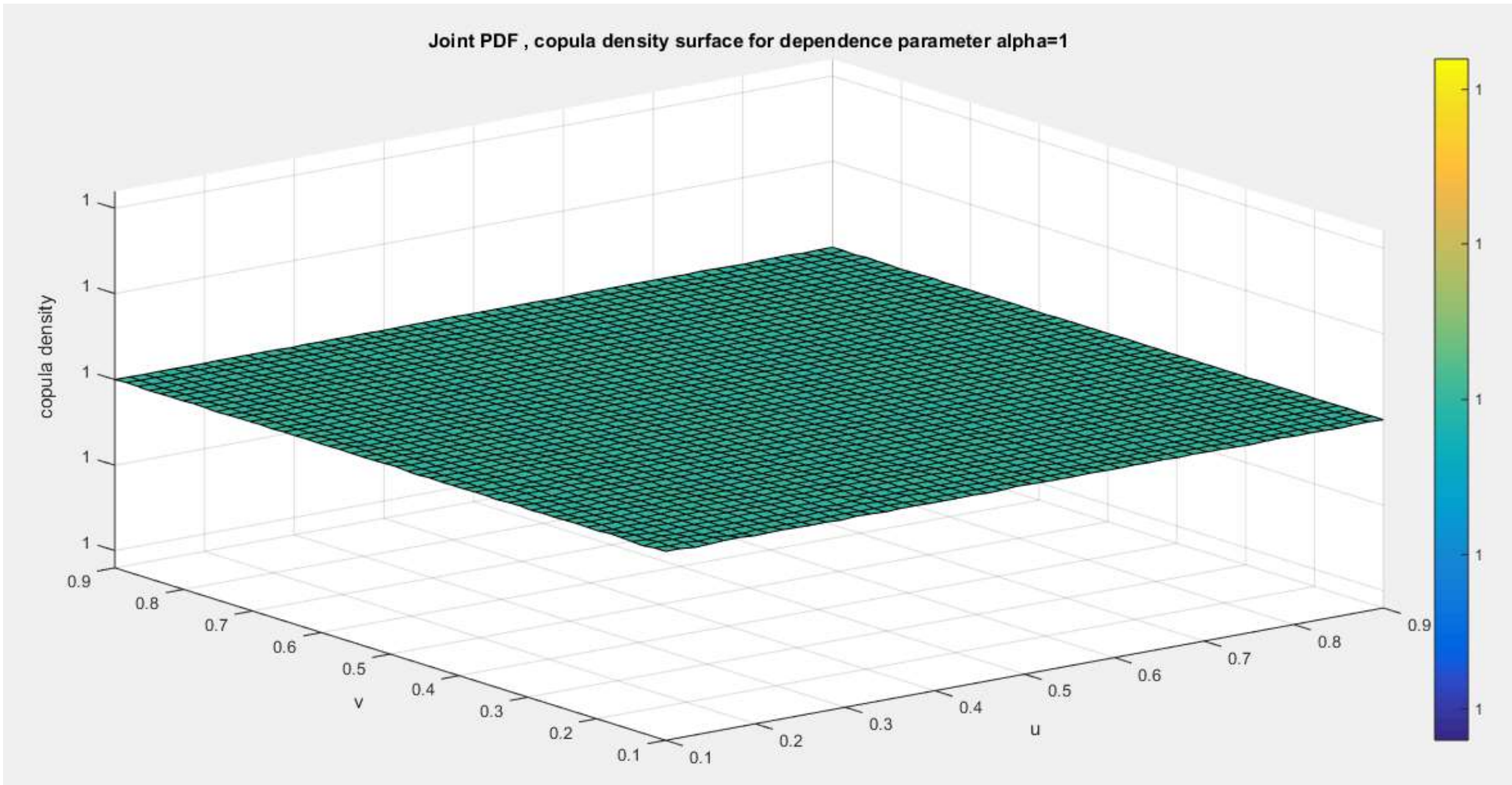

This is the product copula that can be attained at dependency parameter 1

# Second copula

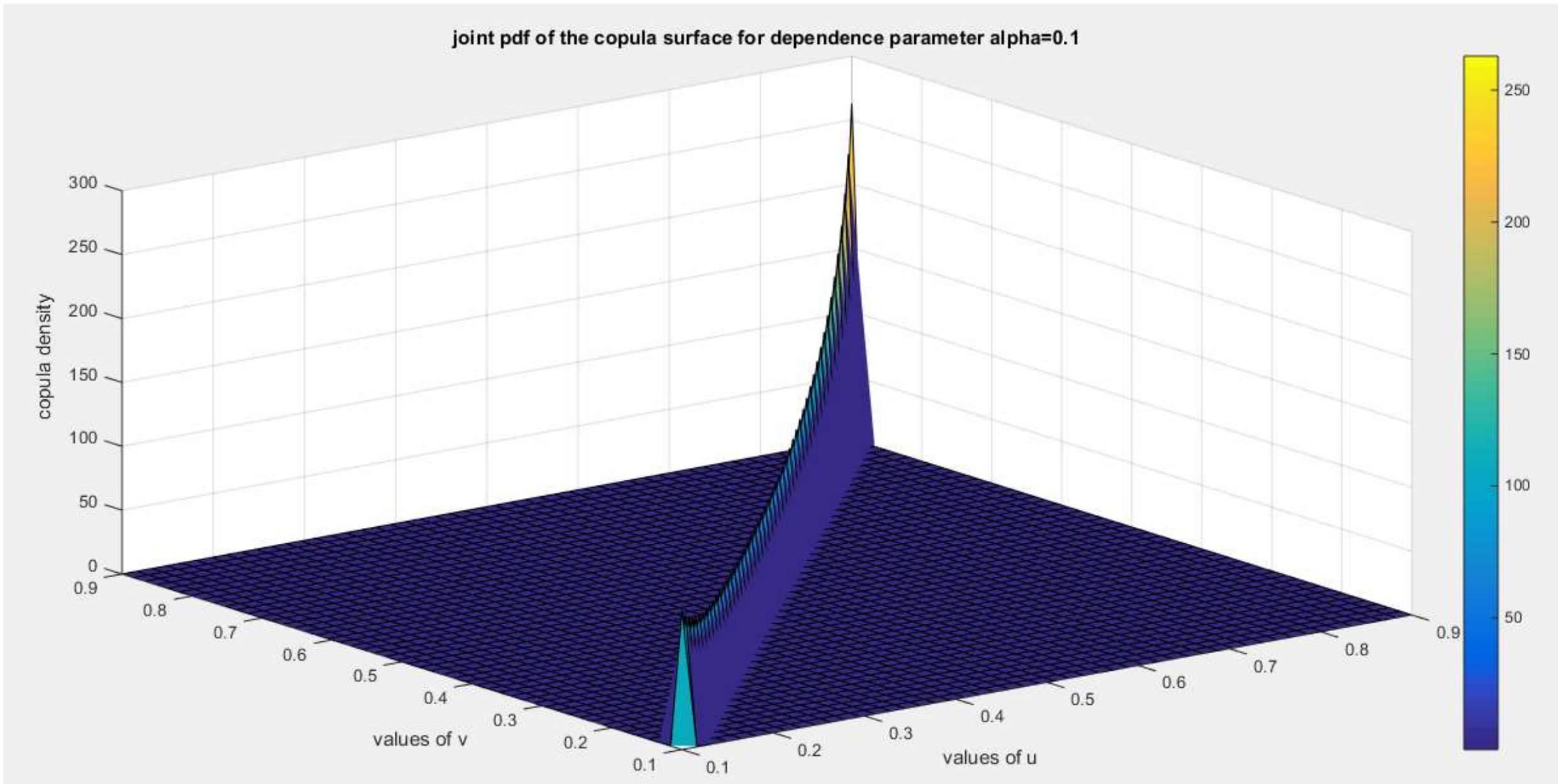

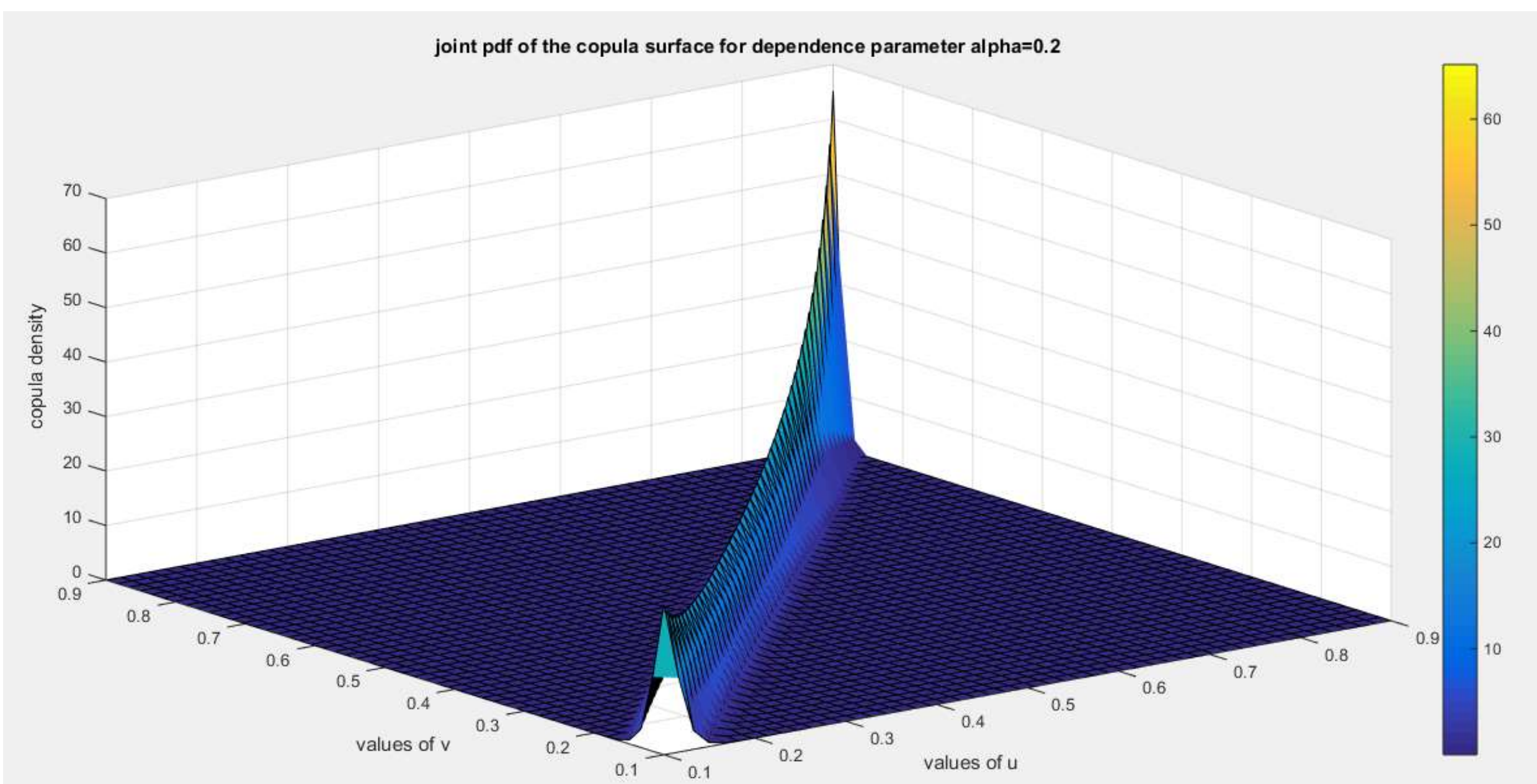

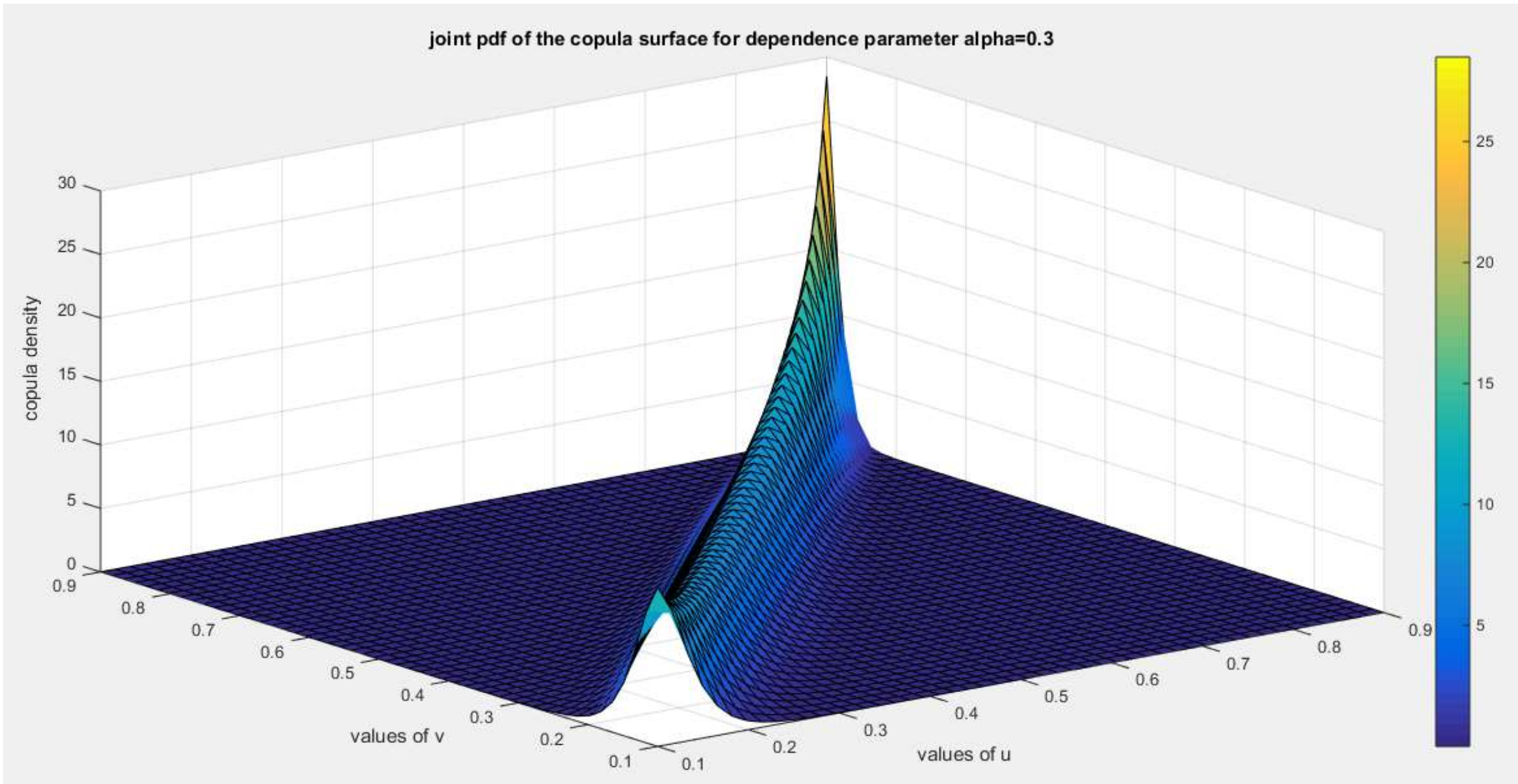
joint pdf of the copula surface for dependence parameter alpha=0.3
copula density
values of v
values of u

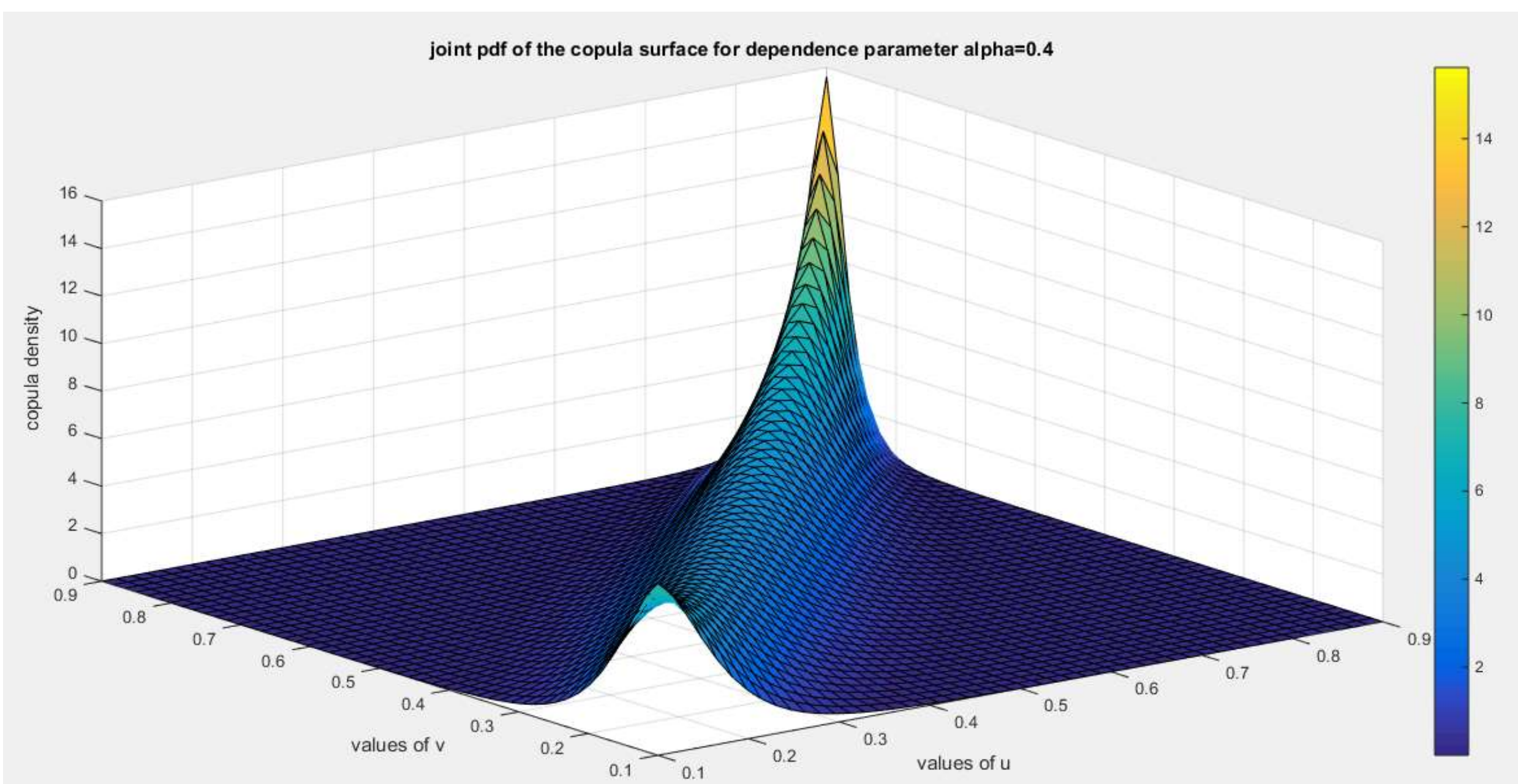
joint pdf of the copula surface for dependence parameter alpha=0.4
copula density
values of v
values of u

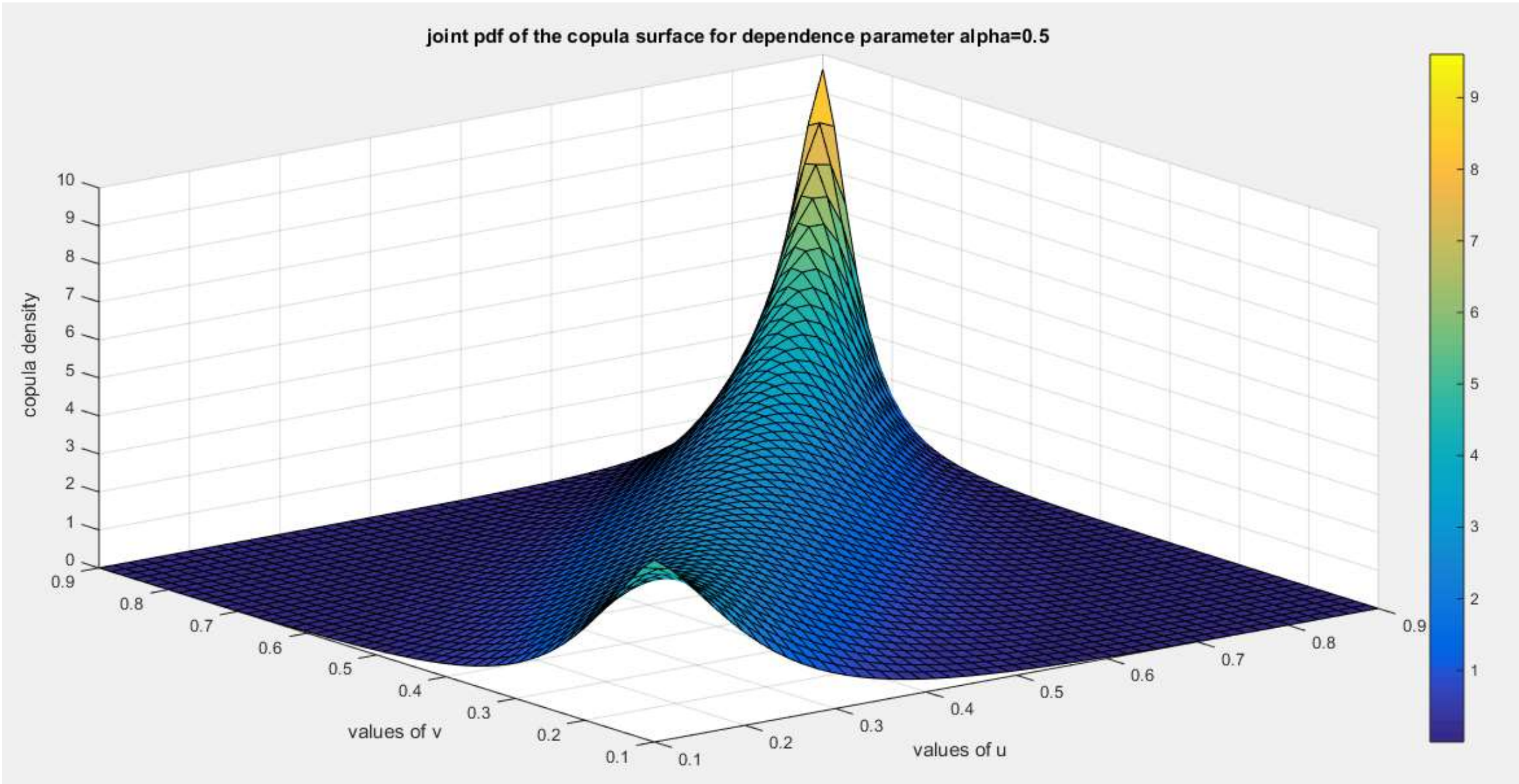
joint pdf of the copula surface for dependence parameter alpha=0.5
copula density
values of v
values of u

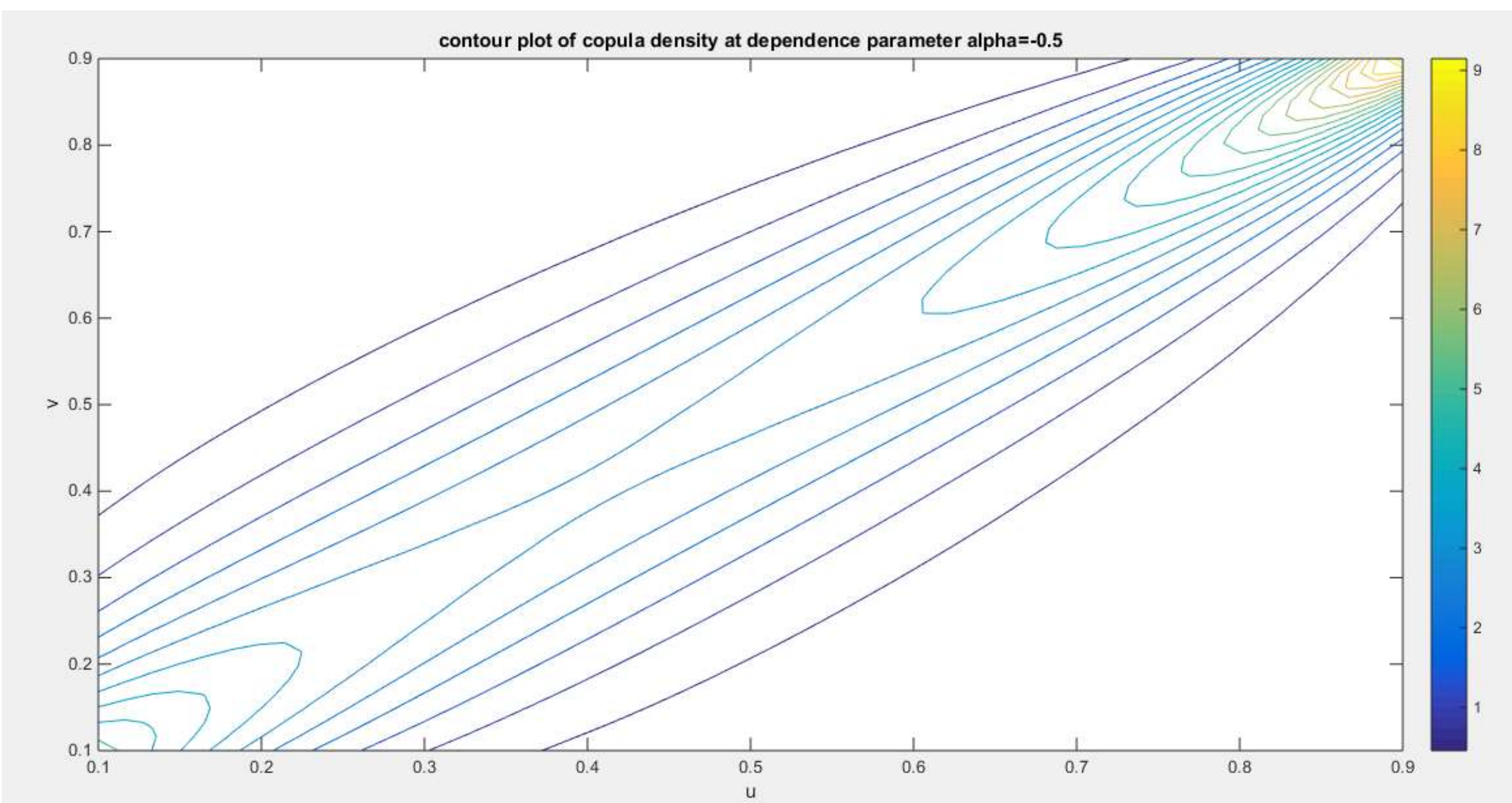
contour plot of copula density at dependence parameter alpha=-0.5
v
u

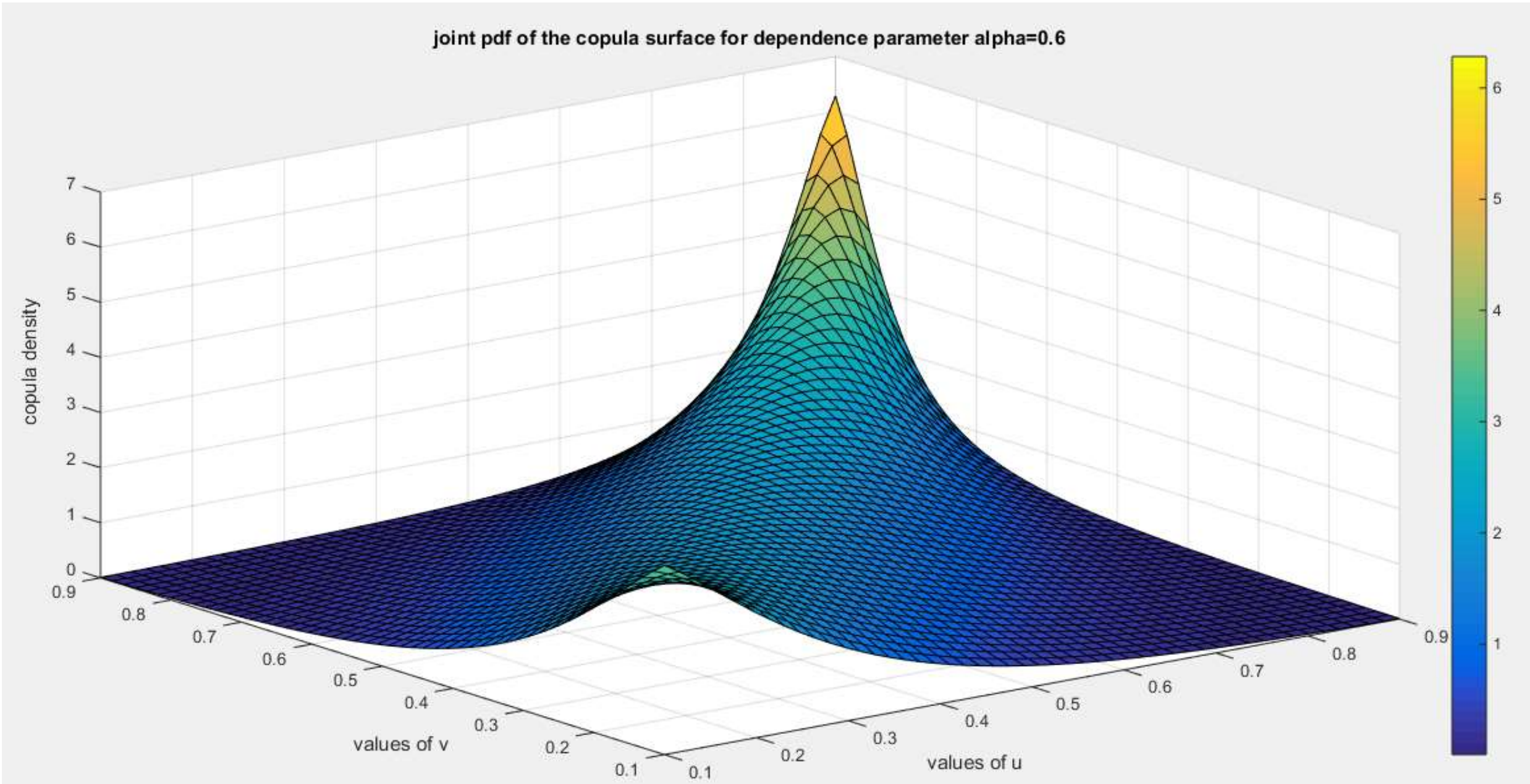
joint pdf of the copula surface for dependence parameter alpha=0.6
copula density
values of v
values of u

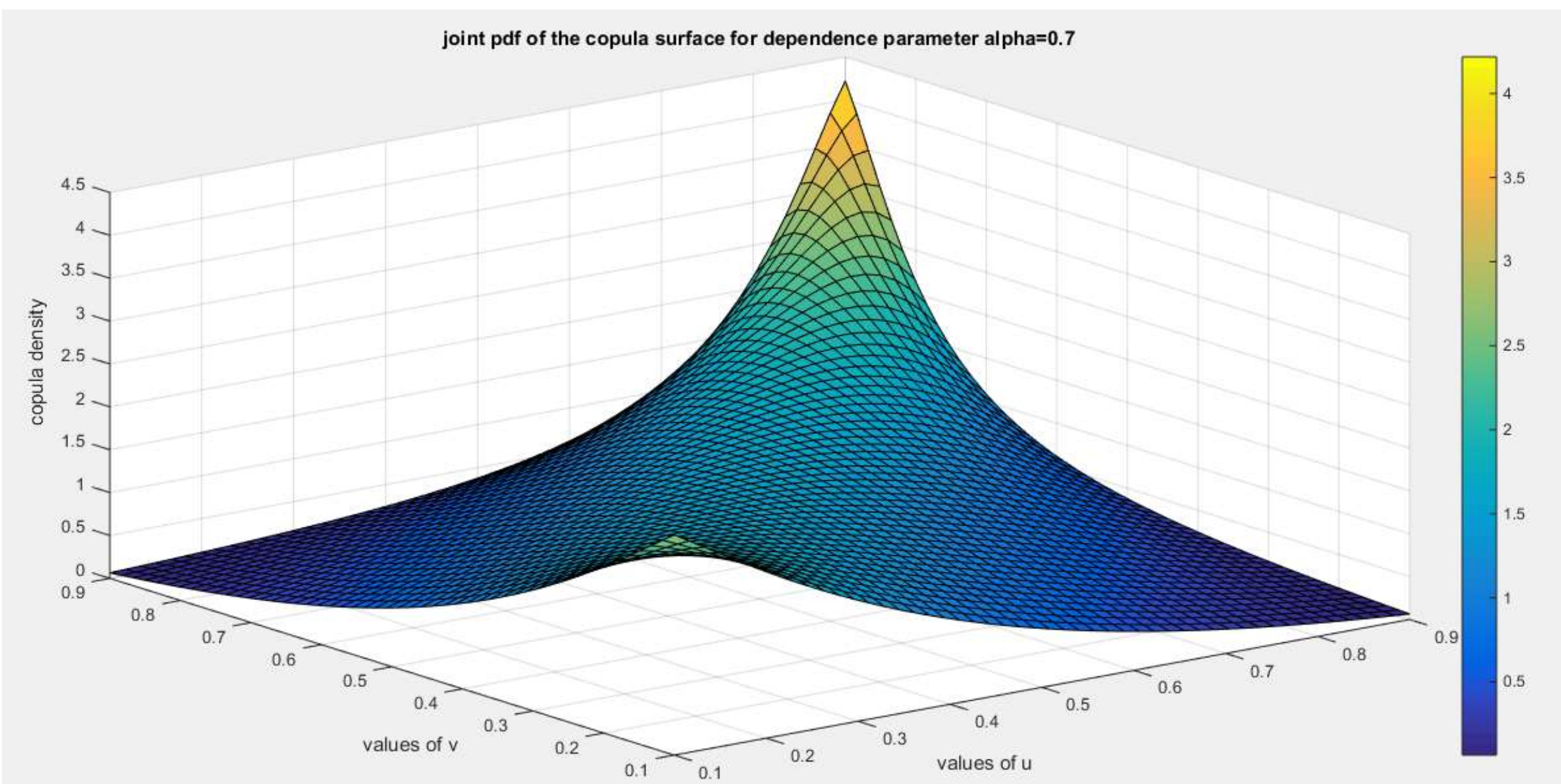
joint pdf of the copula surface for dependence parameter alpha=0.7
copula density
values of v
values of u

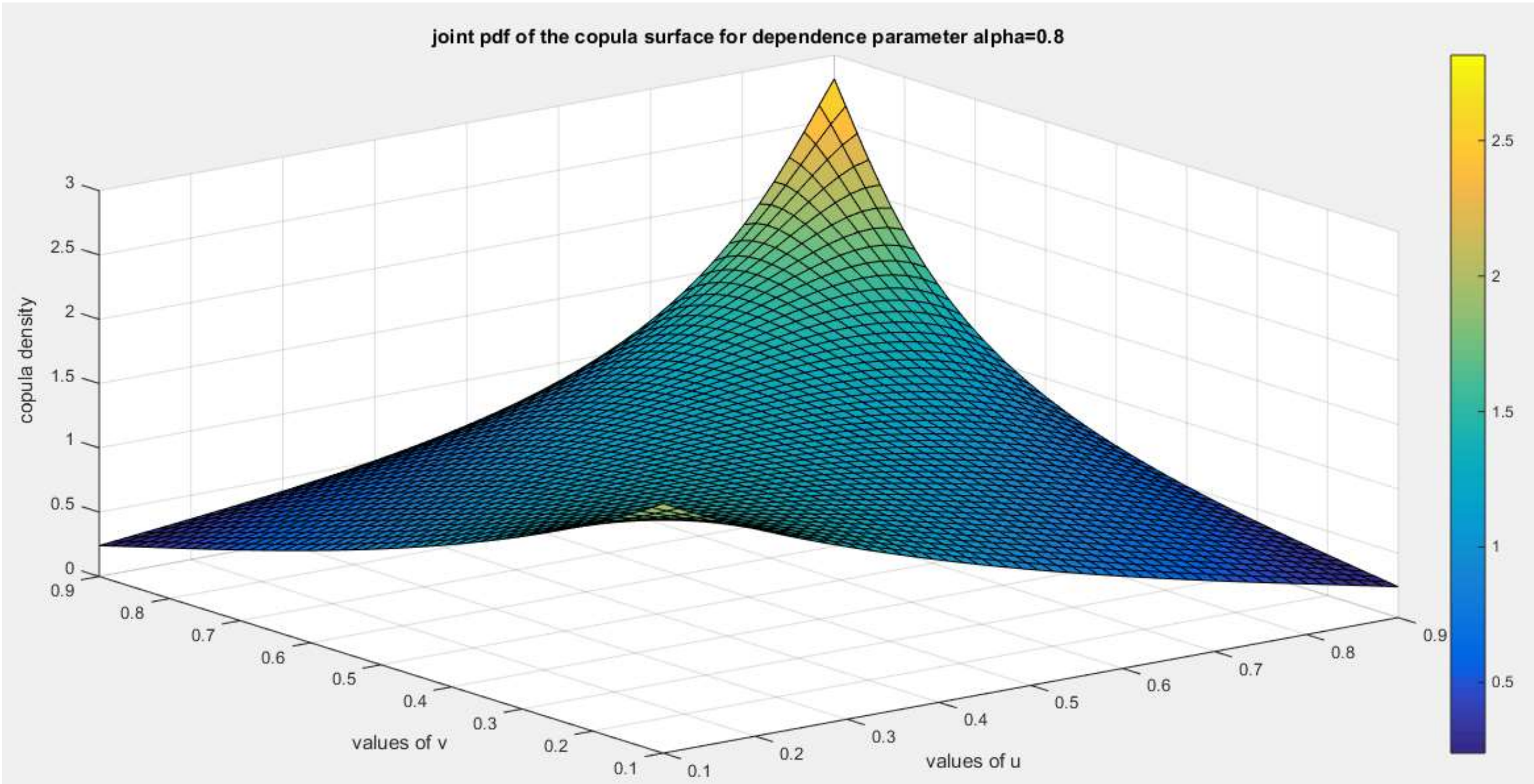
joint pdf of the copula surface for dependence parameter alpha=0.8
copula density
values of v
values of u

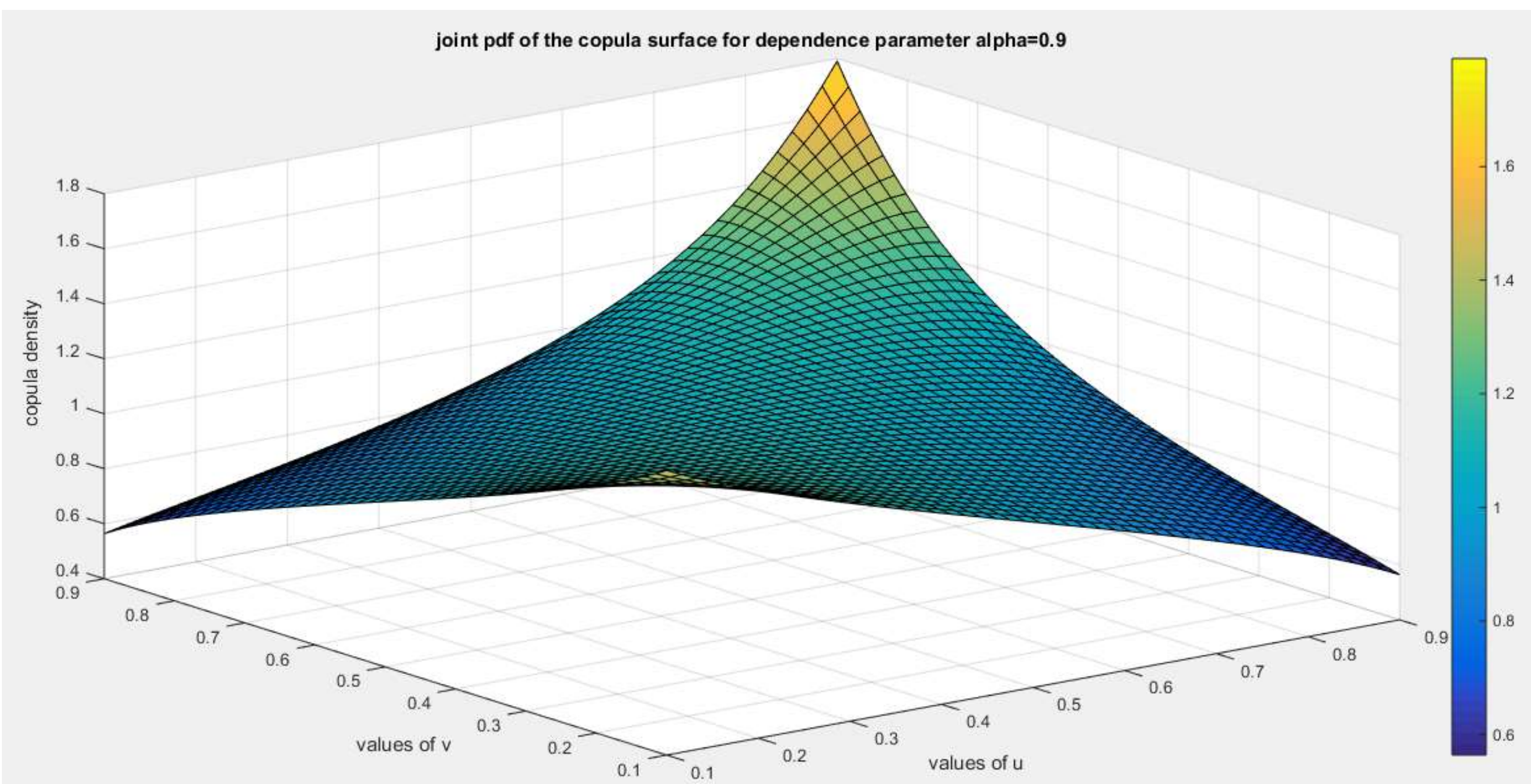
joint pdf of the copula surface for dependence parameter alpha=0.9
copula density
values of v
values of u

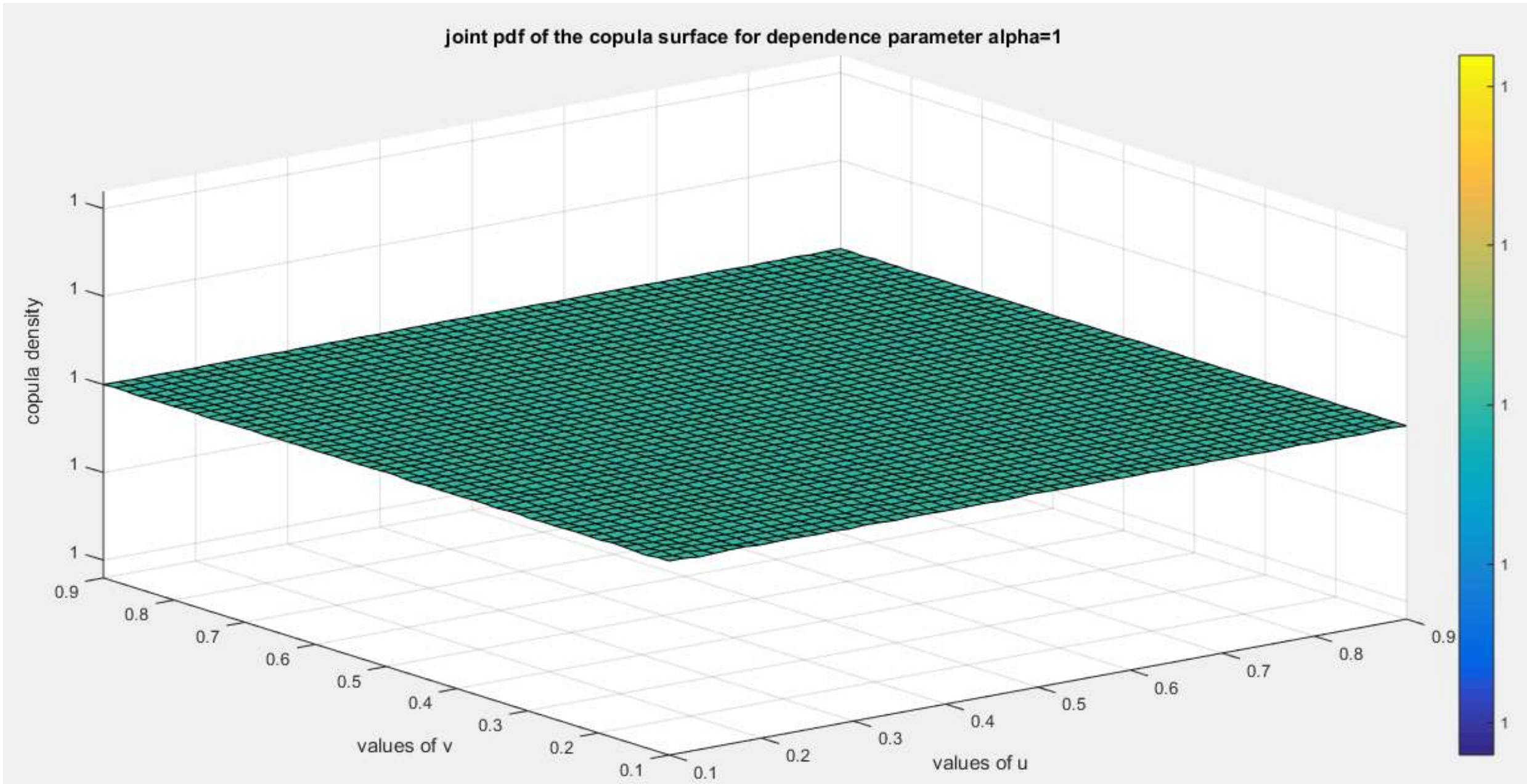

This is the product copula that can be attained at dependency parameter 1 or -1